%% file: bj-arxiv.tex
\documentclass[12pt]{article}
\usepackage{amsfonts}
\usepackage{amsmath}
\usepackage{amssymb}
\usepackage{amsthm}
\usepackage{setspace}
\usepackage{hyperref}
\usepackage[title]{appendix}
\hypersetup{
    colorlinks=true,
    linkcolor=blue,
    filecolor=magenta,      
    urlcolor=blue,
    citecolor= OliveGreen
}
\usepackage[margin = 1.5in, top = 10pt]{geometry}
\usepackage{times}
\usepackage{graphicx}
\usepackage{subcaption}
\usepackage{enumitem}
\providecommand{\keywords}[1]{\textbf{Keywords:} #1}
\usepackage{bbm}
\usepackage{natbib}
\usepackage[dvipsnames]{xcolor}
\usepackage[ruled,vlined]{algorithm2e}
\usepackage{authblk}
\newtheorem{theorem}{Theorem}
\newtheorem{lemma}{Lemma}
\newtheorem{definition}{Definition}
\newtheorem{corollary}{Corollary}
\newtheorem{assumption}{Assumption}
\newtheorem{remark}{Remark}

\newtheorem{proposition}{Proposition}
\theoremstyle{definition}
\input{roy-macro}
\numberwithin{equation}{section}
\numberwithin{theorem}{section}
\numberwithin{definition}{section}
\mathchardef\mhyphen="2D
\newcount\Comments  % 0 suppresses notes to selves in text
\usepackage{color}
\newcommand{\kibitz}[2]{\ifnum\Comments=1\textcolor{#1}{#2}\fi}
\usepackage{color}
\definecolor{darkgreen}{rgb}{0,0.5,0}
\definecolor{purple}{rgb}{1,0,1}

\newcommand{\sapta}[1]{\kibitz{purple}     {[SR: #1]}}

\begin{document}

%\jname{Biometrika}
%% The year, volume, and number are determined on publication
%\jyear{2017}
%\jvol{103}
%\jnum{1}
%% The \doi{...} and \accessdate commands are used by the production team
%\doi{10.1093/biomet/asm023}
%\accessdate{Advance Access publication on 31 July 2016}

%% These dates are usually set by the production team
%\received{2 January 2017}
%\revised{1 April 2017}

%% The left and right page headers are defined here:
%\author{S. Roy, A. Tewari \and Z. Zhu}

%% Here are the title, author names and addresses
\title{High-dimensional variable  selection with heterogeneous signals: A precise asymptotic perspective}

\author{Saptarshi Roy $^\star$\quad Ambuj Tewari $^\star$\quad Ziwei Zhu $^\dagger$\\
$^\star$ University of Michigan, Ann Arbor, USA \\ $^\dagger$ Radix Trading, Chicago, USA}

\maketitle

\begin{abstract}
We study the problem of exact support recovery for high-dimensional sparse linear regression under independent Gaussian design when the signals are weak, rare, and possibly heterogeneous. Under a suitable scaling of the sample size and signal sparsity, we fix the minimum signal magnitude at the information-theoretic optimal rate and investigate the  asymptotic selection accuracy of best subset selection (BSS) and marginal screening (MS) procedures. We show that despite the ideal setup, somewhat surprisingly, marginal screening can fail to achieve exact recovery with probability converging to one in the presence of heterogeneous signals, whereas BSS enjoys model consistency whenever the minimum signal strength is above the information-theoretic threshold. To mitigate the computational intractability of BSS, we also propose an efficient two-stage algorithmic framework called ETS (Estimate Then Screen) comprised of an estimation step and gradient coordinate screening step, and under the same scaling assumption on sample size and sparsity, we show that ETS achieves model consistency under the same information-theoretic optimal requirement on the minimum signal strength as BSS. Finally, we present a simulation study comparing ETS with LASSO and marginal screening. The numerical results agree with our asymptotic theory even for realistic values of the sample size, dimension and sparsity. 
\end{abstract}

\keywords{ Heterogeneous Signals; 
 High-Dimensional Statistics;
 Iterative Hard Thresholding;
  Marginal Screening; 
  Model Consistency;
 Variable Selection. } 

\section{Introduction}
\label{sec: Introduction}
% \zzw{Consider $n$ independent observations $(x_i, y_i)_{i \in [n]}$ of $(x, y)$ that follows the linear regression model:}
Consider $n$ independent observations $(x_i, y_i)_{i \in [n]}$ of a random pair $(x, y)$ drawn from the following linear regression model:
\begin{equation}
    \begin{split}
        (x, \varepsilon)\sim \cP_x\times \cP_\varepsilon,\\
        y = x^\top \beta + \varepsilon,
    \end{split}
    \label{eq: base_model}
\end{equation}
where $\calP_x$ is the $p$-dimensional isotropic Gaussian distribution $\sfN_p(0, \sfI_p)$, and $\calP_\varepsilon$ is the standard Gaussian distribution on $\R$. In matrix notation, the observations can be represented as 
\[
Y = X\beta + E,
\]
%\zzw{Switch from $z$ to $E$.}
where $Y = (y_1, y_2, \ldots, y_n)^\top, X = (x_1,x_2,\ldots, x_n)^\top$ and $E = (\varepsilon_1, \varepsilon_2, \ldots, \varepsilon_n)^\top$. The vector $\beta$ is unknown but sparse in the sense that $\norm{\beta}_0 := \sum_{j =1}^p \ind(\beta_j \neq 0) = s$, which is much smaller than $p$. Denote by $\calS(v)$ the set of non-zero coordinates of a vector $v\in \R^p$. Lastly, we denote by $\pr_{\beta_0} (\cdot)$ and $\bbE_{\beta_0}(\cdot)$ the probability measure and the expectation with $\beta = \beta_0$ respectively. In this paper, we focus on the variable selection problem, i.e., identifying $\calS(\beta)$. We primarily use the 0-1 loss, i.e.,  $\pr_\beta(\calS(\hat{\beta})\neq \calS(\beta))$, to assess the quality of the selected model $\calS(\hat\beta)$. %Another choice of loss function could be Hamming loss which is define as $\bbE_\beta (\vert S(\hat{\beta}) \setminus \calS(\beta) \vert + \vert \calS(\beta) \setminus S(\hat{\beta})\vert)$.
{The isotropic Gaussian design has been widely used to conduct precise analysis of variable selection procedures \citep{fletcher2009necessary, genovese2012comparison, ndaoud2020optimal, su2017false, kowshik2021fundamental}. Specifically, these works either derive the necessary and sufficient condition for exact model recovery \citep{fletcher2009necessary, aeron2010information, rad2011nearly, jin2011limits, akccakaya2009shannon}, or establish tight asymptotic bounds of model selection error \citep{genovese2012comparison, ji2012ups, su2017false}. 
The isotropic Gaussian design is also  used in compressed sensing to generate a measurement matrix \citep{candes2006robust, candes2006near, donoho2006compressed} so that one can sense the sparse signals with few measurements of the high-dimensional signal. \cite{scarlett2016limits, wang2010} considered the variable selection problem and studied the information-theoretic limit of support recovery under non-Gaussian setup. However, they assumed that the entries of $X$ are independent and identically distributed. }
 
% \sapta{
% Apart from independent random design, even simpler models like mean models have gained attention in this area. For example, \cite{Gao2020Fundamental} studied the exact variable selection problem under generalized Gaussian mean model under the same asymptotic regime as \cite{genovese2012comparison}. They showed that simple thresholding procedures can achieve minimax optimal performance when the errors are independent and identically distributed with \emph{log-concave} density and obtained sharp phase transitions in this regime. In addition, a very similar problem like ``sensing'' the sparse signal $\beta$ through fewer linear measurements than the number of features has gained attention in modern applications . In literature, this is known as \textit{compressed sensing}. It is often the case that measurement matrix $X$ is generated from multivariate isotropic Gaussian distribution. \cite{donoho2009counting} studies this problem under noiseless case and showed that $n \asymp 2s \log(p/s)$ dictates the minimum scaling, at which basis pursuit succeeds with high probability. On the other hand, for non-zero noise variance, \cite{wainwright2006sharp} showed that necessary and sufficient condition for the success of LASSO in this setting is $n \asymp 2s \log(n-s) + s + 1$.}

Recently there has been growing interest in the variable selection problem in the presence of \textit{weak} and \textit{rare} signal regimes \citep{genovese2012comparison, ji2012ups} {\color{black} where the active signals are highly sparse with very low magnitude of the order $O(\sqrt{(\log p)/n})$}, which is known to be the information-theoretic optimal rate necessary to achieve model consistency. 
% have studied variable selection problem  under Hamming loss which is defined as $ \bbE_\beta (\vert \calS(\hat{\beta}) \setminus \calS(\beta) \vert + \vert \calS(\beta) \setminus \calS(\hat{\beta})\vert)$. 
This regime is ubiquitous in modern data analytics such as those in Genome-Wide Association Study (GWAS). There the 
genes that exhibit detectable association with the trait of interest can be extremely few with weak effects \citep{wellcome2007genome, marttinen2013genome}. 
Moreover, the number of subjects $n$ typically ranges in thousands, while the number of features $p$ can range from tens of thousands to hundreds of thousands. Such a high dimension further adds to the difficulty of identifying the weak signals. Weak and rare signals also arise in multi-user detection problems \citep{arias2011global} where one typically uses linear model of the form \eqref{eq: base_model}. There the $j$th column of $X$, denoted by $X_j$, is the channel impulse response for user $j$. The signal received from user $j$ is $\beta_j X_j$. Thus $\beta_j =0$ means that $j$th user is not sending any signal. It is a common practice to model the mixing matrix $X$ as random with i.i.d. entries. Under the presence of strong noise, one might be interested in knowing whether information is being transmitted or not. Typically, in some applications it is reasonable to assume that a very few numbers of users are sending signals. Also due to strong noise environment the signals become quite weak, making them harder to detect. Therefore, from an application point of view, understanding variable selection in \textit{weak} and \textit{rare} signal regimes is crucial. {\color{black}Despite its importance, typically most} of the popular methods such as LASSO \citep{tibshirani1996regression}, SCAD \citep{fan2001variable}, adaptive LASSO \citep{huang2008adaptive} have been extensively analyzed in terms of 0-1 loss when the signals are uniformly strong \citep{zhao2006model, guo2015model, zhang2008sparsity, huang2008adaptive, zheng2014high} in the sense that 
\[
a:=\min_{j \in \calS(\beta)}\abs{\beta_j}\gg \left(\frac{\log p}{n}\right)^{1/2}.
\] 
However, {\color{black}\cite{wainwright2006sharp} established a sharp phase transition for LASSO in terms of exact recovery under a general combination of $(n,p,s)$. Under some regularity conditions on the design matrix, the author shows that $a \gtrsim \{(\log p)/n\}^{1/2}$ is necessary and sufficient for the model consistency of LASSO in terms of 0-1 loss. \cite{zhang2008sparsity} proposed MC+ method based on \textit{minimax concave penalty}, which also achieves model consistency under the optimal rate for $a$.
Although  their theory accommodates weak and rare signal regimes, their analysis is only tight up to multiplicative constants. Other works on weak and rare signal regimes include \cite{genovese2012comparison, ji2012ups}, and \cite{jin2014optimality}.}
% However, very little is known about the variable selection problem in the presence of weak and rare signals.

Besides the weakness and rarity of signals, heterogeneity in the signal {\em strength} is another important feature of modern data applications that has not yet received sufficient attention. 
{\color{black} Roughly, \textit{heterogeneity} in the signal allows the magnitude of the active $\beta_j$'s to differ in an arbitrary fashion, whereas \textit{homogeneity} restricts the magnitude of the active signals to be in the same order}. One limitation of the existing literature on variable selection in the \textit{weak} and \textit{rare} signals regime is that it typically assumes that the true signals are homogeneous \citep{genovese2012comparison, ji2012ups, jin2014optimality}. \cite{ji2012ups} refer to this setup as the \textit{Asymptotically Rare and Weak} (ARW) signal regime. 
Many popular approaches have been shown to enjoy satisfactory variable selection properties under the ARW regime.
For instance, \cite{genovese2012comparison} showed that both LASSO and marginal screening enjoy model consistency in terms of Hamming loss under independent random design.                    \cite{ji2012ups} and  \cite{jin2014optimality} investigated the same problem under sparsely correlated design. They proposed two-stage screen and clean algorithms that also exhibit model consistency in terms of Hamming loss. 
However, their theory heavily relies on homogeneous signals and does not extend to the heterogeneous case that is of interest to us. 
In reality, the ARW setup seldom occurs: the signals almost always have different strengths \citep{li2019weak}.   

% It is not yet fully understood whether one can achieve exact support recovery in the presence of weak, rare {\em and heterogeneous} signals. 

To underscore the contrasting effects of homogeneous and heterogeneous signals in terms of exact model recovery, we study the variable selection property of marginal screening (see Section \ref{sec: Thresholding procedures}). We show that under the presence of strong heterogeneity in the signal, marginal screening fails to recover the exact model with probability converging to 1, whereas under homogeneous signal it can recover the exact model asymptotically \citep{genovese2012comparison}. It turns out that due to heterogeneity, the spurious correlations become large and create impediment to selecting the exact model. In correlated design, a different problem known as \textit{unfaithfulness}  \citep{wasserman2009high, robins2003uniform} prevents marginal screening from achieving model consistency. Specifically, due to ``correlation cancellation'', the marginal correlation between $Y$ and $X_j$ becomes negligible even when $\beta_j$ is large and this ultimately leads to \textit{false negatives}. In this paper, we study independent random design model in which correlation cancellation does not occur. Instead, we identify a different source of problem under the presence of signal heterogeneity that affects the exact variable selection performance of marginal screening. {\color{black} Varying effect of signal heterogeneity in variable selection was also identified in the case of LASSO by \cite{su2017false} and \cite{wang2022price} for i.i.d. Gaussian design under a special asymptotic setting. In particular, \cite{su2017false} studied the tradeoff between power and type-I error of LASSO and showed that strong heterogeneity in the signal helps to reduce the false discovery in the LASSO path. The same effect was also analyzed in more detail in \cite{wang2022price}.
These works use approximate message passing (AMP) theory to obtain the exact asymptotic behavior of LASSO estimator in terms of variable selection and show that it is unable to achieve model consistency under \emph{linear sparsity} regime.}
% \sapta{A similar effect of high signals in variable selection for LASSO was also pointed out in \cite{su2017false}. They showed that even under independent design, high signal-to-noise ration introduces early false discovery in the LASSO path when both $s/p$ and $n/p$ approach a positive constant asymptotically. It turns out that the shrinkage of regression coefficients due to $\ell_1$-penalty introduces \emph{shrinkage noise}, whose cumulative effect along a direction of a null variable may dwarf the effect of the strong signal. Hence, the null variables get picked up in the LASSO path, resulting in early false discovery.}

On the computational side, modern methods like LASSO, SCAD, MC+ were initially motivated as alternatives to Best Subset Selection (BSS). BSS is in general an NP-hard optimization problem and was believed to be practically intractable even for $p$ as small as $30$. Thanks to recent advancements in algorithms and hardware, the optimal solution to the BSS problem can now be computed, sometimes with approximations, for some practical settings. \cite{Jain2014iterative} showed that a wide family of iterative hard thresholding (IHT) algorithms can approximately solve the BSS problem, in the sense that they can achieve similar goodness of fit with the best subset with slight violation of the sparsity constraint. \cite{liu2020between} studied the optimal thresholding operator for such iterative thresholding algorithms, which manages to exploit fewer variables than IHT to achieve the same goodness fit as BSS. 
% relates CoSaMP with BSS: given the true sparsity $s$, CoSaMP can find a model of size slightly higher than $s$ that  over the exact best size-$s_0$ subset within a few iterations.
\cite{bertsimas2016best} viewed the BSS problem through the lens of mixed integer optimization (MIO) and showed that for $n$ in 100s and $p$ in 1000s, the MIO algorithm can obtain a near optimal solution reasonably fast. 
% \zzw{double check}
% \sapta{see bullet 3 in pg 6 of Bertismas 2016}. 
{\color{black}\cite{bertsimas2020sparse} developed a new cutting plane method that solves to provable optimality the Tikhonov-regularized \citep{Tikhonov1943OnTS} BSS problem for $n$ and $p$ in the 100,000s. \cite{xie2020scalable} considered solving the Tikhonov-regularized BSS via mixed integer second order cone formulation and the largest problem instance they considered has $p\sim 10^3$. Most recently, \cite{hazimeh2022sparse} developed a Branch-and-Bound method that solves the $\ell_0/\ell_2$-regularized BSS problem for $p\sim 10^7$.} A recent work \citep{zhu2020polynomial} proposed an iterative splicing method called Adaptive Best Subset Selection (ABESS) to solve the  BSS problem. They also showed that ABESS enjoys both statistical accuracy and polynomial computational complexity when the design matrix satisfies sparse Reisz condition and minimum signal strength is of order $\Omega\{(s \log p \log \log n/n)^{1/2}\}$. 

Given these recent advances in solving BSS, there has been growing acknowledgment that BSS enjoys significant statistical superiority over the aforementioned alternative methods. \cite{bertsimas2016best} and \cite{bertsimas2020sparse} numerically demonstrated higher predictive power and lower false discovery rate (FDR) respectively of the BSS solution compared to LASSO.  \cite{guo2020best} and \cite{zhu2021early} reported that the approximate BSS solutions provided by IHT have much fewer false discoveries than LASSO, SCAD and SIS, especially in the presence of highly correlated design. They also theoretically showed that the model selection behavior of BSS does not explicitly depend on the restricted eigenvalue condition for the design \citep{bickel2009simultaneous, van2009conditions}, a condition which appears unavoidable (assuming a standard computational complexity conjecture) for any polynomial-time method \citep{zhang2014lower}. This suggests that BSS is robust against design collinearity in terms of model selection.  
% propose to approximately solve the BSS problem through an iterative hard thresholding (IHT) algorithm with relaxed sparsity constraint. In its core the algorithm is CoSaMP (Compressive Sampling Matching Pursuit), an iterative two-stage hard thresholding algorithm proposed by \cite{needell2009cosamp}. 
% Using this \cite{guo2020best} discovered that CoSaMP (referred to as IHT therein) can inherit sure screening and demonstrates higher True Positive rate (TPR) than existing surrogate methods. A very recent work by  compares BSS with surrogate methods through a new lens: early solution path, i.e., the features entering very early in the BSS solution path. Their study also shows that CoSaMP achieves substantially lower FDR that that of LASSO and SCAD under great variety of design. %These advancements naturally motivates more rigorous theoretical exploration of BSS.

%\sapta{Asymptotic perspective. tight constants}

In this paper, we mainly focus on the \emph{precise}  asymptotic bound, i.e., the bound with the optimal  constant for the minimum signal strength that allows BSS to achieve model consistency. Under a specific asymptotic setup, we show that BSS achieves asymptotic exact recovery of the true model once the minimum signal strength parameter is above the information-theoretic lower bound, meaning that BSS is optimal in terms of the requirement on the signal strength. {\color{black}In contrast, previous works such as  \cite{aeron2010information, wainwright2009information, rad2011nearly} analyze BSS from a sample complexity point of view: they show that BSS can achieve model consistency under the optimal rate of the sample complexity, and under different asymptotic regimes. 
Later \cite{ndaoud2020optimal} showed the existence of a polynomial-time method that achieves model consistency under the same sufficient condition on $n$ as BSS for i.i.d. Gaussian design. For general Gaussian design, \cite{wainwright2009information} showed a similar result for BSS. But the analyses of all these works are tight only up to multiplicative constants. }

%{\color{red} On the other hand, \cite{rad2011nearly, fletcher2009necessary} showed that $\ell_0$-constrained methods enjoy model consistency under i.i.d. Gaussian design when $a\gtrsim \{(\log p)/n\}^{1/2}$. Later \cite{ndaoud2020optimal} showed the existence of a polynomial-time method that achieves model consistency under the same sufficient condition on $a$ as BSS for i.i.d. Gaussian design. For general Gaussian design, \cite{wainwright2009information} showed a similar result for BSS, and under further conditions on the design matrix $X$, \cite{wainwright2009sharp} showed that LASSO achieves models consistency under the same sufficient condition on $a$.}

%Our main contributions are three-fold \zzw{three-fold}. 
The rest of the paper is organized as follows. Section \ref{sec: Ultra rare and  weak minimum signal regime} introduces the Asymptotically Ultra-Rare and Weak Minimum signal (AURWM) regime that accommodates heterogeneous signal strengths.
Section \ref{sec: Thresholding procedures}
shows that in the presence of strong heterogeneity of the signal strength, marginal screening procedures fail to achieve model consistency under the AURWM regime with probability converging to 1. In Section \ref{sec: Analysis of ML decoder}, we derive the asymptotic minimax 0-1 loss under the AURWM regime and show that BSS is optimal in terms of the requirement on the minimum signal strength. 
% \sapta{The proof of this result leverages powerful tools form \cite{fan2018discoveries} which involves delicate analysis of supremum of Gaussian processes. Specifically, in order to obtain the optimal constant, we study a specific maximum spurious correlation process related to the BSS problem and obtain its asymptotic distribution. This is very crucial for the analysis, as this allows us to obtain the sharp asymptotic constants, which is otherwise not possible with routine non-asymptotic probabilistic tools.}
In Section \ref{sec: A novel algorithm}, we propose a computationally tractable two-stage algorithm that also enjoys model consistency under essentially the same condition as BSS. Finally, in Section \ref{sec: simulations}, we carry out simulation studies and numerically demonstrate the superiority of our method over other competing methods. %\ambuj{revise it to say sth less banal than "we corroborate our theory" when the simulation study section has been finalized}.

%\sapta{Bold NOTATION: / don't make a separate section.}

\textbf{Notation.} Let $\R$ and $\R_+$ denote the set of real numbers and the set of non-negative real numbers respectively. Denote by $\R^p$ the $p$-dimensional Euclidean space and by $\R^{p\times q}$ the space of real matrices of order $p\times q$. For a positive integer $K$, denote by $[K]$ the set $\{1, 2, \ldots, K\}$. 

Regarding vectors and matrices, for a vector $v \in \R^p$, we denote by $\norm{v}_2$ the $\ell_2$-norm of $v$. We use $\sfI_{p}\in \R^{p \times p}$ to denote the $p$-dimensional identity matrix.
For a matrix $A \in \R^{p\times p}$, we denote by $A_j$ and $a_j$ the $j$th column and the transposed $j$th row of $A$ respectively.
%We let $\tr(A) = \sum_{j=1}^p A_{jj}$ denote the trace of $A$.

Throughout the paper, let $O(\cdot)$ (respectively $\Omega(\cdot)$) denote the standard big-O (respectively big-Omega) notation, i.e., we say $a_n = O(b_n)$ (respectively $a_n = \Omega(b_n)$) if there exists a universal constant $C>0$, such that $a_n \leq C b_n$ (respectively $a_n\geq C b_n$) for all $n \in \mathbb{N}$. Sometimes for notational convenience, we write $a_n \lesssim b_n$ in place of $a_n = O(b_n)$ and $a_n \gtrsim b_n$ in place of $a_n = \Omega(b_n)$. We write $a_n \asymp b_n$ if $a_n = O(b_n)$ and $a_n = \Omega(b_n)$. We denote by $\Omega_{\pr}$ the big-Omega in probability: for a sequence of random variables $\{Z_n\}_{n\geq 1}$ and a sequence of constants $\{a_n\}_{n\geq 1}$, $X_n = \Omega_{\pr}(a_n)$ means that for any $\varepsilon_0>0$, there exist $C_{\varepsilon_0} >0$ and $n_{\varepsilon_0}\in \mathbb{N}$, both of which depend on $\varepsilon_0$, such that 
\[
\pr(\abs{X_n/a_n} < C_{\varepsilon_0})\leq \varepsilon_0, \quad \forall n\geq n_{\varepsilon_0}.
\]
We use $\overset{\rm p}{\to}$ and $\overset{\rm d}{\to}$ to denote convergence in probability and distribution respectively. Also we say $X \overset{\rm d}{=} Y$ for two random variables $X,Y$ if their distributions are equal. We denote by $\ind(\cdot)$ the indicator function.

Finally, regarding probabilistic distributions, we use $\sfN(0,1)$ to denote the standard Gaussian distribution. We use $\sfN_p(0, \Sigma)$ to denote the $p$-dimensional Gaussian distribution with mean zero and variance-covariance matrix $\Sigma \in \R^{p \times p}$. We denote by $\Ber(\pi)$ the Bernoulli distribution with success probability $\pi\in [0,1]$.  

\section{Ultra rare and  weak minimum signal regime}
\label{sec: Ultra rare and  weak minimum signal regime}

In this section, we focus on a specific asymptotic setup that allows \emph{heterogeneity} among the sparse signals in high dimension. 
% Recall the following model of our interest: 
% \begin{equation*}
%     Y = X \beta + z, \quad z\sim \sfN_p(0,  \sfI_{p}),
%     \label{eq: base_model_2}
% \end{equation*}
% \zzw{A bit repetitive}with the same setup described in model \eqref{eq: base_model}. We are interested in exact recovery of the sparse vector $\beta$ under rare and weak signal regime. Formally, we define the signal class as follows:
Throughout our paper, we consider the following signal class:
$$
\M_{s}^{a}:=\{ \beta \in \R^p: \Vert \beta \Vert_0 =  s, \min_{j \in \calS(\beta)}\vert\beta_j\vert \geq a \}.
$$
Here $a$ denotes the minimum signal strength of $\beta$. Note that the signal class $\M_s^a$ only imposes a lower bound for the minimum signal strength and thus allows arbitrarily large magnitudes across the true signals. This implicitly accommodates heterogeneity in the signal, which is in sharp contrast with the homogeneous signal setup considered by \cite{genovese2012comparison}. 

Now we are in a position to introduce the \textit{Asymptotically Ultra Rare and  Weak Minimum signal} regime (AURWM), in which we mainly consider the signal class above with 
\begin{equation}
    a = \left( \frac{2r \log p}{n} \right)^{1/2} \quad 
    \text{and}\quad s = O(\log p),
    \label{eq:aurwm}
\end{equation}
where the parameter $r$ controls the magnitude of the minimum signal strength. {\color{black} As we will see in Section \ref{sec: Analysis of ML decoder}, the model consistency of BSS will depend on the value of $r$.} 
 Besides, we set the sample size $n$ as
\[
 n = \floor{p^k}, \quad 0 < k<1. 
\]
The assumption that $s\lesssim \log p$ characterizes the ultra-rarity of the signals, which is common in genetic studies such as GWAS \citep{yang2020prioritizing}. Unless stated otherwise, from now on our statistical analysis follows  the scalings of $n, p, s, a$ in this AURWM regime. We say a support estimator $\widehat{\calS}$ achieves \textit{asymptotic consistent recovery} in the AURWM regime if
\begin{equation}
\lim_{p \to \infty} \sup_{\beta \in \M_s^a} \pr_\beta(\widehat{\calS}\neq \calS(\beta) ) = 0.
    \label{eq: consistent 0-1 loss recovery}
\end{equation}
This paper mainly focuses on the criterion \eqref{eq: consistent 0-1 loss recovery} to measure the quality of exact recovery performance for an estimator $\widehat{\calS}$.

It is also worth mentioning that a relevant but different asymptotic setup is studied by \cite{genovese2012comparison} and \cite{ji2012ups}. There the authors assumed a Bayesian model such that all the signals are independent and identicially distrbuted and that the sparsity $s\sim p^{1-\vartheta}$ for some $\vartheta \in (0,1)$. Under such a setup they obtained asymptotically tight phase transition boundaries with respect to Bayesian Hamming risk, which partitions the $r \mhyphen \vartheta$ plane into three regions: (a) Region of exact recovery, (b) Region of almost recovery, (c) Region of no recovery. We skip the details of these results for brevity. The major differences between their setup and ours are twofold: (1) They essentially assume homogeneous signals; (2) They assume $s$ to grow in a polynomial fashion with respect to $p$. 

\section{Marginal screening under heterogeneous signal}
\label{sec: Thresholding procedures}

Marginal screening (MS) is one of the most widely used variable selection methods in practice. It selects the variables with top absolute marginal correlation with the response. Formally, for any $j \in [p]$, write $\mu_j := X_j^\top Y/n$. Given any possibly data-driven threshold $\tau(X,Y)$,
define the marginal screening estimator as follows:
\begin{equation}
\widehat{\calS}_{\tau}:= \{ j\in [p]: \abs{\mu_j} \geq \tau(X,Y)\}. 
\label{eq: thresholding estimator}
\end{equation}
Note that $\mu_j$ is essentially equivalent to the marginal correlation between $X_j$ and $Y$ because of isotropy of $X$. 
Marginal screening has been applied in various fields for feature selection and dimension reduction, including biomedicine \citep{huang2019marginal, lu2005marginal, leisenring1997marginal}, survival data analysis \citep{hong2018conditional, li2016survival}, economics and econometrics \citep{wang2020asset, huang2014feature}. 

%\ambuj{I see. You do have refs on the uses of MS here but they're from Stats journals. Should we look for some influential refs from more domain science journals?}

%Statistical methods for variable selection based on marginal screening on survival data have been studied by \cite{fan2010high}, who extended sure independence screening to survival outcomes based on the Cox model.
Besides the broad applications, marginal screening has been shown to enjoy some desirable statistical properties. \cite{fan2008sure} established the sure screening property of marginal screening under an ultra-high dimensional setup, which serves as theoretical justification for MS to be used for dimension reduction in many applications. Later, \cite{genovese2012comparison} showed that MS enjoys the minimax optimal rate under Hamming loss with homogeneous signals. Nevertheless, as mentioned in Section \ref{sec: Introduction}, precise asymptotic characterization of the 0-1 loss of MS remains fairly underexplored under high dimension, especially in the presence of heterogeneity in signal strength.
%A very recent work by \cite{Gao2020Fundamental} studies the thresholding procedures in well known signal-pulse noise model under general noise setting. Following the same line we also study the properties of marginal screening in the context of high dimensional regression in AURWM regime. 
\subsection{Failure of MS in the AURWM regime}
\label{sec: MS fails}
In this section, we study the 0-1 risk of the MS estimator. Define $\T := \{\widehat \calS_\tau \mid \tau: \R ^ {n \times p} \times \R ^ n \to \R_+ \}$, which is the class of all possible marginal screening estimators. Perhaps surprisingly, under the AURWM regime, we show that MS fails to achieve exact model recovery in the minimax sense.
\begin{theorem}
\label{thm: Thresholding fails in ARMW}
Under the AURWM regime with $n = \floor{p^k}$ for some $k \in (0,1)$, none of the MS estimators of the form \eqref{eq: thresholding estimator} can achieve asymptotic exact recovery, i.e.,
\[
\lim_{p \to \infty}\inf_{\widehat{\calS}_\tau \in \T} \sup_{\beta \in \M_s^a} \pr_\beta(\widehat{\calS}_{\tau} \neq \calS(\beta)) = 1.
\]
\end{theorem}
To understand the main message of this theorem, it is instructive to compare it with the parallel result in \cite{genovese2012comparison} with homogeneous signal. Specifically, \cite{genovese2012comparison} consider a Bayesian setup where all the signal coefficients are independent and identically distributed Bernoulli random variables (up to a universal constant). Under the AURWM regime, $s = O(\log p)$, which implies that $\vartheta =1$ in Theorem 10 of  \cite{genovese2012comparison}. Then Theorem 10 in \cite{genovese2012comparison} says that when $r > 1$, MS enjoys consistency in terms of Hamming risk and thus 0-1 risk too. In contrast, when we broaden the signal class to $\M_s^a$ that embraces possibly heterogeneous signals, the same model consistency fails to hold anymore for MS as shown in Theorem \ref{thm: Thresholding fails in ARMW}. This comparison clearly reveals the curse of signal heterogeneity on MS. {\color{black} However, the above impossibility result does not contradict Theorem 2 in \cite{fletcher2009necessary}. The result therein states that asymptotically $r>(1 + \norm{\beta}_2^2)$ is sufficient for the model consistency of MS. However, in the AURWM regime, $r$ can be smaller than $1+\norm{\beta}_2^2$, which would violate the previous condition. In fact, the proof of the above theorem essentially relies on constructing a sequence of signal patterns that violates the condition $r>(1 + \norm{\beta}_2^2)$ asymptotically. Thus, in a way, the proof techniques of Theorem \ref{thm: Thresholding fails in ARMW} shows that $r > (1 + \Vert \beta\Vert_2^2)$ is also necessary for MS to achieve model consistency in the AURWM regime. Hence, this establishes the sharpness of Theorem 2 of \cite{fletcher2009necessary}, at least in AURWM regime.}

To see how signal heterogeneity hurts MS, for any $j \in [p]$, write $\mu_j$ as
\begin{equation} 
\mu_j =  (\beta_j/n) \norm{X_j}_2^2 + X_j^\top(\sum_{\ell\neq j}  X_\ell \beta_\ell +  E)/n =: \mu_j^{(1)} + \mu_j ^ {(2)}.
\label{eq: alpha_mr_slector_simplified}
\end{equation}
Here $\mu^{(1)}_j = n ^ {-1}\beta_j \norm{X_j}_2 ^ 2$ represents the marginal contribution from $\beta_j$ to $\mu_j$, and $\mu_j^{(2)}$ represents the random error of $\mu_j$ due to the cross covariance between $X_j$ and the other signals and noise. Suppose there are spiky signals among $\{\beta_\ell\}_{\ell \neq j}$. Though $\E(\mu ^ {(2)}_j) = 0$ regardless of the magnitude of $\beta_j$, the spiky signals may incur large variance of $\mu_j ^ {(2)}$ and overwhelm the magnitude of $\mu ^ {(1)}_j$, which is the essential indicator of the significance of $\beta_j$. Consequently, for weak signals, one cannot tell if $\beta_j$ is a true variable based on only $\mu_j$ in the presence of spiky signals. {\color{black} Hence, there is a chance that the true variables associated with weak signals would lose to a noise variable and ultimately leading to false discovery. To rigorously show these claims, we construct a specific example as mentioned before and we study the asymptotic limits of $\max_{j \notin \cS(\beta)} \mu_j^{(2)}$} and $\mu_{j_0}$, where ${j_0}$ denotes the index of a weak signal. While the asymptotic analysis of $\mu
_{j_0}$ is rather straightforward, we borrow some non-trivial results from \cite{fan2018discoveries} to obtain the asymptotic properties of $\max_{j \notin \cS
(
\beta)} \mu_j^{(2)}$. Details of the proof can be found in Section A of the supplementary material.

{\color{black} In contrast, the AMP line of works on LASSO in \cite{su2017false} and \cite{wang2022price} show that under certain asymptotic regime signal heterogeneity actually helps LASSO in terms of variable selection. Specifically, under i.i.d. Gaussian design and  \emph{linear sparsity} regime (i.e. $s/p \to \alpha$ for some constant $\alpha\in(0,1)$), \cite{wang2022price} show that higher signal heterogeneity delays the inclusion of false variables in the LASSO solution path whereas, under signal homogeneity, false discovery occurs in a much earlier stage in the solution path. The effect is somewhat opposite to what we discussed for MS. The reason perhaps lies in the fact that LASSO tries to select the features that are highly correlated with the \emph{shrinkage noise} (see \cite{su2017false}), whereas, MS tries to select the features that have a higher correlation with the response. In the case of LASSO, higher signal heterogeneity makes the magnitudes of the correlations between features and shrinkage noise more distinguishable compared to a homogeneous signal pattern and thus preventing early false discovery in the first case. In the case of MS, higher signal heterogeneity increases the variance of $\mu_j^{(2)}$, which essentially dwarfs the influence of weak signals and leads to false discovery. However, these two phenomenon are not directly comparable as the asymptotic settings are different for the two cases. In fact, when $s = O(\log p)$, the effect of shrinkage noise is much smaller (see Section 3.2 in \cite{su2017false}) and such phenomenon does not occur for LASSO.}

\section{Best subset selection}
\label{sec: Analysis of ML decoder}
Now we shift our focus to BSS, one of the most classical variable selection approaches. With the oracle knowledge of true sparsity $s$, BSS solves for 
%  \begin{equation}
%  \hat{\beta}_{\rm best}(s_0):= \argmin_{\beta \in \R^p, \norm{\beta}_0\leq s_0} \norm{Y - X\beta}_2^2.
%  \label{eq: BSS_optimization}
%  \end{equation}
% If the true sparsity $s$ is known form oracle knowledge then we set $s_0 = s$ and Equation \eqref{eq: BSS_optimization} produces the least square solution

\[
\hat{\beta}_{\rm best}\in \argmin_{\beta \in \R^p, \norm{\beta}_0= s} n^{-1}\norm{Y - X\beta}_2^2.
\]
% Thus the BSS outputs the support $\widehat{\calS}_{\rm best}$ of $\hat{\beta}_{\rm best}$. 
Define $P_{\D}:= X_{\D}(X_{\D}^\top X_{\D})^{-1} X_{\D}^\top$, which is the orthogonal projection operator onto the column space of $X_{\D}$. The BSS above can be alternatively viewed as solving for
\begin{equation}
\widehat{\calS}_{\rm best} :=  \argmin_{\D \subseteq [p]: \abs{\D}=s} n^{-1}{Y^\top (\sfI_n - P_{\D})Y} = \argmax_{\D \subseteq [p]: \abs{\D}=s} n^{-1}{Y^\top P_{\D} Y}.
\label{eq: ML decoder}
\end{equation}
Using a union bound as in \cite{wainwright2009information} or \cite{guo2020best}, one can show that there exists a universal positive constant $\varphi$ (approximately equal to 0.618) such that whenever $r> 4/(1- \varphi)$, BSS  achieves model consistency, i.e.,
\[
\lim_{p \to \infty}\sup_{\beta \in \M_s^a}\pr_\beta(\widehat{\calS}_{\rm best} \neq \calS(\beta)) =0.
\]
% The next theorem shows that when signal strength $r$ is sufficiently large in AURWM regime then oracle BSS can achieve exact asymptotic full recovery.
% \begin{theorem}
% Under AURWM regime if $r > \frac{4}{1 - \varphi}$ where $\varphi$ is an universal constant (approximately equal to 0.618) then the oracle BSS \eqref{eq: ML decoder} achieves asymptotic full recovery, i.e.,
% \[
% \lim_{p \to \infty}\sup_{\beta \in \M_s^a}\pr_\beta(\widehat{\calS}_{\rm best} \neq \calS(\beta)) =0.
% \]
% \label{thm: suffcient condition for ML }
% \end{theorem}
We emphasize that the requirement on $r$ here is more stringent than needed: we will show that BSS achieves model consistency whenever $r > 1$, which turns out to be the minimal requirement for any approach to obtain exact support recovery.

\subsection{Exact support recovery of BSS }
In the following theorem, we show that 
% BSS achieves model consistency under  That is to say that
$r>1$ is sufficient for BSS to achieve asymptotic exact recovery. Recall that $n = \floor{p^k}$ with $0 < k < 1$.
\begin{theorem}
\label{thm: tight phase transition boundary for ML}
Let $r>1$ and write $\delta = r-1$. Then there exists a universal positive constant $C_0$ such that whenever 
$$ s<  C_0 \min \bigg\{2k,    \frac{\delta^2}{\{(1+ 0.75\delta)^{1/2}+ (1+ 0.5\delta)^{1/2}\}^{2}} \bigg\}\log p,$$ we have 
\[
    \lim_{p\to \infty} \sup_{\beta\in \M_s^a} \pr_\beta (\widehat{\calS}_{\rm best} \neq \calS(\beta))=0.
\]
\end{theorem}

In order for BSS to achieve model consistency, we need to ensure that the maximum spurious correlation, i.e., correlation between the spurious variables and the response, is well controlled so that the best subset does not involve any false discovery. 
One important ingredient of our analysis is the asymptotic distribution of the maximum spurious correlation due to \cite{fan2018discoveries}, based on which we can derive the sharp constant in the minimum signal strength for BSS to be model-consistent. 
It is worth emphasizing that pursuing the exact asymptotic distributions is crucial to obtain constant-sharp results; typically, standard non-asymptotic analysis can only yield optimal rates rather than optimal constants. Detailed proof can be found in Section B.1 of the supplementary material.

Note also that Theorem \ref{thm: tight phase transition boundary for ML} requires $s$ to grow slowly. Given that we have at least $\binom{p-s}{s}$ spurious models and that this number increases with respect to $s$ when $s$ is small, a larger $s$ implies higher maximum spurious correlation due to randomness and thus thinner chance for the best subset to remain the true model. 

\begin{remark}\label{remark: BSS sub-Gaussian}
    The result of Theorem \ref{thm: tight phase transition boundary for ML} can be extended to the sub-Gaussian case. In particular, if the coordinates of $x$ follow i.i.d. distribution with \textit{mean-zero} and \textit{unite variance}, and $\epsilon$ is also distributed as a \textit{mean-zero} sub-Gaussian distribution with \textit{unit variance} and independently from $x$, then BSS is model consistent under a similar condition on sparsity $s$. Details can be found in {\color{black}Section B.2} of the supplementary material. 
\end{remark}

In the next section, we show that $r > 1$ is the weakest possible requirement on the minimum signal strength for {\color{black}BSS} to achieve asymptotic model consistency. 

\subsection{Necessary condition for exact recovery \sapta{\textbf{Add BSS impossibility result}}}
{\color{black}In this section, we show that under the AURWM regime, it is impossible for BSS to achieve model consistency if $r\leq 1$, i.e., $r>1$ is necessary for BSS to exactly recover the true support of $\beta$. The following theorem shows that if $r=1$, the 0-1 loss for BSS is strictly bounded away from 0.}
\begin{theorem}
    \label{thm: failure of BSS}
    Under AURWM regime \eqref{eq:aurwm} with $r =1$ and $n = \floor{p^k}$ for some $k \in (0,1)$, BSS is unable to achieve model consistency, i.e.,
    \[
    \lim_{p \to \infty} \sup_{\beta \in \cM_s^a} \pr (\widehat{\cS}_{\rm best} \neq \cS(\beta)) > \frac{1}{10}.
    \]
\end{theorem}
{\color{black} The above theorem shows that when $r=1$, BSS is unable to achieve model consistency asymptotically. Furthermore, using Theorem 1 of \cite{fletcher2009necessary} in our setting yields that whenever $r<1$, i.e., $r<1-\delta_0$ for some $\delta_0 \in (0,1)$, the 0-1 loss of BSS approaches 1 as $p$ grows to infinity. This shows that if $r\leq 1$, BSS is not model consistent. In the regime $r\le 1$, the main difficulty for BSS arises from the fact that it gets confused between $\cS(\beta)$ and its closest competitors. To be precise, let $j_0$ denote the index of a weak signal, i.e., $\beta_{j_0} = \{(2r \log p)/n\}^{1/2}$ with $r \leq 1$. Due to the weak magnitude of $\beta_{j_0}$, it becomes indistinguishable from 0 and as a result, BSS confuses $\cS(\beta)$ with other candidate models $\{\cD \subset [p] : 
\cS(\beta) \setminus \cD = \{j_0\}, \abs{\cD} = s\}$ of size $s$ that differ only at $j_0$ with non-negligible probability. To prove Theorem \ref{thm: failure of BSS}, we also analyze the asymptotic distribution of an appropriate maximum spurious correlation statistics using results from \cite{fan2018discoveries}.
We point the readers to Section B.3 of the supplementary material for further details of the proof.
These results along with Theorem \ref{thm: tight phase transition boundary for ML}, provide a complete characterization of model consistency of BSS in terms of the magnitude of $r$. 

It is worth mentioning that a more general information-theoretic impossibility result is true for the regime $r<1$. In other words, if $r<1$, then no method can achieve model consistency. Towards this end, we consider the minimax 0-1 loss  $$\inf_{\widehat{\calS}}\sup_{\beta \in \M_s^a} \pr_\beta(\widehat{\calS}\neq \calS(\beta)),$$
}
where the infimum is taken over the class of all possible measurable functions $\widehat{\calS}: (X, Y)\to \{ \D\subseteq [p]: \abs{\D}=s\}$.
The next result establishes a lower bound of the above minimax 0-1 loss.
\begin{proposition}
\label{thm: Information theoretic boundary}
Under the AURWM regime with $ n = \floor{p^k}$ for some $k \in (0,1)$, whenever $r < 1$, there exists a universal positive constant $c$ such that
\[
\lim_{p\to \infty}\inf_{\widehat{\calS}} \sup_{\beta \in \M_a^s} \pr_\beta(\widehat{\calS} \neq \calS(\beta)) \geq c.
\]
\end{proposition}

Proposition \ref{thm: Information theoretic boundary} suggests that $r\geq 1$ is a necessary condition for exact support recovery. Combining this with Theorem \ref{thm: tight phase transition boundary for ML} and Theorem \ref{thm: failure of BSS}, we can see that BSS is \textit{almost} optimal in terms of the requirement on the constant $r$ in minimum signal strength to achieve model consistency.
The proof of Proposition \ref{thm: Information theoretic boundary} leverages Theorem 1 in \cite{wang2010} and detailed proof can be found in {\color{black}Section B.4 of the supplementary material}.

{\color{black}
\vspace{5pt}
\hspace{-12pt}
\textbf{Comparison with previous literature:} As pointed out before in Section \ref{sec: Introduction}, there is sharp contrast between the above results and the results in the previous works like \cite{wainwright2009information, rad2011nearly, aeron2010information}, where the authors study the necessary and sufficient conditions for model consistency of BSS in terms of sample complexity under different asymptotic regimes. For example, under \textit{strong-noise} regime, \cite{aeron2010information} showed that the necessary and sufficient conditions for model consistency in terms of 0-1 loss are given by $n = \Omega(s \log(p/s))$ and $a^2 = \Omega(\log(p-s))$, and BSS is optimal in the sense that it achieves exact recovery under these conditions. For the \emph{fixed noise-variance} regime, the 
 results are different. Firstly, \cite{wang2010} showed that the following condition is necessary for any method to achieve exact recovery:
 \begin{equation}
 n = \Omega \left( \frac{s \log (p/s)}{\log (1 + s a^2) } \vee \frac{\log(p-s)}{\log (1 + a^2)}\right),
 \label{eq: necessary conds.}
 \end{equation}
 where $u \vee v := \max\{u,v\}$.
Under the restriction that $a = O(1)$ and $a = \Omega(1/\sqrt{s})$, which represents \textit{strong-signal} regime, \cite{rad2011nearly} showed that BSS achieves model consistency under the necessary condition \eqref{eq: necessary conds.}. In the general case, that is with no assumption on the joint behavior of $(n,p,s, a)$, \cite{wainwright2009information} established that $n = \Omega (\max\{s \log(p/s), a^{-2}\log(p-s)\})$ is a sufficient condition for model consistency of BSS.
% :
% \[
% n = \Omega \left( s \log(p/s) \vee \frac{\log (p-s)}{a^2}\right).
% \]
One can check that the previous condition match with condition \eqref{eq: necessary conds.} under the weak signal regime $a = O(1/\sqrt{s})$. This indicates that BSS is also optimal in this regime in terms of sample complexity. It is interesting to note that the AURWM regime \eqref{eq:aurwm} also falls under this regime as $a = O(\sqrt{(\log p)/n}) \ll 1/\sqrt{s}$, and $n \asymp p^k$. However, all of these results are tight only up to multiplicative constants and do not study the precise dependence on $a$ in terms of sharp requirement on the constant $r$. 
}

\section{Achieving  {\color{black} information-theoretic optimality} with computational efficiency}
\label{sec: A novel algorithm}
\sapta{1. Clarify the meaning of info-theoretic bound. 2. As we are doing asymptotics there is not story of tradeoff that is usually present in the literature with finite sample bound/ sample complexity 3. Mention Ndoud's work which shows BSS and their method achieve similar sample complexity.}
In spite of the optimality of BSS in terms of model selection, its NP-hardness seriously restricts its practical applicability. To address the computational issue, we propose a two-stage algorithm {\color{black}framework} called ETS (Estimate then Screen) that combines an estimation step  with a follow-up coordinate screening step. {\color{black}Under this framework, one has the flexibility to use any sensible algorithm in the first stage that outputs an estimate with a good estimation guarantee for $\beta$. For example, one choice could be the well-known \textit{iterative hard thresholding} (IHT) algorithm \citep{blumensath2009iterative} which is a computational surrogate for BSS and enjoys a desirable estimation guarantee \citep{Jain2014iterative}. Other choices may include algorithms like \textit{pathwise calibrated sparse shooting algorithm} (PICASSO) or \textit{prox-gradient homotopy} (PGH) method that are known to produce good approximate solutions for LASSO (see \cite{Zhao2018Pathwise, xiao2013proximal}) in high-dimensional setup. We show that in the AURWM regime, ETS enjoys the same selection optimality as BSS in terms of the requirement on the minimum signal strength, i.e., ETS asymptotically achieves model consistency whenever $r$ is greater than the information-theoretic threshold 1, which is also the the optimal requirement for BSS to achieve exact recovery.
} 

{\color{black} The above framework is similar to the methodology introduced in \cite{ndaoud2020optimal}. In that paper, the authors used the square-root SLOPE estimator \citep{bogdan2015slope} for the estimation step, and under i.i.d. Gaussian design they showed that their algorithm achieves model consistency under the same \textit{sample complexity} as BSS. However, they do not study optimal dependence on $r$, which is the main focus of our paper. }
%Moreover, ETS uses algorithms such as IHT, CCD, or PGH that are actually used in practice to get approximate solutions for BSS or LASSO and our theory takes both the statistical error and optimization error into account, which is not the case in \cite{ndaoud2020optimal}. }\sapta{Move this in the intro of this section}
% This makes ETS more attractive for practical applications compare to BSS.  

%Performance of ETS crucially relies on the $\ell_2$-estimation guarantee of IHT \citep{Jain2014iterative}. To elaborate, we show that ``good'' $\ell_2$-estimation error of the first stage IHT estimator helps the screening step to correctly identify the true signals (see section \ref{sec: IHT based ETS algorithm}). 

% Recently statistical performance of IHT algorithm has come under spotlight. \cite{Jain2014iterative} has established  $\ell_2$-estimation guarantees for IHT estimates under certain regularity conditions of the loss function in high-dimension regression. Also \cite{guo2020best} has shown that IHT can achieve high TPR (True positive rate) in model selection under general fixed design. In this paper, we show that a further screening step on top of IHT can asymptotically recover the true support. 
%But many questions regarding exact support recovery properties of IHT based algorithms remain fairly unexplored. Addressing this line of problem, we show that under AURWM regime ETS provably achieves asymptotic exact recovery under the same optimal condition on minimum signal strength as in Theorem \ref{thm: tight phase transition boundary for ML}.

\subsection{The ETS algorithm}
\label{sec: IHT based ETS algorithm}
{\color{black}In this section, we introduce our ETS algorithm (Algorithm \ref{alg: ETS}) in detail. Given a partition parameter $0 < \gamma < 1$, ETS first splits the full sample $(x_i, Y_i)_{i \in [n]}$ into two subsamples $\D_1, \D_2$ of respective sizes $n_1 = \floor{\gamma n}$ and $n_2 = n - n_1$. Then ETS performs two main steps on these two sub-samples respectively:
\begin{enumerate}
    \item Given an objective function $f_{n_1}(\cdot; \cD_1)$ and a constraint set $\cC \subseteq \bbR^p$, in the estimation step, ETS procures a close approximation to $\beta$ by solving for an approximate solution to the optimization problem 
    \begin{equation}
    \text{minimize}_{\theta \in \cC} f_{n_1}(\theta; \cD_1)
    \label{eq: optimization problem}
    \end{equation}
    via a suitable iterative algorithm $\cA(\cdot,\cdot)$ that takes the objective function $f_{n_1}(\cdot; \cD_1)$ and a set of tuning parameters $\cT_\cA$ as inputs. In particular, in this step, ETS outputs an estimator $\hat{\beta}: = \cA(f_{n_1}(\cdot; \cD_1), \cT_\cA)$ of the true signal vector $\beta$.

    \item In the second step, ETS performs a coordinatewise screening based on $\D_2$ and $\hat{\beta}$ to select the true variables. 
\end{enumerate}}

\begin{algorithm}[h] 
\SetAlgoLined
Input: Data $\cD = \{(x_i, Y_i)\}_{i=1}^n$, objective function $
f_{n_1}(\cdot; \cD_1)$, partition parameter $\gamma$ , threshold parameter $\varsigma$\;
 1. Randomly partition the whole dataset $\cD$ into two disjoint subsets $\cD_1= (X^{(1)}, Y^{(1)})$ and $\cD_2 = (X^{(2)}, Y^{(2)})$ \;
 2.Apply the algorithm $\cA$ to compute an approximate solution $\hat{\beta}$ of the optimization problem \eqref{eq: optimization problem} \;
 3. Construct the statistics $\{\Delta_i\}_{i=1}^p$ and thresholds $\{\kappa_\varsigma(X_i^{(2)})\}_{i=1}^p$ using \eqref{eq:delta}-\eqref{eq: iht_hat_eta_threshold}\;
 4. Finally compute the selector $\hat{\eta}(X,Y)$\;
 Output: The selector $\hat{\eta}(X.Y)$.
 \caption{ETS}
 \label{alg: ETS}
\end{algorithm}

{\color{black}To elaborate more on the method, for $\ell \in \{1, 2\}$,  let $X^{(\ell)} \in \R^{n_\ell \times p}$ and $Y^{(\ell)} \in \R ^ {n_{\ell}}$ denote the design matrix and the response vector of the $\ell$th sub-sample respectively. ETS computes $\hat{\beta}$ based on the first sub-sample $\cD_1 := (X^{(1)},Y^{(1)})$ by finding an approximate solution the optimization problem \eqref{eq: optimization problem} via algorithm $\cA$. In practice, there could be several choices for both the objective function $f_{n_1}(\cdot; \cD_1)$ and the algorithm $\cA$. For example, one of the most common choices is to consider the $\ell_0$-constrained squared-error loss, i.e.,  $f_{n_1}(\theta; \cD_1) = n_1^{-1} \Vert Y^{(1)} - X^{(1)} \theta\Vert_2^2$ with $\cC = \{\theta\in \bbR^p : \norm{\theta}_0 \leq s\}$. In this case, a natural choice for $\cA$ is the IHT algorithm which is basically a projected gradient descent method. Another popular choice for the objective function is the well-known $\ell_1$-regularized LASSO objective function $f_{n_1}(\theta;\lambda,  \cD_1) = n_1^{-1} \Vert Y^{(1)} - X^{(1)} \theta\Vert_2^2 + \lambda \norm{\theta}_1$ with $\cC = \bbR^p$ and one can choose either PICASSO,  PGH or the composite gradient method proposed in \cite{agarwal2012fast} as the algorithm $\cA$. Besides these, another choice could be to solve  the square-root LASSO problem \citep{bogdan2015slope} via proximal-gradient descent algorithm proposed in \cite{li2020fast}.  }

% \begin{equation}
% \label{eq: empiriacal_l1_penalized_loss}
%     f_{n_1, \lambda}(\theta; X^{(1)}, Y^{(1)}):=  n_1^{-1}\Vert Y^{(1)} - X^{(1)} \theta\Vert_2^2 + \lambda \norm{\theta}_1.
% \end{equation}

Next comes the screening step of ETS. For each $i \in [p]$, define
\begin{equation}
    \label{eq:delta}
    \Delta_i := \frac{
    X_i^{(2)\top} \left( Y^{(2)} - \sum_{j\neq i} X_{ j }^{(2)} \hat{\beta}_j  \right) }{ \Vert X_i^{(2)} \Vert_2}
\end{equation}
and
%\vspace{-3pt}
\begin{equation}
    \kappa_{\varsigma}(u) := \frac{a\norm{u}_2}{2} + \frac{\varsigma^2 \log p}{a \norm{u}_2}, \quad \forall u \in \R^{n_2}, 
    \label{eq: iht_hat_eta_threshold}
\end{equation}
where $\varsigma > 0$ is specified later. 
ETS selects the $i$th variable if and only if $\abs{\Delta_i}> \kappa_{\varsigma}(X_i^{(2)})$. 
To see why we can screen variables based on $\{\Delta_i\}_{i \in [p]}$, note that
\begin{equation}
\Delta_i = \beta_i \Vert X_i^{(2)}\Vert_2 + \frac{X_i^{(2)\top}\bigl( \sum_{j \neq i} X_j^{(2)} (\beta_j - \hat{\beta}_j) + E \bigr)}{\Vert X_i^{(2)}\Vert_2}.
\label{eq: alpha_distribution}
\end{equation}
A straightforward argument shows that conditioned on $\D_1$ and $X_i^{(2)}$, $\Delta_i$ is distributed as: 

$$ \Delta_i \,\bigr\vert\, \big(\D_1, X_{i}^{(2)}\big) \overset{\rm d}{=} \beta_i \big\Vert X_i^{(2)}\big\Vert_2 + \bigg\{1 + \sum_{j \neq i}(\beta_j - \hat{\beta}_j)^2 \bigg\}^{1/2} g_i, $$
{\color{black}where $g_i \sim \sfN(0,1)$ and is independent of $X_i^{(2)}$}. If estimation method performs well in the sense that  $\Vert\hat{\beta} - \beta\Vert_2$ is small, then for all $i \in \calS(\beta)$, $\beta_i \Vert X_i^{(2)}\Vert_2$ becomes the dominant term in $\Delta_i$. In contrast, for all $i \notin \calS(\beta)$, $\beta_i \Vert X_i^{(2)}\Vert_2 = 0$ and we thus expect $\Delta_i$ to be small. This suggests the existence of a threshold $t(\cdot)$ on $(\Delta_i)_{i \in [p]}$ that distinguishes the true support $\calS(\beta)$ from the irrelevant variables. We follow \cite{ndaoud2020optimal} to choose the threshold function in \eqref{eq: iht_hat_eta_threshold}, which is shown to be a reasonable choice to identify the true variables.

% Next we see that definition \eqref{eq: iht_hat_eta_entry} involves the following random variable:
% \begin{equation}
% \alpha_i := \frac{X_i^{(2)\top} \left( Y^{(2)} - \sum_{j\neq i} X_{ j }^{(2)} \hat{\beta}^{\rm iht}_j  \right)}{\Vert X_i^{(2)}\Vert_2} = \beta_i \Vert X_i^{(2)}\Vert_2 + \frac{X_i^{(2)\top}\left( \sum_{j \neq i} X_j^{(2)} (\beta_j - \hat{\beta}_j^{\rm iht}) + z \right)}{\Vert X_i^{(2)}\Vert_2}.
% \label{eq: alpha_distribution}
% \end{equation}
% Right hand side of Equation \eqref{eq: alpha_distribution} and simple probability argument shows that conditioned on $\D_1$ and $X_i$ teh quantity $\alpha_i$ has the following distribution:  $$ \alpha_i \mid \D_1, X_{i}^{(2)} \overset{\rm d}{=} \beta_i \Vert X_i^{(2)}\Vert_2 + \{1 + \sum_{j \neq i}(\beta_j - \hat{\beta}_j^{\rm iht})^2 \}^{1/2} g_i, $$ where $g_i \sim \sfN(0,1)$ and independent of $X_i^{(2)}$. 

For each $i \in [p]$, define $\hat{\eta}_i(X,Y) := \ind \{ 
    \abs{\Delta_i}> \kappa_{\varsigma}(X_i^{(2)}) 
\}$
and write $$\hat{\eta}(X,Y) := (\hat{\eta}_1(X,Y), \ldots, \hat{\eta}_p(X,Y)) ^ \top.$$ 
The selector $\hat{\eta}(X,Y)$ is the final estimate of the support $\cS(\beta)$ produced by the ETS algorithm. Algorithm \ref{alg: ETS} shows the detailed steps of the ETS algorithm. 

% \begin{algorithm}[h] 
% \SetAlgoLined
% Input: Data $\cD = \{(x_i, Y_i)\}_{i=1}^n$, objective function $
% f_{n_1}(\cdot; \cD_1)$, partition parameter $\gamma$ , threshold parameter $\varsigma$\;
%  1. Randomly partition the whole dataset $\cD$ into two disjoint subsets $\cD_1= (X^{(1)}, Y^{(1)})$ and $\cD_2 = (X^{(2)}, Y^{(2)})$ \;
%  2.Apply the algorithm $\cA$ to compute an approximate solution $\hat{\beta}$ of the optimization problem \eqref{eq: optimization problem} \;
%  3. Construct the statistics $\{\Delta_i\}_{i=1}^p$ and thresholds $\{\kappa_\varsigma(X_i^{(2)})\}_{i=1}^p$ using \eqref{eq:delta}-\eqref{eq: iht_hat_eta_threshold}\;
%  4. Finally compute the selector $\hat{\eta}(X,Y)$\;
%  Output: The selector $\hat{\eta}(X.Y)$.
%  \caption{ETS}
%  \label{alg: ETS}
% \end{algorithm}

\subsubsection{Model consistency of ETS}
In this section, we establish theoretical guarantees for ETS-IHT. First, we introduce a technical assumption that concerns how fast algorithm $\cA$ can generate a good approximation of the true signal $\beta$.
{\color{black}
\begin{assumption}
    \label{assumption: iteration count}
    %There exists a sequence $\{\alpha_p\}_{p\geq 1} \subseteq [0,\infty)$ converging to $0$, such that the following holds with probability at least $1 - \alpha_p$:
    The following holds with a probability converging to 1 as $p$ diverges to infinity:
    
   \hspace{-.5cm} For any given tolerance level $\epsilon>0$, there exists a suitable set of deterministic tuning parameters $\cT_\cA$ such that the algorithm $\cA$ requires no more than $T(\epsilon, p,  \beta)$ iterations to produce a solution  $\hat{\beta}:= \cA(f_{n_1}(\cdot;\cD_1), \cT_\cA)$ such that $\norm{\hat{\beta} - \beta}_2^2 \leq \epsilon$.
   %inside the $\ell_2$-ball $\{\theta : \norm{\theta - \beta}_2^2 \leq \epsilon\}$.
\end{assumption}
The above assumption essentially tells that the $T(\epsilon, p, \beta)$th iterate of algorithm $\cA$ is already $\epsilon$-close to $\beta$ in squared $\ell_2$-distance with high probability for large enough $p$.} If $\epsilon$ is small, then $\hat{\beta}$ is a good estimate of $\beta$ and we can use it in the screening step to select the variables. Typically, as $\epsilon$ decreases towards 0, the iteration counts $T(\epsilon, p, \beta)$ increases to infinity as higher accuracy generally demands more computation. However, in many examples, as we will see in {\color{black} Section \ref{sec: ETS examples}}, $T(\epsilon, p, \beta)$ depends only poly-logarithmically on $\epsilon^{-1}$, which alleviates the computational cost. 

Next, we define the binary decoder of the true support $\calS(\beta)$ as
$\eta_\beta := (\ind\{\beta_1 \neq 0\}, \ldots, \ind\{ \beta_p \neq 0\}) ^ \top$. 
The following theorem shows that ETS can achieve exact recovery under suitable choices of tuning parameters.

% \begin{algorithm}[h] 
% \SetAlgoLined
% Input: Objective function $f$, sparsity level $\hat s$, step size $h$ \;
%  $\beta^{(0)} =0$\; $t =0$ \;
%  \While{not converged}{
%  $ \beta^{(t+1)} = P_{\hat{s}}^0(\beta^{(t)} - h \nabla_\theta f(\beta^t))$\;
%  $t\leftarrow t+1$
%  }
%  Output: $\hat{\beta}^{\rm iht} = \beta^{(t)}$.
%  \caption{IHT}
%  \label{alg: two-stage IHT}
% \end{algorithm}
% \begin{equation}
% \hat{\eta}_i(X,Y) = \ind \left\{ 
% \abs{ 
% X_i^{(2)\top} \left( Y^{(2)} - \sum_{j\neq i} X_{ j }^{(2)} \hat{\beta}^{\rm iht}_j  \right) }> t(X_i^{(2)}) \norm{X_i^{(2)}}_2
% \right\},
% \label{eq: iht_hat_eta_entry}
% \end{equation}
% for $i \in [p]$. The threshold $t(\cdot)$ in Equation \eqref{eq: iht_hat_eta_entry} is defined as
% \begin{equation}
% t(u) = \kappa_{\varsigma}(u) := \frac{a\norm{u}_2}{2} + \frac{\varsigma^2 \log p}{a \norm{u}_2}, \quad u \in \R^{n_2}.
% \label{eq: iht_hat_eta_threshold}
% \end{equation}

% \sapta{I think we need to keep some of the things that we moved to appendix here. This is needed to clearly state Theorem 4. For example $\epsilon$ needs to be clarified which was previously present here.}
% \zzw{$8\delta_0 \to \delta_0$, set $\epsilon$ to be sth explicit in the proof.}
\begin{theorem}
% Let the assumptions of Theorem \ref{thm: two_stage_IHT_sub_optimality} be satisfied. 
{\color{black}Assume the condition in Assumption \ref{assumption: iteration count} hold and the sample size $n = \floor{p^k}$ for some $k \in (0,1)$.}
Let $r> 1$ and write $\delta = r - 1$. Then, under AURWM regime \eqref{eq:aurwm}, there exist universal positive constants $A_1,A_2$ such that with overall iteration count {\color{black}no more than} $ T(A_1 \delta, p, \beta)$ for algorithm $\cA$, $\gamma \in (0, \delta/(8+ 8\delta))$ and $\varsigma = (1 + A_2\delta)^{1/2} $,  we have that $\lim_{p \to \infty} \sup_{\beta\in \M_{a}^s} \pr_\beta (\hat{\eta}\neq \eta_\beta) =0$.

%\[\sup_{\beta\in \M_{a}^s} \pr_\beta (\hat{\eta}\neq \eta_\beta) \lesssim  p^{-\Omega\left(\frac{\delta^2}{1 + \delta^2}\right)} + \alpha_p\to 0 \quad \text{as $p\to \infty$}.\]

%Let $\hat{\eta}$ be the selector \eqref{eq: iht_selectot_hat_vector}-\eqref{eq: iht_hat_eta_entry} with the threshold $t(\cdot)= t_{\sqrt{1 + \delta_p^2} }(\cdot)$ defined in \eqref{eq: iht_hat_eta_threshold}.  If $r> 1$ then there exists universal positive constants $c_0,  c_2,c_3, A$ such that for ideally chosen partition parameter $\gamma$ and sub-optimality gap $\epsilon$ in Algorithm \ref{alg: two-stage IHT} we have the following: %{\color{red} (Obtain the correct expression for rate. Also finish the proof.)}

% \[
% \lim_{p\to \infty}\sup_{\beta\in \M_{a}^s} \pr_\beta (\hat{\eta}\neq \eta_\beta) \lesssim L_p \; p^{ - \frac{ \left(\frac{r(1- \gamma)}{\sqrt{1 + A \epsilon}} -  \sqrt{1 + A \epsilon} \right)^2}{4r(1- \gamma)}} + \exp(-c_0 n_1) + \exp(- c_4 n_1)+ c_2 p^{-c_3} \overset{p \to \infty}{\longrightarrow}0.
% \]
% Here $L_p$ is a $\text{poly-log(p)}$ term.
\label{thm: IHT-scrrening exact recovery rate}
\end{theorem}

Note that as the signal strength parameter $r$ approaches the information-theoretic boundary, i.e., as $\delta$ approaches 0, ETS {\color{black}may require} more iterations to achieve model consistency as $T(A_1 \delta, p , \beta)$ generally increases as $\delta$ decreases to 0. This is not surprising: intuitively, weaker signals are harder to identify than strong ones. 
%Moreover, the statistical power of ETS methods also suffer as $\delta$ approaches 0 as suggested by the bound in the above display.  

Besides, ETS does not require the knowledge of the true sparsity $s$, but requires the knowledge of $a$ in the second stage for accurate screening. If the true sparsity $s$ is known, then we can enforce ETS to select exactly $s$ features as follows: Let $\abs{\Delta}_{(m)}$ denote the $m$th largest value of  $\{\abs{\Delta_i}\}_{i \in [p]}$. For each $i \in [p]$, define 
\begin{equation}
\hat{\eta}_i (X, Y; s) = \ind\{\abs{\Delta_i} \geq \abs{\Delta}_{(s)}\}.
\label{eq: selector_ETS_tilde}
\end{equation}
 Hence $\hat{\eta}(s): = (\hat{\eta}_1(s), \ldots, \hat{\eta}_p(s))^\top$ selects exactly $s$ features and the knowledge of $a$ is not required in this case. 
The following corollary shows that under the same conditions of Theorem \ref{thm: IHT-scrrening exact recovery rate}, $\hat\eta(s)$ achieves model consistency. 
\begin{corollary}
{\color{black}Assume the condition in Assumption \ref{assumption: iteration count} holds and the sample size $n = \floor{p^k}$ for some $k \in (0,1)$.}
Let $A_1$ be the same universal constant as in Theorem \ref{thm: IHT-scrrening exact recovery rate}, $r> 1$ and write $\delta = r - 1$. Then, under AURWM regime \eqref{eq:aurwm}, with overall iteration count {\color{black}no more than} $ T(A_1 \delta, p, \beta)$ for algorithm $\cA$ and $\gamma \in (0, \delta/(8+ 8\delta))$ ,  we have that
$
     \lim_{p \to \infty} \sup_{\beta \in \M_s^a} \pr_\beta (\hat{\eta}(s) \neq \eta_\beta) = 0.
     %\label{eq: model_consistecy_ETS_tilde}
$
 \label{cor: model_consistecy_ETS_tilde}
\end{corollary}
{\color{black}
\begin{remark}
\label{remark: adaptive threshold}
The algorithm can be made adaptive to $a$ in some certain regime of $r$. In particular, if there exists a known positive constant $\delta_*$ such that $r> 1 + \delta_*$, then a threshold as in \eqref{eq: iht_hat_eta_threshold} can be constructed without the knowledge of $a$ or $r$ so that ETS still enjoys model consistency. In this case, $\delta_*$ can be arbitrarily small and as long as $\delta_*$ is known, an adaptive choice of threshold exists.
\end{remark}
}

Detailed proofs of Theorem \ref{thm: IHT-scrrening exact recovery rate}, Corollary \ref{cor: model_consistecy_ETS_tilde} and Remark \ref{remark: adaptive threshold} can be found in  {\color{black} Section C of the supplementary materials}. Next, we will discuss some concrete examples of ETS methods that enjoys model consistency.
%Besides, the sparsity level $\hat{s}$ and the step-size $h$ are chosen based on certain geometric properties of $f_{n_1}(\theta)$ and details can be found in Section 9.2 of the supplementary materials.
%It is worth mentioning that $\delta_*$ can be arbitrarily close to 0. 

% \begin{remark}
% In this note we remark that one can use the threshold $t_{\sqrt{1 + A\epsilon}}(\cdot)$ instead of the threshold in Theorem \ref{thm: IHT-scrrening exact recovery rate}. It will also give a similar rate of convergence.

% \end{remark}

\subsubsection{Examples of ETS methods}
\label{sec: ETS examples}
{\color{black}In this section, will present a few examples of ETS methods. In particular, we will consider the ETS methods with different choices for the base algorithm $\cA$: (1) the IHT algorithm which solves the $\ell_0$-constrained optimization problem, (2) the PICASSO and PGH algorithm which solves the $\ell_1$-regularized optimization problem. We will show that Assumption \ref{assumption: iteration count} is met in all these cases and we will explicitly derive the dependence of $T(\epsilon, p, \beta)$ on $(\epsilon, p, \beta)$. Hence, this will automatically establish the model consistency of these three variants of the ETS method due to the result in Theorem \ref{thm: IHT-scrrening exact recovery rate}. For clarity, depending on the algorithm used in the estimation step, we will refer to these methods as ETS-IHT, ETS-PICASSO, and ETS-PGH. To be self-contained, we describe the steps of IHT in Algorithm \ref{alg: two-stage IHT}. However, we do not add the description of PICASSO and PGH as those are too involved to add in this paper. Detailed description of PICASSO and PGH can be found in \cite{Zhao2018Pathwise} and \cite{xiao2013proximal} respectively. We remind the readers that
throughout the discussion in this section, we will consider the AURWM regime defined in \eqref{eq:aurwm} with sample size $n = \floor{p^k}$ for some $k \in (0,1)$.  More details and proofs related to the examples can be found in Section C.4 of the supplementary material.} 
% Complete descriptions of PICASSO and PGH can be found in \cite{Zhao2018Pathwise} and \cite{xiao2013proximal} respectively. The following theorem establishes the model consistency of these two variants of ETS.

\textbf{Solving $\ell_0$-constrained problem:}
 As discussed in Section \ref{sec: IHT based ETS algorithm}, In this case, the optimization problem \eqref{eq: optimization problem} takes the form 
        \begin{equation}
        \label{eq: empirical_1_square_error_loss}
        \text{minimize}_{\theta: \norm{\theta}_0 \leq s} n_1^{-1} \norm{Y^{(1)} -  X^{(1)}\theta}_2^2.
        \end{equation}
        We consider the ETS-IHT in this case which uses IHT (Algorithm \ref{alg: two-stage IHT}) to obtain an approximate solution to the above optimization problem. In this case $\cT_\cA = \{\hat{s}, h\}$, where $\hat{s}$ is the sparsity level and $h$ is the gradient step-size. Following the discussion of Section 4 in \cite{Jain2014iterative}, in particular, using Theorem 3 of that paper we have that the final output $\hat{\beta}$ of IHT satisfies $\norm{\hat{\beta} - \beta}_2^2 \leq \epsilon$ with probability converging to 1, when  $\hat{s} = 2592 s, h \leq 8/27$ and $T(\epsilon, p, \beta) = O( \log p + \log((1 +  \norm{\beta}_\infty)/\epsilon))$. Hence, the conditions in Assumption \ref{assumption: iteration count} hold. Moreover, Theorem 3 of \cite{Jain2014iterative} suggests that if $f_{n_1}(\hat{\beta}; \cD_1) - \min_{\theta: \norm{\theta}_0 \leq s} f_{n_1}(\theta; \cD_1) \leq (\epsilon/16)$, then for large values of $p$, we have $\Vert\hat{\beta} - \beta\Vert_2^2 \leq \epsilon$. Hence, it is enough to output an estimator $\hat{\beta}$ which incurs a sub-optimality gap of the order $O(\epsilon)$. 
\begin{algorithm}[h] 
\small
\SetAlgoLined
Input: Objective function $f$, sparsity level $\hat s$, step size $h$ \;
 $\beta^{(0)} =0$\; $t =0$ \;
 \While{not converged}{
 $ \beta^{(t+1)} = P_{\hat{s}}^0(\beta^{(t)} - h \nabla_\theta f(\beta^t))$, \text{where $P^0_{\hat{s}}(v) = \argmin_{z: \norm{z}_0 = \hat{s}} \norm{v - z}_2$}\;
 $t\leftarrow t+1$
 }
 Output: $\hat{\beta} = \beta^{(t)}$.
 \caption{IHT}
 \label{alg: two-stage IHT}
\end{algorithm}

\textbf{Solving $\ell_1$-regularized problem:}
 In this case, the objective function is \[f_{n_1}(\theta; \lambda, \cD_1) = n_1^{-1} \norm{Y^{(1)} - X^{(1)}\theta}_2^2 + \lambda \norm{\theta}_1,\] and $\cC = \bbR^p$. For brevity of discussion, we only consider PICASSO and PGH as the candidate methods for solving the above optimization problem.
 We omit the details of the tuning parameters for these algorithms in this paper, but details of those can be found in \cite{Zhao2018Pathwise} and \cite{xiao2013proximal} respectively.

 \begin{itemize}
     \item (ETS-PICASSO): For ETS-PICASSO, Theorem 3.12 of \cite{Zhao2018Pathwise} yields that with $T(\epsilon, p, \beta) \lesssim (\log p)^3 (\log p + \log \norm{\beta}_\infty) (\log \log p + \log (\epsilon^{-1} \vee C_1))^2 $ and the regularization parameter $\lambda = C_2 \{(\log p)/n_1\}^{1/2}$ for appropriate absolute constants $C_1, C_2>0$, the approximate solution $\hat{\beta}$ has the property $f_{n_1}(\hat{\beta}; \lambda , \cD_1) -  f_{n_1}(\hat{\beta}_L; \lambda , \cD_1)  = O(\epsilon)$ with probability converging to 1, where $$\hat{\beta}_L := \argmin_\theta f_{n_1}(\theta; \lambda , \cD_1) .$$ Then, the strong convexity property of the Gram matrix $X^{(1)\top} X^{(1)}/n_1$, and the good estimation property of $\hat{\beta}_L$ yields that $\Vert \hat{\beta} - \beta \Vert_2^2 \leq \epsilon$.

     \item (ETS-PGH):
      For ETS-PGH, Theorem 3.2 of \cite{xiao2013proximal} yields that with $T(\epsilon, p, \beta ) = O((\log p + \log \norm{\beta}_\infty) \log\log p + \log (\epsilon^{-1}\vee \tilde{C}_1)) $ and $\lambda = \tilde{C}_2 \{(\log p)/n_1\}^{1/2}$ for appropriate absolute constants $\tilde{C}_1,\tilde{C}_2>0$, the approximate solution $\hat{\beta}$ satisfies $f_{n_1}(\hat{\beta}; \lambda , \cD_1) -  f_{n_1}(\hat{\beta}_L; \lambda , \cD_1)  = O(\epsilon)$ with probability converging to 1. Then, again by the strong convexity property of the Gram matrix $X^{(1)\top} X^{(1)}/n_1$, and the good estimation property of $\hat{\beta}_L$, it follows that $\Vert \hat{\beta} - \beta \Vert_2^2 \leq \epsilon$.
     
 \end{itemize}
 It is worth mentioning that the choice of PICASSO or PGH is not special for solving the $\ell_1$-regularized problem. As long as the base algorithm $\cA$ outputs $\hat{\beta}$ which enjoys a sub-optimality gap of the order $O(\epsilon)$ in the functional value, it follows that $\norm{\hat{\beta} - \beta}_2^2 \leq \epsilon$. We formalize this result in the next proposition.

 \begin{proposition}
 \label{prop: lasso optimality-gap}
     Consider the AURWM regime in \eqref{eq:aurwm} and let the sample size $n = \floor{p^k}$ for some $k \in (0,1)$. Then, there exists a positive universal constant $C_3$ such that for all $\epsilon \in(0, C_3)$, the following holds with probability at least $1 - 3p^{-0.5} $ for large enough $p$ and $\lambda  = 8 \{(\log p)/n_1\}^{1/2}$:
     \[
     f_{n_1}(\hat{\beta}; \lambda , \cD_1) -  f_{n_1}(\hat{\beta}_L; \lambda , \cD_1) \leq \epsilon/C_3 \quad \text{implies} \quad \norm{\hat{\beta} - \beta}_2^2 \leq \epsilon.
     \]
 \end{proposition}
The proof of the above result is deferred to {\color{black}Section C.4} of the supplementary material. The above proposition basically shows that in the case of solving the LASSO problem, if $\hat{\beta}$ can be produced efficiently, then an optimality-gap of $C_3^{-1}\epsilon$ in the functional value is enough to guarantee that $\hat{\beta}$ falls inside the $\epsilon^{1/2}$ neighborhood of the true parameter $\beta$, i.e., the conditions in Assumption \ref{assumption: iteration count} hold with probability at least $1 - O(p^{-0.5})$. Hence, this provides us the flexibility to use any sensible algorithm for the $\ell_1$-regularized problem such as the composite gradient method proposed in \cite{agarwal2012fast}.
\sapta{\textbf{Fix the scaling of $\epsilon$ in the proofs}}

\sapta{\textbf{Add supplementary section number}}

{\color{black}
\subsection{Discussion on information-theoretic optimality, statistical accuracy and computational efficiency}
In the previous section, we have shown that ETS achieves the model consistency under the same information-theoretic optimal requirement on $r$ with computational expediency.
 However, the computational efficiency of ETS heavily depends on the magnitude of $r$. Theorem \ref{thm: IHT-scrrening exact recovery rate} suggests that as $r$ approaches the information-theoretic boundary 1, the demand on the number of iterations in the estimation step increases. It could be possible that the computational load of ETS surpasses the computational load of BSS when $r$ is extremely close to 1 for a fixed ambient dimension.
 Hence, even though ETS is able to recover the weak signals \textit{asymptotically} (as $n$ approaches infinity) under the optimal requirement on $r$, it may suffer from high computational costs.  Moreover, as $r$ approaches 1, it turns out that the decaying rate of the error probability worsens. This fact can be verified from the rates obtained in the proof of Theorem \ref{thm: IHT-scrrening exact recovery rate} in the supplementary material and we do not include those in the main theorem for conciseness.
 This suggests that weak signals also hurt the statistical power or accuracy of the ETS methods. Hence, both statistical accuracy and computational efficiency suffer as the signals get weaker. 
 %In comparison, a different line works  \citep{wainwright2009information, wang2010, ndaoud2020optimal} have highlighted the trade-off phenomenon between sample complexity and statistical efficiency highlighted in most of the works like, etc under \textit{finite sample regime}. The main message of these papers is that a smaller sample size worsens the model recovery probability and vice versa. However, in this paper, we consider the asymptotic setting when sample size $n$ approaches infinity, and we are only concerned about the effect of $r$ on the asymptotic model recovery performance of ETS.

}

\section{Numerical experiments}
\label{sec: simulations}
%{\color{red} \textbf{Complete the simulation with marginal screening and ETS}.}
In this section, we first numerically investigate the probability for MS to achieve exact recovery of the true model with growing ambient dimension $p$ under both homogeneous and heterogeneous signal setups. Our results show that while MS exhibits model consistency under the homogeneous signal regime, it completely fails to do so under the heterogeneous signal regime, which is consistent with Theorem \ref{thm: Thresholding fails in ARMW}.  
We then conduct simulation experiments to demonstrate the superiority of ETS methods over competing methods including LASSO and MS as signal strength grows or signal heterogeneity grows. {\color{black} For ETS methods, we only include ETS-IHT and ETS-PICASSO methods. For ETS-PICASSO we used \texttt{picasso} package in R which uses the PICASSO method for solving the LASSO problem}.
To this end, we mention that we do not numerically compare {\em exact} BSS in this section mainly due to computational issues. In most of our simulation setups, {\color{black}we consider $p$ in thousands and exact BSS suffers from high computational costs in such regimes, which is also a limitation of  the commercial solver Gurobi \citep{hastie2020best}. \cite{bertsimas2020sparse, hazimeh2022sparse, xie2020scalable} have made efforts to overcome this computational bottleneck by considering different methods for solving approximate versions of BSS. In particular, they all consider different regularized versions of BSS which is beyond the scope of this paper, and hence we do not include those in the numerical experiments.} Instead, we focus on LASSO and ETS, both of which are two different computational surrogates of the BSS problem.
\sapta{probably need to change this line. should we add approximate BSS in the simulations?}
%\zzw{Why not exact BSS; LASSO as convex relaxation of BSS}
\subsection{Exact recovery performance of MS}
\begin{figure}[h]
\centering
  
    \begin{subfigure}{0.47\linewidth}
\includegraphics[width=\linewidth]{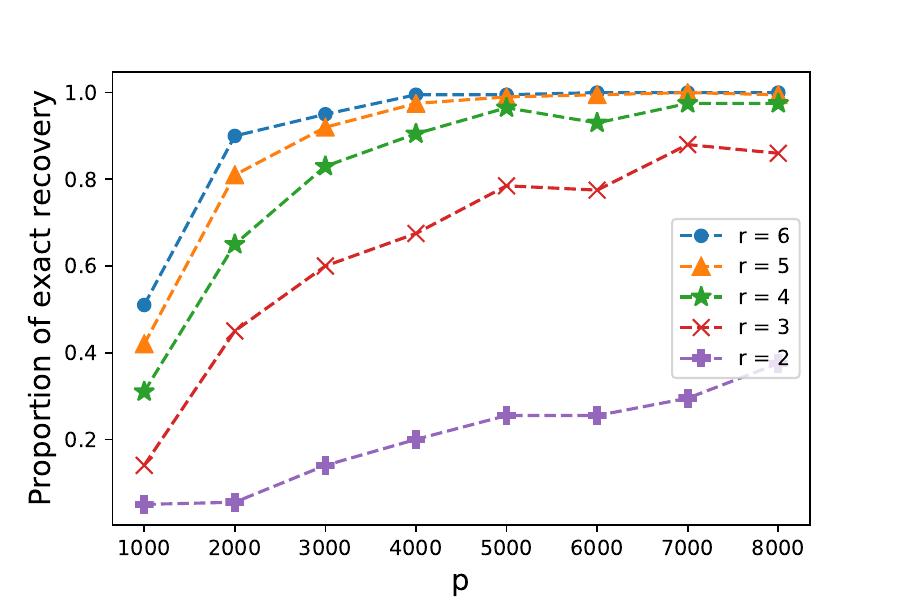}
    \caption{Proportion of exact recovery under weak homogeneous signal. ($0.19 \leq \text{SNR} \leq 2.15$)}
%\label{fig:1a}
    \end{subfigure}\hfill
    \begin{subfigure}{0.47\linewidth}
\includegraphics[width=\linewidth]{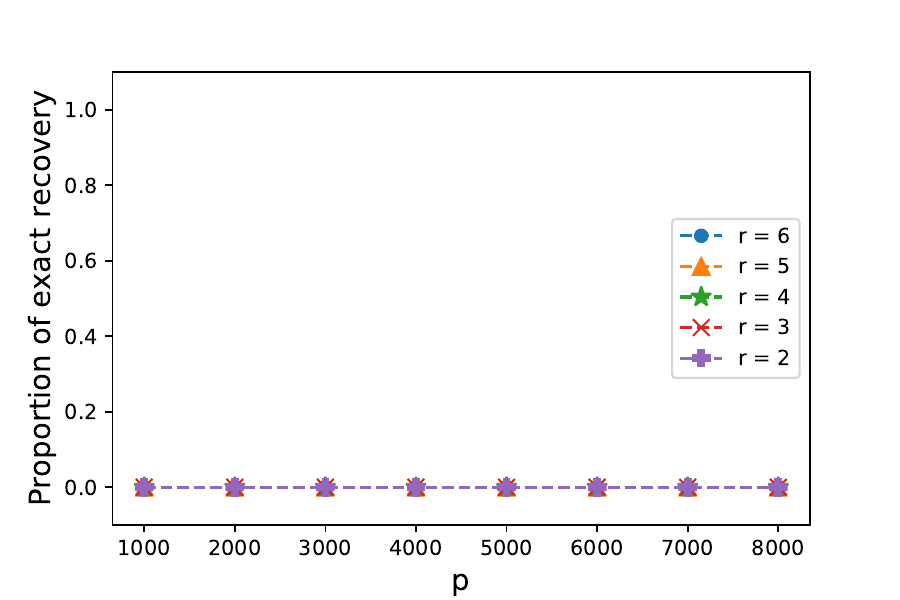}
    \caption{Proportion of exact recovery under heterogeneous signal. (SNR = 10)}
%\label{fig:1a}
    \end{subfigure}
    
\caption{Asymptotics of MS with growing dimension $p$.}
    \label{fig: prop_recovery_path_plot_ms}
    \end{figure}

In Figure \ref{fig: prop_recovery_path_plot_ms}, we demonstrate the asymptotics of MS under both homogeneous and heterogeneous signal patterns. We consider $p \in \{ 1000, 2000, \ldots, 8000\}$ and signal strength parameter $r \in \{2,3,4,5,6\}$. We set $s = \floor{2 \log p}$ and $n = \floor{p^{0.9}}$. We let $\tau(X, Y)$ in \eqref{eq: thresholding estimator} be equal to the $s$th largest value of $\{\abs{\mu_1}, \ldots, \abs{\mu_p}\}$, so that MS always chooses a model of size $s$.
For the homogeneous signal setup, we consider $\beta$ with $\norm{\beta}_0 =s$ and $\beta_j =  a$ for all $j \in \calS(\beta)$, where $a$ is defined in \eqref{eq:aurwm}. This implies that the SNR varies between 0.19 and 2.15 across different choices of $(r,p)$. For the heterogeneous signal setup, we consider $\beta$ with $(s-1)$ active coordinates equal to $a$ and one ``spiky'' coordinate equal to $\{10 - (s-1)a^2\}^{1/2}$. This ensures that the SNR is fixed at 10 for all choices of $r,p$. 

Figure \ref{fig: prop_recovery_path_plot_ms}(a) shows that under homogeneous signal MS is able to recover the exact model with probability converging to 1 as $p$ grows. In contrast, Figure \ref{fig: prop_recovery_path_plot_ms}(b) shows that under heterogeneous signal MS never achieves exact model recovery: plots for all values of $r$ are at level 0. Such a contrast corroborates Theorem \ref{thm: Thresholding fails in ARMW}: signal spikes can give rise to substantial spurious correlation and jeopardize the accuracy of MS. 

\subsection{Effect of growing signal strength}
Here we numerically compare the probability of exact support recovery of ETS with those of LASSO and MS as signal strength parameter $r$ grows. We investigate both homogeneous and heterogeneous signal patterns. We set $p =2000$, $s \in \{13, 52\}$ and $n = \floor{p^{0.9}} = 935$. We set signal strength parameter $r \in \{1.5, 2, 2.5, \ldots, 9\}$ in \eqref{eq:aurwm}. The support $\calS$ is chosen uniformly over all the size-$s$ subsets of $[p]$, and each support coordinate of $\beta$ is chosen as follows:
\[
\beta_j =(1 - b_j) (1 + Z_j^2/n)^{1/2} a + b_j r^{1/2}  \quad \forall j\in \calS,
\]
where $(Z_j)_{j \in \calS} \overset{i.i.d.}{\sim} \sfN(0,1)$, and where $(b_j)_{j \in \calS} \overset{i.i.d.}{\sim} \Ber(\pi)$ with $\pi \in \{0,0.2\}$.  $\pi =0$ corresponds to the homogeneous signal pattern, and $\pi =0.2$ corresponds to the heterogeneous signal pattern,  where spiky signals are present with probability $0.2$. Each entry $x_{ij}$ of the design matrix $X$ is generated independently from $\sfN(0,1)$.

%The error variance is set to be $\sigma =1$ and 
In this experiment, we grant all the approaches with the knowledge of $s$, so that the comparison is fair. Using this oracle knowledge, we only look at the solutions of the aforementioned three methods with sparsity exactly equal to $s$. Specifically, for LASSO, we look at the solution path and select the model of size exactly equal to $s$. For MS, we just select the top $s$ variables corresponding to the largest absolute values of $\mu$'s. 
% set the threshold $\tau(X,Y)$ to be $\abs{\mu}_{(s)}$ where $\abs{\mu}_{(s)}$ denotes the $s$th largest value of $\{\abs{\mu_1}, \ldots, \abs{\mu_p}\}$.  
For ETS, we do not split data for estimation and screening separately; instead we use the full data in both steps. Specifically, we replace $X ^ {(1)}$ and $Y ^ {(1)}$ with $X$ and $Y$ respectively in \eqref{eq: empirical_1_square_error_loss} and replace $X ^ {(2)}$ and $Y ^ {(2)}$ with $X$ and $Y$ respectively in \eqref{eq:delta}. We set gradient step size $h=0.5$ in IHT. We choose projection size $\hat{s}$ by cross validation in terms of mean squared prediction error. Lastly, for selecting exactly $s$ features we use \eqref{eq: selector_ETS_tilde} in the screening stage of ETS. It is worthwhile to mention that from an application point of view, incorporating data splitting in ETS is not necessary as we are only interested in identifying the active signals, which is akin to point estimation. Also, given the fact that $n \ll p$ in high dimensional regime, using full sample in both the estimation and screening step delivers greater sample efficiency and provides better inference. %\zzw{Justify why to drop sample splitting}
% In the screening step, we use the full data version of selector \eqref{eq: selector_ETS_tilde} to perform support recovery, i.e., we replace  and then compute $\tilde \eta$ based on \eqref{eq: selector_ETS_tilde}.
  
\begin{figure}[t]
\tiny
\centering
    \begin{subfigure}{0.44\linewidth}
\includegraphics[width=\linewidth]{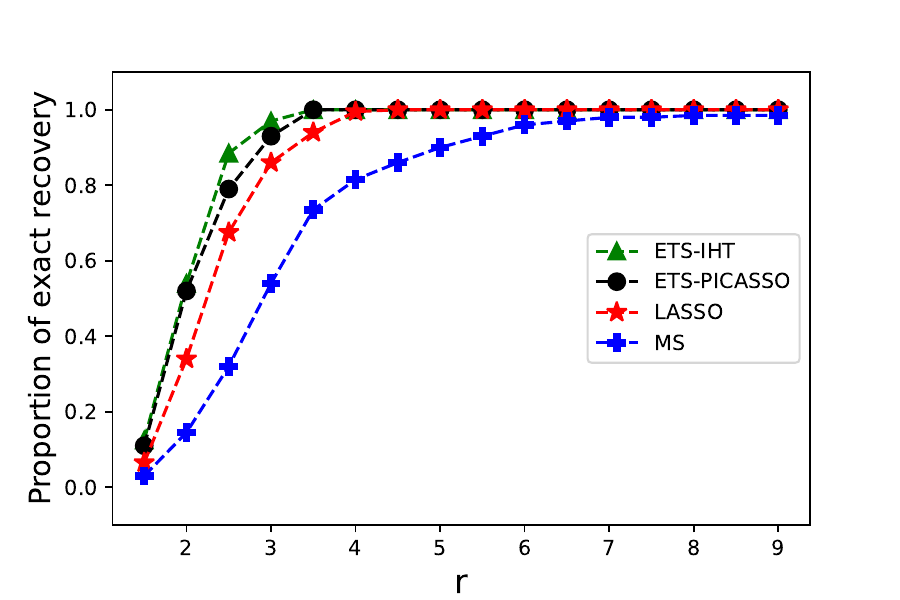} 
    \caption{Homogeneous signal pattern $(\pi =0, s= 13)$}
%\label{fig:1a}
    \end{subfigure}\hfill
    \begin{subfigure}{0.44\linewidth}
\includegraphics[width=\linewidth]{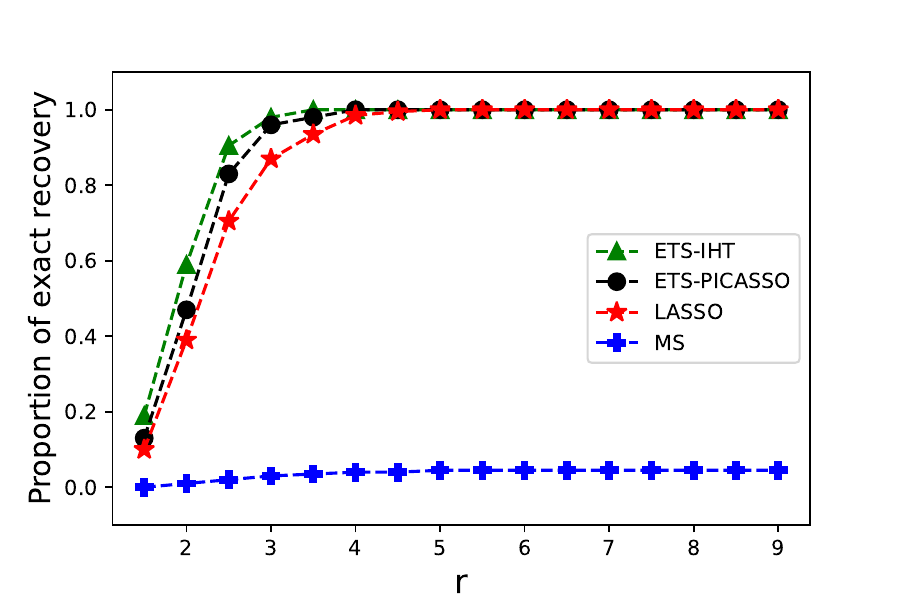}
    \caption{Heterogeneous signal pattern $( \pi =0.2, s=13)$}
%\label{fig:1a}
    \end{subfigure}\hfill
    
    \begin{subfigure}{0.44\linewidth}
\includegraphics[width=\linewidth]{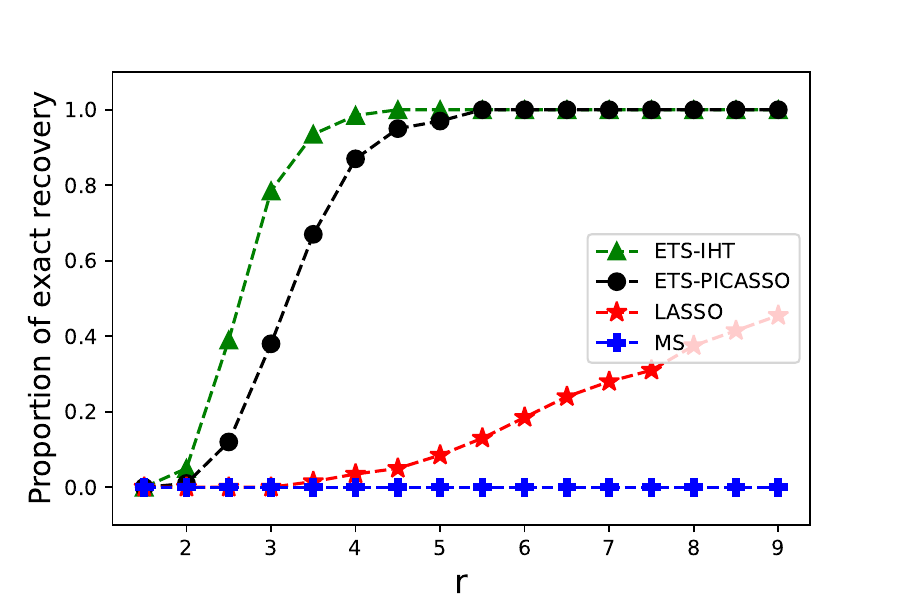}
    \caption{Homogeneous signal pattern $(\pi =0, s= 52)$}
%\label{fig:1a}
    \end{subfigure}\hfill
    \begin{subfigure}{0.44\linewidth}
\includegraphics[width=\linewidth]{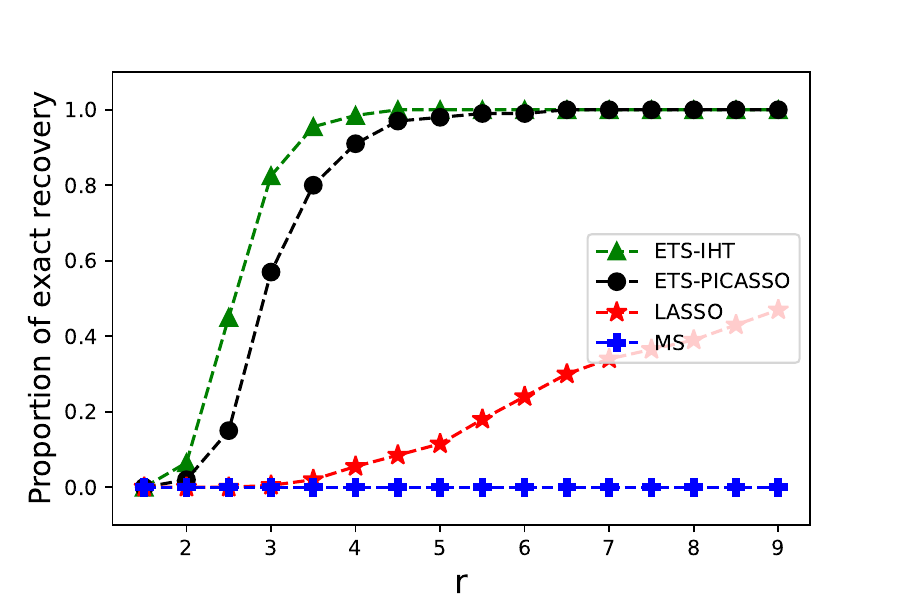}
    \caption{Heterogeneous signal pattern $(\pi =0.2, s= 52)$}
%\label{fig:1a}
    \end{subfigure}
\caption{Plot of the proportion of exact recovery for varying $r$. \sapta{Change the labels of figures}}
    \label{fig: propo_recovery_path_plot}
    \end{figure}

Next, for each choice of $ r$, we run LASSO, MS and ETS over 200 independent Monte Carlo experiments to compute the empirical probability of exact recovery. Figure \ref{fig: propo_recovery_path_plot} presents the results. We make the following important observations: 
\begin{enumerate}
    \item All three methods enjoy a higher chance of exact support recovery as the signal strength grows. 
    \item MS completely fails to achieve exact support recovery when $s$ becomes large (compare panels (a) and (c)) or the signal becomes heterogeneous (compare panels (a) and (b)).
    \item LASSO and ETS algorithms are insensitive to the heterogeneity of the signal. However, LASSO suffers from larger sparsity, while ETS algorithms are much more robust against it.
    \item {\color{black}Overall, ETS is the best among all the three methods in terms of exact support recovery. However, ETS-IHT is somewhat better than ETS-PICASSO, and the difference between their performance is more prominent when $s$ is large.}
\end{enumerate}

{\color{black}
\begin{remark}
    \cite{genovese2012comparison} established model consistency of LASSO under rare and weak signal regime for i.i.d. Gaussian design under a different asymptotic setting and the scaling $s = O(\log p)$ is a special case of their setting. However, they assume homogeneous signal, and our simulation results in Figure \ref{fig: propo_recovery_path_plot}(a) concur with their theoretical findings. However, the performance of LASSO degrades significantly for larger sparsity, even under homogeneous signal (see Figure \ref{fig: propo_recovery_path_plot}(b)), which perhaps shows the limitations of the results in \cite{genovese2012comparison} for realistic values of $p$ and $s$.  
    In contrast, \cite{wainwright2006sharp} obtains the sharpest possible results for model consistency of LASSO under a very general setting. To be precise, under i.i.d. Gaussian design and for general combination of $(n,p,s)$, the paper shows that LASSO achieves model consistency for $a = \Omega(\lambda)$, where $\lambda$ is the regularization parameter and $\lambda\gtrsim \{(\log p)/n\}^{1/2}$. Hence, it is also valid for rare and weak signal regimes and also accommodates heterogeneity in the signal when $\lambda \asymp \{(\log p)/n\}^{1/2}$. However, those are tight only up to multiplicative constants. Moreover, as pointed out in Section B of \cite{wainwright2006sharp}, the model consistency of LASSO depends on whether or not the following is achieved: 
    \[
    \frac{n}{s \log(p-s)} > 1+ \frac{1}{s\lambda^2 }.
    \]
    The above condition is harder to satisfy if $s$ becomes large keeping other parameters fixed, which could be a possible explanation for the phenomenon observed in the third point of the prior observations.
\end{remark}
}
\subsection{Effect of growing heterogeneity}
In this numerical experiment, we study the effect of growing heterogeneity on ETS-IHT, LASSO and MS. We set $p = 2000, s = 13$, $n = \floor{p^{0.9}} = 935$ and $r\in \{2,6\}$ in \eqref{eq:aurwm}. Next, we introduce $\nspike$, the number of ``spiky'' signals in $\calS$. We vary $\nspike$ in $\{0\} \cup [6]$. The case $\nspike = 0$ corresponds to the homogeneous signal setup where the true signals are set as $a$ uniformly. For $\nspike>0$, we randomly set $(s- \nspike)$ signals in $\calS$ to be equal to $a$ and the remaining signals to be equal to $\aspike$, which is defined as 
$$
\aspike  := \bigg\{ \frac{(2 - sa^2)}{\nspike} + a^2\bigg\}^{1/2}.
$$ 
Such a choice of $\aspike$ ensures that the SNR always equals 2 whenever $\nspike > 0$. We perform ETS-IHT, LASSO, and MS over 200 Monte Carlo simulations for each choice of $r$ and $\nspike$ to obtain the empirical probability of exact support recovery. Similarly to the previous sections, we assume that the true sparsity $s$ is known and we apply the three methods in the same fashion as before.

\begin{figure}[h]
\centering
    \begin{subfigure}{0.47\linewidth}
\includegraphics[width=\linewidth]{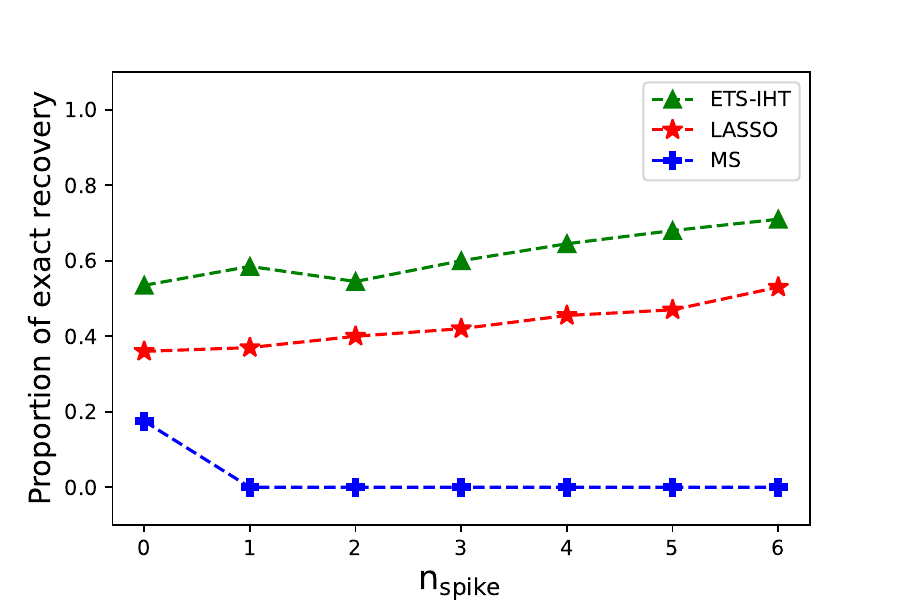} 
    \caption{$r=2$}
%\label{fig:1a}
    \end{subfigure}\hfill
    \begin{subfigure}{0.47\linewidth}
\includegraphics[width=\linewidth]{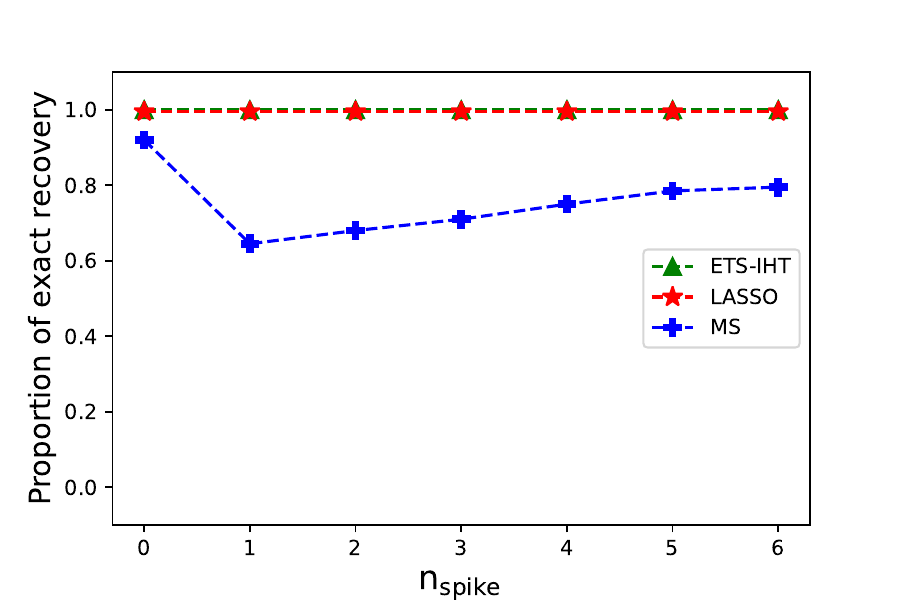}
    \caption{$r=6$}
%\label{fig:1a}
    \end{subfigure}\hfill
    
\caption{Plot of proportion of exact recovery with varying $\nspike$.}
    \label{fig: prop_recovery_path_plot_varying_spike}
    \end{figure}

Figure \ref{fig: prop_recovery_path_plot_varying_spike} shows again the detrimental effect of heterogeneity on MS in terms of exact recovery. In both panels we see a significant drop in the proportion of exact recovery for MS when $\nspike$ changes from 0 to 1. This is consistent with the theory in Section \ref{sec: Thresholding procedures}. However, in Figure \ref{fig: prop_recovery_path_plot_varying_spike}(b) we see that the proportion of exact recovery is slowly increasing as $\nspike$ grows from 1 to 6. This is because as $\nspike$ increases, $\aspike$ monotonically decreases, so that the signals become more homogeneous. MS is then able to recover the exact model more frequently. We do not see a similar phenomenon in Figure \ref{fig: prop_recovery_path_plot_varying_spike}(a) because $\aspike$ is too large. Another important observation is that while ETS-IHT and LASSO are both performing nearly perfectly when $r =6$, ETS-IHT significantly outperforms both LASSO and MS when $r = 2$. Therefore, ETS-IHT is again the overall winner.

\section{Conclusion}
In this paper, we study exact support recovery in high-dimensional sparse linear regression with independent Gaussian design. We focus on the AURWM regime that not only accommodates \textit{rare} and \textit{weak} signals as the ARW regime does, but also allows \textit{heterogeneity} in the signal strength. 
% Under this regime, we do a comparative study of marginal screening, BSS and ETS from an asymptotic viewpoint and present several numerical experiments to support our theory.
Our first theoretical result (Theorem \ref{thm: Thresholding fails in ARMW}) shows that marginal screening fails to achieve exact support recovery under the AURWM regime. 
% We identify \textit{signal heterogeneity} as the main reason behind such failure: 
The main reason is that the presence of ``spiky'' signals increases the maximum spurious marginal correlation, thereby blinding the marginal screening procedure to weak signals. Therefore, one needs to be cautious with usage of marginal screening for variable selection in practice.

% We then move on to study the performance of BSS under the AURWM regime. We obtain the precise asymptotic bound on the signal strength parameter above which BSS achieves model consistency.

In contrast, we show that BSS is robust to signal spikes and is able to achieve model consistency under the AURWM regime with the optimal requirement on signal strength (Theorem \ref{thm: tight phase transition boundary for ML}, \ref{thm: failure of BSS} and Proposition \ref{thm: Information theoretic boundary}). The primary reason behind this is that unlike MS, BSS takes into account multiple features simultaneously and thus selects variables based on their capability of fitting the residualized responses given the other variables rather than the responses themselves. Therefore, spiky signals do not affect BSS: They are very likely to be in plausible candidate models in the first place and their effect on the response has been removed in the residualization procedure. Given the recent computational advancements in solving BSS, our positive result on BSS makes it more appealing from an application point of view. 

However, it is worth mentioning that even with modern advances in optimization, BSS suffers from high computational costs when the ambient dimension is {\color{black}extremely} high.
To address this issue, we propose a computationally tractable two-stage method ETS that delivers essentially the same optimal exact recovery performance as BSS (Theorem \ref{thm: IHT-scrrening exact recovery rate}). Similar to BSS, ETS seeks for the features that exhibit high explanation power for the residuals from the model that excludes these features themselves (see \eqref{eq:delta}). Therefore, ETS is robust to spiky signals. This fact together with the slowly growing sparsity condition in \eqref{eq:aurwm} yields the optimal exact recovery accuracy of ETS.

Our work naturally raises several important questions for future research. One question is whether similar optimality results hold for BSS and ETS when the sparsity $s$ grows faster than $\log p$. The same question can also be asked for \textit{correlated} random design. Another direction of our interest is studying the problem of exact recovery in a \textit{distributed} setting where data are stored at different places and communication between them is restricted.

\bibliographystyle{apalike}  %Style BST file (imsart-number.bst or imsart-nameyear.bst)
\bibliography{bj-ref}       % Bibliography file (usually '*.bib')

% % or include bibliography directly:
% \begin{thebibliography}{}
% \bibitem[\protect\citeauthoryear{???}{???}]{b1}
% \end{thebibliography}

\include{bj-arxiv_supp}
\end{document}

%% file: roy-macro.tex
\usepackage{mathrsfs}

\def\abs#1{\left|#1\right|}

\def\argmin{{\arg\min}}
\def\argmax{{\arg\max}}

\def\bbE{\mathbb{E}}

\def\bbP{\mathbb{P}}

\def\bbR{\mathbb{R}}

\def\sfI{{\sf I}}

\def\sfN{{\sf N}}

\def\Ber{{\sf Ber}}

\def\cA{\mathcal{A}}

\def\cC{\mathcal{C}}
\def\cD{\mathcal{D}}
\def\cE{\mathcal{E}}

\def\cH{\mathcal{H}}
\def\cI{\mathcal{I}}

\def\cM{\mathcal{M}}

\def\cP{\mathcal{P}}

\def\cS{\mathcal{S}}
\def\cT{\mathcal{T}}

\def\norm#1{\left\|#1\right\|}
\def\Norm#1{\Vert#1\Vert}
\def\innerprod#1{\left\langle#1\right\rangle}

\def\pr{{\bbP}}

\def\op{{\rm op}}

\def\emptyset{\varnothing}

\def\betapic{{\widehat{\beta}^{\sf{picasso}}}}
\def\betapgh{{\widehat{\beta}^{\sf{pgh}}}}

\def\Ber{{\sf {Ber}}}

\def\sfN{{\sf N}}
\def\sfI{{\sf I}}
\def\nspike{{\mathsf{n}_{\rm spike}}}
\def\aspike{{a_{\rm spike}}}

\def\barPhi{\overline{\Phi}}

 \DeclareMathOperator{\E}{\mathbb{E}}

 \DeclareMathOperator{\T}{\mathcal{T}}
 \DeclareMathOperator{\calP}{\mathcal{P}}
\def\M{{\mathcal{M}}}
\def\R{{\mathbb{R}}}
 \DeclareMathOperator{\calE}{\mathcal{E}}
\DeclareMathOperator{\V}{\mathcal{V}}
\DeclareMathOperator{\D}{\mathcal{D}}

\DeclareMathOperator{\I}{\mathcal{I}}
\DeclareMathOperator{\J}{\mathcal{J}}
\DeclareMathOperator{\his}{\mathcal{H}}
 \DeclareMathOperator{\calS}{\mathcal{S}}
\usepackage{bbm}
\def\ind{{\mathbbm{1}}}

\newcommand{\floor}[1]{\left\lfloor#1\right\rfloor}
\mathchardef\mhyphen="2D

%% file: bj-arxiv_supp.tex
\begin{appendices}
\section*{Supplementary material}
This section collects the proofs of all the main theorems, propositions, and corollaries. We also discuss about  Remark 5.4 and  prove some important results related to the examples of ETS. 
% Specifically, we present the proofs of Theorem \ref{thm: Thresholding fails in ARMW} and Theorem \ref{thm: tight phase transition boundary for ML}. The first one shows that marginal screening fails to achieve exact recovery in the AURWM regime, whereas BSS achieves model consistency under optimal requirement on the signal strength. 
\section{Proof of Theorem 3.1}
\label{sec: Proof of thm: Thresholding fails in ARMW }
Consider a MS procedure $\widehat{\calS}_\tau \in \T$. Now, there are mainly two steps of the proof:
\begin{enumerate}
    \item Upper bound the probability of recovery in terms of the probability of an event depending only on $\max_{j \in \cS^c} \abs{\mu_j}$ and $\min_{j \in \cS} \abs{\mu_j}$.
    \item Find the asymptotic limits of the above random variables and find the limiting probability of the aforementioned event.
\end{enumerate}

To start with note that
\[
\pr_\beta(\widehat{\calS}_\tau = \calS(\beta) ) = \pr_\beta\left(\max_{j \in \cS^c} \abs{\mu_j} < \tau(X,Y) \leq \min_{j \in \cS} \abs{\mu_j} \right)\leq  \pr \left(\max_{j \in \cS^c} \abs{\mu_j} < \min_{j \in \cS} \abs{\mu_j} \right).
\]
Recall that for for $j \in [p]$ we have 
$$
\mu_j =
\begin{cases}
\beta_j \norm{X_j}^2_2/n + \omega_j  \norm{X_j}_2 g_j/n & \text{if $j \in \calS(\beta)$}\\
 \omega_j \norm{X_j}_2 g_j/n & \text{if  $j \notin \calS(\beta)$}, 
\end{cases}
$$
where $\omega_j^2  = 1 + \sum_{k \neq j}\beta_k^2 $ and $g_j= X_j^\top(\sum_{k\neq j} X_k \beta_k + E
) /(\omega_j \norm{X_j}_2)\sim \sfN(0,1) $. Thus we have 

\begin{equation}
\pr_\beta(\widehat{\calS}_\tau = S) \leq \pr_\beta\left(\max_{j \in S^c}{\omega_j \norm{X_j}_2 \abs{g_j}}/(2n \log p )^{1/2} < \min_{j\in S} \vert\beta_j \norm{X_j}_2^2 + \omega_j \norm{X_j}_2 g_j\vert/( 2 n \log p )^{1/2} \right). 
\label{eq: upper bound accuracy}
\end{equation}

Right hand side of Equation \eqref{eq: upper bound accuracy} does not depend on $\widehat{\calS}_\tau$ hence the above inequality is valid uniformly over the class $\T$. Now choose a sequence $\{c_p\}_{p=1}^\infty$ such that $\lim_{p \to \infty}c_p^2/r \geq 1$ . %$c_p \uparrow\infty$.
Next construct a sequence of $\beta^{(p)}$ in the following manner:
\begin{itemize}
    \item Consider the set $\calS_0 =  \{1,\ldots, s\} \subseteq [p]$ with $s=O(\log p)$. 
    \item Set $\beta_{1}^{(p)} = c_p$. For all other $i\in \calS_0\setminus \{1 \}$ set $\beta_i^{(p)} = a = (2r (\log p)/n)^{1/2}$.
    
    \item Set $\beta_i^{(p)} =0$ if $i \notin \calS_0$.
\end{itemize}
In this setup we have $\omega_j \sim (1 +c_p^2)^{1/2}$ for all $j \neq 1$.
Now fix $k_0\in \calS_0\setminus\{1\} $ (say $k_0 = 2$). From Equation \eqref{eq: upper bound accuracy} it can be concluded that
\begin{align*}
&\sup_{\widehat{\calS}_\tau \in \T}\pr_{\beta^{(p)}}(\widehat{\calS}_\tau = \calS_0) \\
& \leq \pr_{\beta^{(p)}}\left(\max_{j \in S^c}{\omega_j \norm{X_j}_2 \abs{g_j}}/(2n \log p )^{1/2} < \vert \beta_{k_0} \norm{X_{k_0}}_2^2 + \omega_{k_0} \norm{X_{k_0}}_2 g_{k_0} \vert/(2 n \log p )^{1/2} \right). 
\end{align*}
 Also note that $\omega_j > (1+c_p^2)^{1/2}$ for all $j \in \calS_0^c$. Using these facts and  lemma 3 from \cite{fletcher2009necessary} we get 
\[
\frac{1}{(1 + c_p^2)^{1/2}}\max_{j \in \calS_0^c} \omega_j \frac{\norm{X_j}_2 \abs{g_j}}{(2n \log p)^{1/2}}\geq \min_{j \in \calS_0^c}  \frac{\norm{X_j}_2}{n^{1/2}} \max_{j \in \calS_0^c} \frac{\abs{g_j}}{(2 \log p)^{1/2}}. %\quad (\text{as $\omega_j / (1+c_p^2)^{1/2} >1$ for all $j \in S^c$})
\]
Now, note that for $j \in \cS_0^c$, we have all $\omega_j^2 = 1+ \norm{\beta}_2^2$ and 
\[
g_j = \frac{X_j^\top z}{\norm{X_j}_2},
\]
where $z  = (\sum_{k \in \cS_0} X_k \beta_k + E)/\omega_j \sim \sfN_n(0,\sfI_n)$ and it is independent of $\{X_j\}_{j \in \cS_0^c}$. Thus, $g_j^2  = \norm{P_j z}_2^2$, where $P_j$ is the orthogonal projection operator onto the subspace $\text{span}\{X_j\}$. This shows that the random quantity $\max_{j \in \cS_0^c} g_j^2$ is the scaled version of maximum spurious correlation (defined in Section 7.2 of \cite{fan2018discoveries}) between $\{X_j\}_{j \in \cS_0^c}$ and the noise $z$ with sparsity level 1, i.e., $\max_{j \in \cS_0^c} g_j^2 = n \widehat{R}_n^2(1,p-s)$, where

\begin{equation*}
    \widehat{R}_n(1,p-s):= \sup_{\alpha: \norm{\alpha}_2 = 1, \norm{\alpha}_0 =1} \frac{1}{n} \sum_{i=1}^n \frac{\alpha^\top (z_i x_{i, \cS_0^c})}{(\alpha^\top \widehat{\boldsymbol{\Sigma}}_{n,\cS_0^c} \alpha)^{1/2}},
\end{equation*}
with $\widehat{\boldsymbol{\Sigma}}_{n,\cS_0^c} = n^{-1} \sum_{i =1}^n x_{i,\cS_0^c}x_{i, \cS_0^c}^\top$.
 Thus, following the arguments of \cite{fan2018discoveries}, in particular, using Theorem 3.1 and Remark 3.3 of \cite{fan2018discoveries} we have,
 %\newpage
\[
\max_{j \in \calS_0^c} \frac{\abs{g_j}}{(2 \log p)^{1/2}} \overset{p}{\to} 1.
\]
This gives us
\[
\frac{1}{(1 + c_p^2)^{1/2}}\max_{j \in \calS_0^c} \omega_j \frac{\norm{X_j}_2 \abs{g_j}}{(2n \log p)^{1/2}}\geq \min_{j \in \calS_0^c}  \frac{\norm{X_j}_2}{n^{1/2}} \max_{j \in \calS_0^c} \frac{\abs{g_j}}{(2 \log p)^{1/2}} \overset{p}{\to} 1. %\quad (\text{as $\omega_j / (1+c_p^2)^{1/2} >1$ for all $j \in S^c$})
\]
Now, recall that $\beta_{k_0} = \{2r (\log p)/n\}^{1/2}$ and and define $u_r:= \frac{1}{2}(1 + \frac{r^{1/2}}{\sqrt{1+r}}) <1$ .Thus we have the following:
$$ \pr( \vert \beta_{k_0} \norm{X_{k_0}}_2^2 + \omega_{k_0} \norm{X_{k_0}}_2 g_{k_0}\vert / \{(1 + c_p^2) (2 n \log p )\}^{1/2} < u_r) \to 1 .$$ This tells that,
\[
\lim_{p \to \infty} \sup_{\widehat{\calS}_\tau \in \T} \inf_{\beta \in \M_{s}^a} \pr_{\beta} (\widehat{\calS}_\tau = \calS(\beta))\leq \lim_{p\to \infty}\sup_{\widehat{\calS}_\tau \in \T}\pr_{\beta^{(p)}}(\widehat{\calS}_\tau = \calS_0) = 0.
\]
This finishes the proof.

\section{Best subset selection}
In this section, we prove the main results related to BSS.

\subsection{Proof of Theorem 4.1}
\label{sec: proof of ind-design}
In this proof, we reparametrize $\delta$ by $8\delta_0$ for algebraic convenience. The main result can be salvaged by back substituting $8\delta_0$ by $\delta$ in all the main equations in this section. Also, for brevity of notation, in this proof we use $\widehat{\calS}$ and $\calS$ to denote that oracle BSS estimator and $\calS(\beta)$ respectively. We highlight the three main steps of the proof:
\begin{enumerate}
    \item Convert the BSS problem in the problem of selecting the model with maximum spurious correlation.

    \item Use results from \cite{fan2018discoveries} to find the asymptotic distribution of the maximum spurious correlation statistics.

    \item Use the asymptotic distribution along with non-asymptotic concentration inequalities to upper bound the error probability.
\end{enumerate}

Recall that BSS is defined as
\vspace{-.8mm}
\[
\widehat{\calS} = \argmax_{\D, \abs{\D}=s} \norm{P_{\D} Y}_2^2 = \argmin_{\D:\abs{\D}=s} \norm{(\sfI_n-P_{\D})Y}_2^2.
\]
Thus from the above definition we have the following equality:
\begin{align*}
    \pr(\widehat{\calS}\neq \calS) & = \pr \left(\norm{P_{\calS}  Y }_2^2 < \max_{\D\neq \calS} \norm{P_{\D} Y}_2^2 \right).
\end{align*}
Now we will try to understand how the quantity $\norm{P_{\calS} Y}_2^2$ behaves asymptotically. First it is easy to see that $P_{\calS} Y  = \sum_{j \in \calS}X_j \beta_j + P_{\calS} E$. Note that $\sum_{j \in \calS} X_j \beta_j = \norm{\beta}_2 \tilde{\varepsilon} $ where $\tilde{\varepsilon} \sim \sfN_n(0, \sfI_{ n})$ and independent of the noise $z$. Hence we  up with the following:
\begin{align*}
    \norm{P_{\calS} Y }_2^2& = \Big\Vert\sum_{j \in S}X_j \beta_j \Big\Vert_2^2 + 2 E^\top P_{\calS} \left( \sum_{j \in S} X_j \beta_j\right) + E^\top P_{\calS} E\\
    & = \Big\Vert\sum_{j \in S}X_j \beta_j \Big\Vert_2^2 + 2 E^\top \left(\sum_{j \in S} X_j\beta_j \right) + E^\top P_{\calS} E\\
    & = \norm{\beta}_2^2 \norm{\tilde{\varepsilon}}_2^2 + 2 \norm{\beta}_2 E^\top \tilde{\varepsilon} + E^\top P_{\calS} E.
\end{align*}
Recall that $\abs{\beta_j} \geq \{2r (\log p) /n\}^{1/2}$ for all $j \in \calS$. This is presumably the hardest setup as increasing signal strength can only decrease the error probability. Then $\norm{\beta}_2^2 \geq  (2rs \log p)/n$. also note that $E^\top P_{\calS} E\sim \chi^2_s$. Hence we have ,
\[
\frac{\norm{P_{\calS} Y }_2^2}{s \log p} \geq  2r \frac{\norm{\tilde{\varepsilon}}_2^2}{n} + 2 \left(\frac{2r}{s \log p}\right)^{1/2} \frac{\tilde{\varepsilon}^\top E}{n^{1/2}} + \frac{E^\top P_{\calS} E}{s \log p} \overset{\rm p}{\longrightarrow} 2r . 
\]
Thus $\lim_{p \to \infty}\pr(\widehat{\calS}\neq \calS)\leq \pr \left ( 2r \leq \limsup_{p\to \infty}\max_{\D \neq \calS} \norm{P_{\D} Y}_2^2/(s \log p) \right) $. The limiting behaviour of the obtained maximal process turns out to be very challenging to analyze and hence we do not directly study this maximal process. Instead we focus on a related maximal process (will be defined shortly) derived from the earlier one and we use the results from \cite{fan2018discoveries} to study its asymptotic behaviour. Now let us denote the set $\calS \cap \D$ by $\I_0$ and $\D \setminus \calS$ by $\I_1$, i.e., $\D = \I_0 \cup \I_1$. Next define the class $\J_{\I_0} = \{ \I_1 \subseteq [p]: \I_1 \cap S = \emptyset, \abs{\I_1 \cup \I_0} = s   \}$ for each $\I_0\subset S$. Note that $0\leq \abs{\I_0}\leq s-1$ from the construction (if $\abs{\I_0}=s$ then $\D= \calS$). The random variable of interest can be rewritten as follows:

\[
\max_{\D\neq \calS} \frac{\norm{P_{\D} Y}_2^2}{s \log p} = \max_{\I_0: \I_0 \subset \calS}\; \max_{\I_1: \I_1 \in \J_{\I_0}} \frac{\norm{P_{\I_0 \cup \I_1} Y}_2^2}{s \log p}.
\]

Using union bound we get,
\begin{equation}
    \pr(\widehat{\calS}\neq \calS)\leq \sum_{\I_0\subset \calS} \pr \left(\max_{\I_1: \I_1\in \J_{\I_0}} \frac{\norm{P_{\I_0 \cup \I_1} Y}_2^2}{s \log p} > \frac{\norm{P_{\calS} Y}_2^2}{s \log p} \right).
    \label{eq: 0-1_loss_decompositon_constant_sparsity}
\end{equation}
Now fix a subset $\I_0$ of the true support $\calS$. Similar to previous section define $\tilde{Y} = Y - X_{\I_0}\beta_{\I_0}$ and this independent of the features in $\I_0 \cup \I_1$. Also we have 
\begin{align*}
&\norm{P_{\I_0 \cup \I_1} Y }_2^2 = \norm{P_{\I_0 \cup \I_1} \tilde{Y}}_2^2 + \norm{X_{\I_0} \beta_{\I_0}}_2^2 + 2 \beta_{\I_0}^\top X_{\I_0}^\top \tilde{Y},\\
&\hspace{.78cm}\norm{P_{\calS} Y }_2^2 = \norm{P_{\calS} \tilde{Y}}_2^2 + \norm{X_{\I_0} \beta_{\I_0}}_2^2 + 2 \beta_{\I_0}^\top X_{\I_0}^\top \tilde{Y}.
\end{align*}
Thus the summands in the right hand side of \eqref{eq: 0-1_loss_decompositon_constant_sparsity} can be written as the probability of the event $\{\max_{\I_1: \I_1\in \J_{\I_0}} \norm{P_{\I_0 \cup \I_1} g}_2^2/(s\log p) > \norm{P_{\calS} g}_2^2/(s \log p)  \}$, where $g := (1 + \norm{\beta_{\calS \setminus\I_0}}_2^2)^{-1/2} \tilde{Y}$. Note that $g\sim \sfN_n(0,\sfI_n)$ and is independent of the features in $\cD$.
Now fix a specific $\I_0$. In the analysis we encounter the maximal process
\[
\max_{\I_1: \I_1\in \J_{\I_0}} \frac{\norm{P_{\I_0\cup \I_1} g}_2^2}{s\log p},
\]
Now consider the set of indices $F_{\I_0} = (\{ 1,\cdots, p\}\setminus \calS)\cup \I_0$. Hence it is easy to see that $\tilde{p}:=\abs{F_{\I_0}} = p -s + \abs{\I_0}$. Without loss of generality, let $F_{\I_0} = \{1, \ldots, \tilde{p}\}$. Also define $\tilde{s}:= s - \abs{\I_0}$. Let the set $\V_{\I_0} = \{\alpha \in \R^p: \norm{\alpha}_0= s, \norm{\alpha}_2=1, \I_0\subseteq \calS(\alpha), \alpha_{F_{\I_0}^c} = 0\}$. Here $\alpha_J$ denotes the sub-vector of $\alpha$ corresponding to the indices in $J \subseteq [p]$. Next we will focus on the random variable,
\begin{equation}
    \widehat{L}_n:=\widehat{L}_n(\tilde{s},\tilde{p})= \sup_{\alpha \in \V_{\I_0}} \frac{1}{n^{1/2}} \sum_{i=1}^n \frac{\alpha^\top (g_i x_i)}{(\alpha^\top \widehat{\boldsymbol{\Sigma}}_n \alpha)^{1/2}},
    \label{eq: hat_Ln}
\end{equation}
here $\widehat{\boldsymbol{\Sigma}}_n =\frac{1}{n} \sum_{i=1}^n x_i x_i^\top$. Now recall that $\D = \I_0 \cup \I_1$ for all $\D\neq \calS$ with $\abs{\D} = s$. To see the connection, first note that the above optimization problem can be viewed as the following:
\begin{align*}
\widehat{L}_n &= \max_{\I_1\in \J_{\I_0}}\max_{\alpha \in \V_{\I_0 \cup \I_1}} \frac{\alpha_{\D}^\top (\sum_{i=1}^n g_i x_{i,\D}/n^{1/2}) }{\{\alpha_{\D}^\top (\widehat{\boldsymbol{\Sigma}}_{n,\D \D}) \alpha_{\D}\}^{1/2} }\\
&= \max_{\I_1\in \J_{\I_0}} \left\{(\sum_{i=1}^n g_i x_{i,\D}/n^{1/2})^\top  \widehat{\boldsymbol{\Sigma}}^{-1}_{n,\D \D}(\sum_{i=1}^n g_ix_{i,\D}/n^{1/2}) \right\}^{1/2}\\
&= \max_{\I_1\in \J_{\I_0}} \{g^\top X_{\D} (X_{\D} ^\top X_{\D})^{-1} X_{\D}^\top g\}^{1/2}\\
& = \max_{\I_1\in \J_{\I_0}} \norm{P_{\I_0\cup \I_1
} g}_2.
\end{align*}
Thus it is essential to study the asymptotic property of $\widehat{L}_n$. Now we define the standardized version of $\widehat{L}_n$ as follows
\[
L_n :=L_n(\tilde{s},\tilde{p}) =  \sup_{\alpha \in \V_{\I_0}} \frac{1}{n^{1/2}} \sum_{i=1}^n \alpha^\top (g_i x_i).
\]
Let $\boldsymbol{Z}=(Z_1,\cdots, Z_{\tilde{p}})$ be $\tilde{p}-$variate Gaussian random variable with covariance matrix $I_{\tilde{p}\times \tilde{p}}$ and define the random variable $T^*:= T^*(\tilde{s},\tilde{p}) = \sup_{\alpha \in \V_{\I_0}} \alpha_{F_{\I_0}}^\top \boldsymbol{Z}$.

\begin{lemma}
There exists universal constants $K_0,K_1$ such that for any $\delta_1\in(0, K_0 K_1]$,
\begin{equation}
    \abs{L_n - T^*}\lesssim n^{-1} c_n^{1/2}(\tilde{s},\tilde{p}) + K_0 K_1 n^{-3/2} c_n^2(\tilde{s},\tilde{p}) +\delta_1
    \label{eq: approximating-L_n&T*}
\end{equation}
holds with probability at least $1- C\Delta_n(s, \tilde{p};\delta_1)$ where $c_n(s,\tilde{p})= s \log(e \tilde{p}/s)\vee \log n$ and
\[
\Delta_n(\tilde{s},\tilde{p};\delta_1) = (K_0 K_1)^3 \frac{\{ \tilde{s} b_n(\tilde{s},\tilde{p})\}^2}{\delta_1^3 n^{1/2}} +  (K_0 K_1)^4 \frac{\{ \tilde{s} b_n(\tilde{s},\tilde{p})\}^5}{\delta_1^4 n}
\]
with $b_n(\tilde{s}, \tilde{p}) = \log(\tilde{p}/\tilde{s})\vee \log n$.
\label{lemma: approximating-L_n&T*}
\end{lemma}

\begin{lemma}
Assume that the sample size satisfies $n\geq C_1(K_0\vee K_1)^4 c_n(\tilde{s},\tilde{p})$. then with probability at least $1- C_2 n^{-1/2} c_n^{1/2}(\tilde{s},\tilde{p})$,

\begin{equation}
    \abs{\widehat{L}_n - L_n}\lesssim (K_0\vee K_1)^2 K_0 K_1 n^{-1/2} c_n(\tilde{s},\tilde{p}),
    \label{eq: approximating-L_n&L_n_hat}
\end{equation}
where $c_n(\tilde{s},\tilde{p}) = \tilde{s}\log(e \tilde{p}/\tilde{s})\vee \log n$.
\label{lemma: approximating-L_n&L_n_hat}
\end{lemma}

Proof of the above two lemmas are omitted as it is in the same line of the proofs of \cite{fan2018discoveries}.
Now applying Lemma \ref{lemma: approximating-L_n&T*} and \ref{lemma: approximating-L_n&L_n_hat} with $$\delta_1 = \delta_1(s,\tilde{p}) = (K_0 K_1)^{3/4} \min[1, n^{-1/8}\{ \tilde{s} b_n(\tilde{s},\tilde{p})^{3/8}\}]$$ yields that  with probability at least $1- C (K_0K_1)^{3/4} n^{-1/8} \{s b_n(s,\tilde{p})\}^{7/8}$,
\[
\abs{\widehat{L}_n - T^*}\lesssim  (K_0K_1)^{3/4} n^{-1/8} \{\tilde{s} b_n(\tilde{s},\tilde{p})\}^{3/8}.
\]
 Together with Lemma 2.3 from \cite{chernozhukov2014gaussian} we can conclude that 
 \begin{equation}
     \sup_{t\in \R} \abs{\pr(\widehat{L}_n \leq t) - \pr(T^*\leq t)}\lesssim C (K_0K_1)^{3/4} n^{-1/8} \{\tilde{s} b_n(\tilde{s},\tilde{p})\}^{7/8}.
     \label{eq: KS_distance-L_n&T*}
 \end{equation}
 Next by the definition of $T^*$ it follows that 
 $$T^{*2} = \max_{\I_1\in \J_{\I_0}} \norm{\boldsymbol{Z}_{\I_0 \cup \I_1}}_2^2=\sum_{j \in \I_0} Z_j^2 + \max_{\I_1\in \J_{\I_0}} \sum_{k \in \I_1} Z_k^2.$$
 
 Let $\boldsymbol{W}\sim \sfN_{(p-s)}(0, \sfI_{(p-s)})$ be a Gaussian vector  independent of $\boldsymbol{Z}$. Thus it follows that
 \[
 T^{*2}\overset{d}{=} \sum_{j \in \I_0} Z_j^2 + \sum_{k= p-2s+\abs{\I_0}+1}^{p-s} W^2_{(k:p-s)}\leq \sum_{j \in \I_0} Z_j^2 + (s- \abs{\I_0}) W^2_{(p-s:p-s)}.
 \]
 From Equation \eqref{eq: KS_distance-L_n&T*} it also follows that
 \begin{equation}
     \sup_{t\geq 0} \abs{\pr(\widehat{L}^2_n \leq t) - \pr(T^{*2}\leq t)}\lesssim C (K_0K_1)^{3/4} n^{-1/8} \{\tilde{s} b_n(\tilde{s},\tilde{p})\}^{7/8}.
     \label{eq: KS_distance- square_L_n&T*}
 \end{equation}
 
 Now from the assumption, we have $r = 1+ 8 \delta_0$. Assume that 
 \begin{equation}
 \label{eq: sparsity_upper_bound}
 s\leq 0.5\min\{\delta_0, \frac{2\delta_0^2}{\{(1+ 6\delta_0)^{1/2}+ (1+ 4\delta_0)^{1/2}\}^2}\} \log p.
 \end{equation}
 Hence $\abs{\I_0}\leq s -1 \leq 0.5\delta_0 \log p < \delta_0 \log p$. Thus we have,
 \begin{align*}
 \scriptstyle
     &\quad \pr(T^{*2} > 2(1+ 4\delta_0) (s-\abs{\I_0}) \log p) \\
     &= \pr(\sum_{j \in \I_0} Z_j^2 + (s- \abs{\I_0})W^2_{(p-s:p-s)}> 2(1+ 4\delta_0) (s-\abs{\I_0}) \log p)\\
     &\leq \pr \left(\sum_{j \in \I_0} Z_j^2 > 4\delta_0 (s-\abs{\I_0}) \log p \right) + \pr \left((s- \abs{\I_0}) W^2_{(p-s:p-s)} > 2(1+ 2\delta_0)(s-\abs{\I_0}) \log p \right)\\
     & \overset{(a)}{\leq} \pr \left( \frac{\sum_{j \in \I_0} Z_j^2 - \abs{\I_0}}{\abs{\I_0}} > (4  \delta_0\log p - \abs{\I_0})/\abs{\I_0} \right)  + \pr\left( W^2_{(p-s:p-s)} > 2(1+ 2\delta_0)\log p   \right)\\
     & \overset{(b)}{\leq} \pr \left( \frac{\sum_{j \in \I_0} Z_j^2 - \abs{\I_0}}{\abs{\I_0}} > (3 \delta_0\log p)/\abs{\I_0} \right)  + \pr\left( W^2_{(p-s:p-s)} > 2(1+ 2\delta_0)\log p \right)\\
     & \lesssim  \exp(- 0.75\delta_0 \log p) + (p - s)\pr\left( W^2_1 > 2(1+ 2\delta_0)\log p \right)\\ % 1- \left(1- \sqrt{\frac{1}{\pi (1+ 2\delta_0) \log p}} \frac{1}{p^{1+2\delta_0}} \right)^{p-s}\\
     & \lesssim p^{-0.75\delta_0} + C \frac{p^{-2\delta_0}}{\sqrt{\log p}}.
 \end{align*}
 Inequality $(a)$ uses $s - \abs{\I_0}\geq 1$ and inequality $(b)$ uses $\abs{\I_0}<s< \delta_0 \log p$ (Condition \eqref{eq: sparsity_upper_bound}). Also, the first probability bound in (b) follows from Equation (56) in \cite{wainwright2009information} and the fact that $(3 \delta_0 \log p)/\abs{\cI_0} \geq 6>4$. The last inequality in (b) follows from tail bound of standard Gaussian distribution.
 Now define the event $ \calE_{\I_0} := \{  \norm{P_{\calS} {g}}_2^2/(s\log p) > 2(1 + 4\delta_0)R_{\I_0}\}$ where $R_{\I_0}:= (s- \abs{\I_0})/s$.  Recall that 
 \begin{equation}
\begin{split}
    \frac{\norm{P_{\calS} {g}}_2^2}{s\log p} &\geq \frac{\norm{ \frac{ \sum_{j \in \calS \setminus
     \I_0} X_j \beta_j + P_{\calS} E }{(1 + \Vert\beta_{\calS \setminus \I_0}\Vert_2^2)^{1/2}} }_2^2}{s \log p}\\
     & \geq \frac{ \left\{  \frac{ \norm{\sum_{j \in \calS \setminus
     \I_0} X_j \beta_j }_2}{(1 + \Vert\beta_{\calS \setminus \I_0}\Vert_2^2)^{1/2}} -
      \frac{\norm{ P_{\calS} E}_2}{(1 + \Vert\beta_{\calS \setminus \I_0}\Vert_2^2)^{1/2}} \right\}^2   }{s \log p} = (T_1^{1/2}- T_2^{1/2})^2.
\end{split}
\label{eq: E_I_0 lower bound}
\end{equation}
where $$ T_1:=\frac{ \norm{\sum_{j \in \calS \setminus
     \I_0} X_j \beta_j }_2^2}{(1 + \norm{\beta_{\calS \setminus \I_0}}_2^2)s\log p}\geq \frac{2r R_{\I_0}}{1 + 2r (s- \abs{\I_0}) (\log p)/n} \frac{V_n}{n},$$ and $V_n:= \frac{\Vert \sum_{j \in \calS \setminus
     \I_0} X_j \beta_j \Vert_2^2} {\Vert \beta_{\calS \setminus \I_0}\Vert_2^2}$ is an $\chi^2_n$ random variable. Also we have $$ T_2:=\frac{\norm{ P_{\calS} E}^2_2}{(1 + \norm{\beta_{\calS \setminus \I_0}}_2^2) s\log p}\leq V_s/(s\log p)$$ where $V_s:= \norm{P_{\calS} E}_2^2$ is an $\chi^2_s$ random variable independent of $X_S$. Next, we state the following simple algebraic relationship:
     \[
     (1+ 6\delta_0)^{1/2} - \frac{\delta_0}{(1+ 6\delta_0)^{1/2}+ (1+ 4\delta_0)^{1/2}}\geq (1+ 4\delta_0)^{1/2}.
     \]
     In light of Equation \eqref{eq: E_I_0 lower bound} and using the above algebraic inequality we have the following: 
     \[
     \calE_{\I_0}^c \subseteq \{ T_1 \leq 2(1 + 6\delta_0) R_{\I_0} \} \bigcup \{T_2 \geq 2 \delta_0^2 \{(1+ 6\delta_0)^{1/2}+ (1+ 4\delta_0)^{1/2}\}^{-2} R_{\I_0}\}
     \]
     Next, we have
\[
         \pr(T_1 \leq  2(1+ 6\delta_0)R_{\I_0}) \leq \pr \left(\frac{V_n}{n} \leq \frac{1 + 6\delta_0}{1 + 8\delta_0}(1 + 2rs \log p/n)\right).
\]
Now choose large $n$ such that $(1 + 6\delta_0)(1 + 2rs \log p/n))< (1+ 7\delta_0)$. Then for large $n$ we have,
$$
\pr(T_1 \leq  2(1+ 6\delta_0)R_{\I_0})\leq \pr\left(\abs{V_n/n-1}\geq \frac{\delta_0}{1+ 8\delta_0}  \right) \lesssim \exp \left\{-C^* \frac{\delta_0^2}{(1+ 8\delta_0)^2} n \right\},
$$
where $C^*$ is a universal constant.
Now we analyze the quantity $T_2$. 
%\newpage
We have the following inequalities:
\begin{align*}
&\pr(T_2 \geq \frac{2 \delta_0^2}{\{(1+ 6\delta_0)^{1/2}+ (1+ 4\delta_0)^{1/2}\}^2} R_{\I_0})\\
& \leq \pr(V_s/s\geq \frac{2 \delta_0^2}{\{(1+ 6\delta_0)^{1/2}+ (1+ 4\delta_0)^{1/2}\}^2} R_{\I_0}\log p)\\
&\leq \pr(V_s\geq  \frac{2 \delta_0^2}{\{(1+ 6\delta_0)^{1/2}+ (1+ 4\delta_0)^{1/2}\}^2} \log p)\\
&\leq \pr(\abs{V_s/s-1}\geq 0.5 \frac{2 \delta_0^2}{\{(1+ 6\delta_0)^{1/2}+ (1+ 4\delta_0)^{1/2}\}^2} (\log p)/s)\\
& \leq \exp(- C^{\prime} \frac{2 \delta_0^2}{\{(1+ 6\delta_0)^{1/2}+ (1+ 4\delta_0)^{1/2}\}^2} \log p) \quad (\text{Using Condition \eqref{eq: sparsity_upper_bound}})\\
& = p^{-C^{\prime} \frac{2 \delta_0^2}{\{(1+ 6\delta_0)^{1/2}+ (1+ 4\delta_0)^{1/2}\}^2}} \quad (\text{$C^\prime>0$ is universal constant}).
\end{align*}
 Ultimately it shows that $\pr(\calE_{\I_0}^c)\lesssim \exp \left\{-C^* \frac{\delta_0^2}{(1+ 8\delta_0)^2} n \right\}+  p^{-C^{\prime} \frac{2 \delta_0^2}{\{(1+ 6\delta_0)^{1/2}+ (1+ 4\delta_0)^{1/2}\}^2}}$. 
%  This statement uses the fact that 
%      \[
%      \sqrt{1+6\delta_0} - \frac{\delta_0}{\sqrt{1+ 6\delta_0}+ \sqrt{1+ 4\delta_0}}\geq \sqrt{1+4\delta_0}.
%      \]
     Now we are ready to show that the error probability goes to $0$.
     \begin{align*}
    &\pr_\beta(\widehat{\calS}\neq \calS) \leq \sum_{\I_0\subset S} \pr_\beta \left(\max_{\I_1: \I_1\in \J_{\I_0}} \frac{\norm{P_{\I_0 \cup \I_1} Y}_2^2}{s \log p} > \frac{\norm{P_{\calS} Y}_2^2}{s \log p} \right)\\
    &\leq \sum_{k =0}^{s-1}\sum_{\I_0: \abs{\I_0}=k} \pr_\beta \left(\max_{\I_1: \I_1\in \J_{\I_0}} \frac{\norm{P_{\I_0 \cup \I_1} Y}_2^2}{s \log p} > \frac{\norm{P_{\calS} Y}_2^2}{s \log p} \right)\\
    &\leq 
    \sum_{k=0}^{s-1} \sum_{\I_0: \abs{\I_0}=k} \pr_\beta \left(\max_{\I_1: \I_1\in \J_{\I_0}} \frac{\norm{P_{\I_0 \cup \I_1} Y}_2^2}{s \log p} > \frac{\norm{P_{\calS} Y}_2^2}{s \log p} ,\calE_{\I_0}\right) + \pr(\calE_{\I_0}^c)\\
    & 
    \leq 
    \sum_{k=0}^{s-1} \sum_{\I_0: \abs{\I_0}=k} \pr_\beta\left(\max_{\I_1: \I_1\in \J_{\I_0}} \frac{\norm{P_{\I_0 \cup \I_1} g}_2^2}{s \log p} > \frac{\norm{P_{\calS} g}_2^2}{s \log p} ,\calE_{\I_0}\right) + \pr(\calE_{\I_0}^c)\\
    &  \leq 
    \sum_{k=0}^{s-1} \sum_{\I_0: \abs{\I_0}=k} \pr_\beta\left(\max_{\I_1: \I_1\in \J_{\I_0}} \frac{\norm{P_{\I_0 \cup \I_1} g}_2^2}{s \log p}  > 2(1 + 4\delta_0)R_{\I_0}\right) + \pr(\calE_{\I_0}^c)\\
    &  \overset{(a)}{\lesssim}
    \sum_{k=0}^{s-1} \sum_{\I_0: \abs{\I_0}=k}\left[ \pr\left(T^{*2} > 2(1 + 4\delta_0)s R_{\I_0} \log p\right) + \pr(\calE_{\I_0}^c) + C (K_0K_1)^{3/4} n^{-1/8} \{s b_n(s,p)\}^{7/8}\right]\\
    & 
    \lesssim \sum_{k=0}^{s-1} \sum_{\I_0: \abs{\I_0}=k} p^{-0.75\delta_0} + C \frac{p^{-2\delta_0}}{\sqrt{\log p}} + \exp \left\{-C^* \frac{\delta_0^2}{(1+ 8\delta_0)^2} n \right\} +  p^{-C^{\prime} \frac{2 \delta_0^2}{\{(1+ 6\delta_0)^{1/2}+ (1+ 4\delta_0)^{1/2}\}^2}}\\
    & + n^{-1/8} \{s b_n(s,p)\}^{7/8} \\
    & \lesssim \sum_{k=0}^{s-1} \binom{s}{k} \left[ p^{-0.75\delta_0} + \exp \left\{-C^* \frac{\delta_0^2}{(1+ 8\delta_0)^2} n \right\} +  p^{-C^{\prime} \frac{2 \delta_0^2}{\{(1+ 6\delta_0)^{1/2}+ (1+ 4\delta_0)^{1/2}\}^2}}
     + n^{-1/8} \{s b_n(s,p)\}^{7/8}\right]\\
    & \lesssim 2^s \left[  p^{-0.75 \delta_0} + \exp \left\{-C^* \frac{\delta_0^2}{(1+ 8\delta_0)^2} n \right\} +  p^{-C^{\prime} \frac{2 \delta_0^2}{\{(1+ 6\delta_0)^{1/2}+ (1+ 4\delta_0)^{1/2}\}^2}}+ n^{-1/8} \{s b_n(s,p)\}^{7/8}\right].
\end{align*}
Inequality $(a)$ uses $ \tilde{s} b_n(\tilde{s},\tilde{p})\leq s b_n(s,p)$ for large $p$. Thus if 
\begin{align*} s & \lesssim \left( \delta_0  \wedge \frac{2 \delta_0^2}{\{(1+ 6\delta_0)^{1/2}+ (1+ 4\delta_0)^{1/2}\}^2} \wedge \frac{k}{16} \right)\log p\\
&  = \left( \frac{2 \delta_0^2}{\{(1+ 6\delta_0)^{1/2}+ (1+ 4\delta_0)^{1/2}\}^2} \wedge \frac{k}{16} \right)\log p,
\end{align*}
then error probability goes to 0 uniformly over $\beta \in \M_s^a$.

\subsection{Model consistency of BSS for sub-Gaussian model}
\label{sec: BSS result sub-Gaussian}
In this section, we will show that Theorem 4.1 also holds beyond the Gaussian model. We assume that the entries of the design matrix $X$ are i.i.d. \textit{mean-zero} and \textit{sub-Gaussian} with \textit{unit variance}. We also assume that the entries of $E$ are also i.i.d. \textit{mean-zero} and \textit{sub-Gaussian} with \textit{unit variance} and independent of $X$.

In this setup, the results of Theorem 3.1 in \cite{fan2018discoveries} are also valid and the proof steps follow exactly the same steps as the proof of Theorem 4.1 until the introduction of the random variables $V_n$ and $V_s$.

We note that Gaussianity was only used to characterize the distributions of $V_n$ and $V_s$. In particular, due to Gaussianity, we have $V_n \sim \chi^2_n$ and $V_s \sim \chi^2_s$. However, under the sub-Gaussian case, $V_n$ is the sum of $n$ independent sub-Exponential random variables with \textit{unit-mean}. Hence, the probability bound for $T_1$ shown in the original proof is also valid. 

Next, for $V_s$, we can use Theorem 1.1 of \cite{rudelson2013hanson}. Note that, there exists a constant $K_{\psi_2}>1$ such that $\norm{\varepsilon}_{\psi_2}\leq K_{\psi_2}$, where $\norm{\varepsilon}_{\psi_2}:= \inf_{t >0}\{t : \bbE \exp(\varepsilon^2/t^2) \leq 2\}$. By Theorem 1.1 of \cite{rudelson2013hanson}, we can obtain the same probability bound for $T_2$ if 
$$s\leq (0.5/K_{\psi_2}^2)\min\{\delta_0, \frac{2\delta_0^2}{\{(1+ 6\delta_0)^{1/2}+ (1+ 4\delta_0)^{1/2}\}^2}\} \log p.$$
Hence, the rest of the proof is verbatim to the proof in the Gaussian case.

\subsection{Proof of Theorem 4.3} % 4.2
In this section, we will show that BSS fails to recover the exact support when $r=1$. 
We highlight three main steps of the proof:
\begin{enumerate}
    \item Convert the BSS problem in the problem of selecting the model with maximum spurious correlation.

    \item Use results from \cite{fan2018discoveries} to find the asymptotic distribution of the maximum spurious correlation statistics.

    \item Use the exact form of the asymptotic distribution along with scaling and centering parameters to approximate the recovery probability.
\end{enumerate}
Recall the linear model $Y = X \beta + E$ with $\cS:= \cS(\beta)$ as the set of active features and $s = \abs{\cS} = O(\log p)$. As $r = 1$, there exists $j_0 \in \cS$ such that $\beta_{j_0} = \{(2 \log p)/n\}^{1/2}$. WLOG
, let us assume that $j_0 = 1$ and define $\cI_0 = \cS \setminus \{1\}$. In order for BSS to recover the exact support $\cS$, it is necessary that 
\begin{align*}
   & \max_{j \notin \cS } \norm{P_{\cI_0 \cup \{j\}} Y}_2^2 < \norm{P_\cS Y}_2^2\\
   & \Leftrightarrow \max_{j \notin \cS} \norm{P_{\cI_0 \cup \{j\}} \Tilde{Y}}_2^2 < \norm{P_\cS \Tilde{Y}}_2^2 \quad \left(\text{where $\Tilde{Y} =  X_1 \beta_1 + E$}\right)\\
   & \Leftrightarrow \max_{j \notin \cS}\norm{P_j^\perp \Tilde{Y}}_2^2 < \norm{P_1^\perp \Tilde{Y}}_2^2,
\end{align*}
where $P_j^\perp$ is the orthogonal projection operator onto the sub-space $\text{span} \{ \tilde{X}_j\}$ for all $j \in \cS^c \cup \{1\}$, where $\Tilde{X}_j = (\sfI_n - P_{\cI_0}) X_j$. Due to Gaussianity, it follows that $\norm{\tilde{X}_j}_2^2 \sim \chi^2_{n-s+1}$.

Now, note that 
\begin{align*}
\norm{P_1^\perp \Tilde{Y}}_2^2 &= \beta_1^2 \norm{\tilde{X}_1}_2^2 + 2 \beta_1 \tilde{X}_1^\top E + \norm{P_1^\perp E}_2^2\\
& = (2 \log p) \frac{\norm{\tilde{X}_1}_2^2}{n} + 2 (2\log p)^{1/2} \frac{\tilde{X}_1^\top E}{\sqrt{n}} + \norm{P_1^\perp E}_2^2.
\end{align*}

Next, define the events
\[
\cE_1 := \left\{ \abs{\frac{\norm{\tilde{X}_1}_2^2}{n - s + 1} - 1} \leq 1/\sqrt{\log p}\right\}, \quad \cE_2 := \left\{ \frac{\tilde{X}_1^\top E}{\norm{\Tilde{X}_1}_2} \leq  -1\right\}, \quad \cE_3:= \left\{ \norm{P_1^\perp E}_2^2  \leq 16\right\}.
\]
When $p>4$, by Bernstein's inequality, we have $\pr(\cE_1^c) \leq e^{- c_1 \frac{n}{\log p}}$, where $c_1>0$ is a universal constant. Next, note that $\frac{\tilde{X}_1^\top E}{\norm{\Tilde{X}_1}_2} \sim \sfN(0,1)$. Next, we introduce a useful lemma.
\begin{lemma}[\cite{gordon1941values}]
\label{lemma: gaussian tail bounds}
    Let $\Phi(\cdot)$ denote the cumulative distribution function of standard Gaussian distribution. Then for all $x\geq 0$, the following inequalities are true:
    \[
    \left(\frac{x}{1+x^2}\right) \frac{e^{-x^2/2}}{\sqrt{2\pi}} \leq 1- \Phi(x) \leq \left(\frac{1}{x}\right)\frac{e^{-x^2/2}}{\sqrt{2\pi}}.
    \]
\end{lemma}
By the above lemma we can conclude $\pr(\cE_2^c) \leq 1 - \frac{ e^{-1/2}}{ 2\sqrt{2 \pi}}$. Finally, as $\norm{P_1^\perp E}_2^2\sim \chi^2_1$, we have $\pr(\cE_3^c) \leq 2 e^{-8}$. Since $n + 4> 4s$ for large p, we have the following under $\cE_1 \cap \cE_2$:
\[
\frac{\tilde{X}_1^\top E}{\sqrt{n}} = \frac{\tilde{X}_1^\top E}{\norm{\Tilde{X}_1}_2} \times \frac{\norm{\Tilde{X}_1}_2}{\sqrt{n-s+1}} \times \sqrt{\frac{n-s+1}{n}} \le -\frac{\sqrt{3}\{1 - (\log p)^{-1/2}\}}{2}.
\]
Here we used the fact that $\sqrt{1 - (\log p)^{-1/2}} > 1 - (\log p)^{-1/2}$.
We define the event $\cE:= \cap_{i =1}^3 \cE_i$. Then, we have 
%\newpage
\begin{align*}
    &\pr(\widehat{\cS}_{\rm best} = \cS)\\
    & \leq \pr \left(\max_{j \notin \cS}\norm{P_j^\perp \Tilde{Y}}_2^2 < \norm{P_1^\perp \Tilde{Y}}_2^2\right)\\
    & \leq \pr \left(\max_{j \notin \cS}\norm{P_j^\perp \Tilde{Y}}_2^2 < \norm{P_1^\perp \Tilde{Y}}_2^2, \cE\right) + \pr(\cE^c)\\
    & \leq  \pr \left\{\max_{j \notin \cS}\norm{P_j^\perp \Tilde{Y}}_2^2 < 2 \log p \left( 1  + \frac{1}{\sqrt{\log p}} \right) -  (6 \log p)^{1/2}\left(1 - \frac{1}{\sqrt{\log p}}\right) + 16 \right\} + \pr(\cE^c).
\end{align*}

We further note that $g := \tilde{Y}/(1 + \beta_1^2)^{1/2}$ follows a standard isotropic Gaussian distribution. Using Theorem 3.1 of \cite{fan2018discoveries} and the fact that $1 + \beta_1^2>1$, we get 
\begin{align*}
   &\pr(\widehat{\cS}_{\rm best} = \cS)\\
   & \leq \pr \left\{
   Z^2_{(p-s:p-s)} \leq 2 \log p \left( 1  + \frac{1}{\sqrt{\log p}} \right) -  (6 \log p)^{1/2}\left(1 - \frac{1}{\sqrt{\log p}}\right) + 16 
   \right\} + \pr(\cE^c) + o(1)\\
   & \leq \pr \left\{
Z^2_{(p-s:p-s)} - 2 \log (p-s) + \log \log (p-s)  \leq 
\underbrace{- \left( \sqrt{6} - 2 \right) (\log p)^{1/2} + \log \log p +   O(1)}_{:= t_p}
\right\}\\
& \quad + \pr(\cE^c) + o(1).
\end{align*}

where $Z^2_{(p-s:p-s)}$ is the maximum order statistics of $\{Z_j^2\}_{j \in [p-s]}$ with $\{Z_j\}_{j \in [p-s]}$ being i.i.d. standard Gaussian. Finally from Remark 3.3 of \cite{fan2018discoveries}, we know 
\[
Z^2_{(p-s:p-s)} - 2 \log (p-s) + \log \log (p-s) \overset{\rm d}{\to} \Lambda,
\]

where $\pr (\Lambda \leq t)  = \exp(- \pi^{-1/2} \exp(- t/2))$. As $t_p \to - \infty$, we have 

\[
\lim_{p \to \infty} \pr (\widehat{\cS}_{\rm best} = \cS) \leq 1 - \frac{ e^{-1/2}}{ 2\sqrt{2 \pi}} + 2 e^{-8} < 0.9.
\]
In other words, for $a = \{2 (\log p)/n\}^{1/2}$ we have 
\[
\lim_{p \to \infty} \sup_{\beta \in \cM_s^a}\pr (\widehat{\cS}_{\rm best} \neq \cS)  > \frac{1}{10}.
\]

\subsection{Proof of Proposition 4.4} %4.3
We first present a result form \cite{wang2010} which is gives us necessary condition for asymptotic exact recovery.

\begin{theorem}[\cite{wang2010}]
Consider the model (1) with the design matrix $X\in \R^{n \times p}$ be drawn with i.i.d elements from any distribution with zero mean and unit variance. Let $a:= \min_{j \in \calS(\beta)}\abs{\beta_j}$, i.e.,  it denote the minimum signal strength of $\beta$. Define the function 
\[
f_m(p,s, a) := \frac{\log \binom{p-s+ m}{m} -1}{\frac{1}{2} \log\left( 1 + m a^2 (1 - \frac{m}{p-s+m})\right)},\quad 1\leq m \leq s.\] 
Then $n \geq  \max\{ f_1(p,s,a), \ldots, f_s(p,s,a),s\}$ is necessary for asymptotic exact recovery.
\label{thm: wang information bound}
\end{theorem}
In the light of Theorem \ref{thm: wang information bound} the proof of Proposition 4.4 follows immediately. To see this note that if $r<1$ then there exists $\alpha \in (0,1)$ such that $r= 1- \alpha$. Also recall that $a = \{2r (\log p)/n\}^{1/2}, s= O(\log p)$ and $n = \floor{p^k}$. Note that $f_1(p,s,a)/n \sim \frac{1}{1-\alpha}$. This shows that asymptotically the necessary condition in above theorem is violated and hence $r\geq 1$ is necessary.

\section{Results related to ETS}
\subsection{Proof of Theorem 5.2} %5.1
\label{sec: Proof of main ETS}
We first briefly describe the main steps of the proof:
\begin{enumerate}
    \item We first establish the $\ell_2$-error bound of $\hat{\beta}$, i.e., we show that $\norm{\hat{\beta} - \beta}
_2 \leq \epsilon^{1/2}$.
    \item Next, we upper bound the 0-1 loss $\pr_\beta (\hat{\eta} \neq \eta)$ by decomposing it across the coordinates.
    \item We analyze each of the terms separately and use $\ell_2$-error bound along with Gaussian tail inequalities to establish model consistency.
\end{enumerate}
Now we are ready for the main proof.
Due to Assumption 5.1, there exists a sequence $\{\alpha_p\}_{p\geq 1} \subseteq [0, \infty)$ converging to 0 such that with probability $1 - \alpha_p$ the following is true:
$\cA$ requires no more than $T(\epsilon, p, \beta)$ iterations to output $\hat{\beta}$ that satisfies  $\norm{\hat{\beta} - \beta}_2 \leq \epsilon^{1/2}$. 
%Also, recall under the working assumptions, $s \log p\over n_1$ tends to 0 as $p \to \infty$.

For notational brevity, we write $\eta$ instead of  $\eta_\beta$. Define the event $\his = \left\{ \Vert \hat{\beta} - \beta \Vert_2 \leq \epsilon^{1/2} \right\}$ and we have $\pr (\cH^c)\leq \alpha_p$. Note that $\hat{\beta}$ is based on the subsample $\D_1$. 
%From Theorem 3 of \cite{Jain2014iterative}, it follows that $\pr(\his^c) \leq \exp(-c_0 n_1) + \exp(- c_4 n_1) + c_2 p^{-c_3}$. Now, we reparameterize $\epsilon_0$ by $\epsilon/A$. 

Next, for algebraic convenience we again reparametrize $\delta$ as $8 \delta_0$ and set $\epsilon = 6 \delta_0, \varsigma = (1 + \epsilon)^{1/2}$. Now note that for any $\beta \in \M_{a}^s$, we have

\begin{align*}
 &\pr_\beta (\hat{\eta} \neq \eta \vert \D_1)\leq \sum_{j:\beta_j=0} \pr( \hat{\eta}_j =1, \his \vert \D_1) + \sum_{j: \beta_j \neq 0} \pr_\beta (\hat{\eta}_j \neq 1, \his\vert \D_1) + \pr(\his^c\vert \D_1)\\
& = \sum_{j: \beta_j =0 } \pr_\beta (\abs{\Delta_j} > \kappa_\varsigma(X_j^{(2)}), \his\vert \D_1) +   \sum_{j: \beta_j \neq 0 } \pr_\beta (\abs{\Delta_j} \leq \kappa_\varsigma(X_j^{(2)}), \his\vert \D_1) + \pr(\his^c \vert \D_1),
\end{align*}

%Here, $\kappa_\varsigma(X_j^{(2)})\geq 0$ as $s\leq p/2$. 
Using the fact that conditionally on $\hat{\beta}$ and $X_j$, the random variable $\Delta_j$ has the same distribution as the random variable in (11) in the main paper, we conclude that for all $j \notin \calS(\beta)$ ,
\begin{align*}
    \pr(\eta_j = 1, \his \vert \D_1)
    & \leq   \pr \left( (1 +  \epsilon)^{1/2}\abs{g_j}> \frac{a \norm{X_j^{(2)}}_2}{2} + \frac{(1+  \epsilon) \log p}{a \norm{X_j^{(2)}}_2} , \his \big \vert \D_1 \right)
    \\
    & = 2 \E \left \{ \barPhi\left ( \frac{a \norm{X_j^{(2)}}_2}{2 (1 + \epsilon)^{1/2}} + \frac{(1+  \epsilon)^{1/2} \log p}{a \norm{X_j^{(2)}}_2}\right) \right\}.
    %\\
   %& \lesssim \frac{1}{\sqrt{\log p}} p^{- \frac{ \left(\frac{r(1-\gamma)}{\sqrt{1 + A \epsilon}} +  \sqrt{1 + A \epsilon} \right)^2}{4r(1- \gamma)}},
\end{align*}
Here $\barPhi(\cdot)$ denotes the survival function of the standard Gaussian random variable. Now note that for each $j$ we have $\norm{X_j^{(2)}}_2^2 \overset{\rm d}{=} V_{n_2}$, where $V_{n_2}$ is a chi-squared random variable with $n_2$ degrees of freedom. Thus we have 
$$\pr (\eta_j = 1, \his \vert \D_1) \leq 2 \E \left \{ \barPhi\left ( \frac{a V_{n_2}^{1/2}}{2 (1 + \epsilon)^{1/2}} + \frac{(1+ \epsilon)^{1/2} \log p}{a V_{n_2}^{1/2}}\right) \right\}.$$

Analogous argument and the fact that $\abs{\beta_j}\geq a$ for all $\beta_j \neq 0$, leads to the fact that for all $j\in \calS(\beta)$,
\[
\pr(\eta_j \neq 1, \his \vert \D_1)\leq 2 \bbE \left \{ \barPhi\left ( \max \left \{\frac{a V_{n_2}^{1/2}}{2 (1 + \epsilon)^{1/2}} - \frac{(1+ \epsilon)^{1/2} \log p}{a V_{n_2}^{1/2}} , 0 \right \}\right) \right\}.
\]

Now recall that $\epsilon =  6\delta_0$ and $\gamma \in ( 0,\frac{\delta_0}{1 + 8\delta_0})$. With this choice of tuning parameters it is easy to see that $r  (1-\gamma ) / (1 + \epsilon) \geq \frac{1 + 7\delta_0}{ 1 + 6\delta_0}>1$ and hence as $p \to \infty$ we have

\begin{align*}
W_{n_2} & := \frac{1}{(\log p)^{1/2}}\left(\frac{a V_{n_2}^{1/2}}{2 (1 + \epsilon)^{1/2}} - \frac{(1+ \epsilon)^{1/2} \log p}{a V_{n_2}^{1/2}}\right)\\
&
\overset{\rm p}{\longrightarrow} \frac{1}{(2r)^{1/2}} \left\{ r  \left(\frac{1 - \gamma}{ 1 + \epsilon}\right) ^{1/2} - \left( \frac{1 + \epsilon}{ 1 - \gamma}\right)^{1/2} \right\}> 0.
\end{align*}
The above display uses the fact that $n_2/n \to 1 - \gamma  $ and $V_{n_2}/n_2 \overset{\rm p}{\to} 1$ as $p \to \infty$.
Next let is define the following quantity $q$: $$ q: = q( \epsilon, \delta_0, \gamma) =  \frac{1}{\{2 (1+ 8 \delta_0)\}^{1/2}} \left\{ (1+ 8 \delta_0)  \left(\frac{1 - \gamma}{ 1 + \epsilon}\right)^{1/2} - \left(\frac{1 + \epsilon}{ 1 - \gamma} \right)^{1/2} \right\}.$$
Due to choice of $\epsilon$ and $\gamma$ it is easy to show $q>0$.
 Now define the event $G_{n_2}: = \{ W_{n_2} > q/2 \}$. Before we proceed it is useful to note the following:
 \[
 W_{n_2} = \frac{1}{(2 r)^{1/2}} \left(  \frac{r \{V_{n_2}/(n(1-\gamma))\}^{1/2} (1-\gamma)^{1/2}}{(1 + \epsilon)^{1/2}} - \frac{(1 + \epsilon)^{1/2}}{\{V_{n_2}/(n(1-\gamma))\}^{1/2}(1-\gamma)^{1/2}}\right).
 \]
 Next define the function
 \[
 H (u)  := \frac{1}{( 2+ 16 \delta_0)^{1/2}} \left\{ u (1+ 8 \delta_0)  \left(\frac{1 - \gamma}{ 1 + \epsilon}\right)^{1/2} - \frac{1}{u}\left( \frac{1 + \epsilon}{ 1 - \gamma}\right)^{1/2} \right\}, \quad u >0.
 \]
 As $r = 1 + 8\delta_0$ we have $W_{n_2} = H( \{V_{n_2}/(n(1-\gamma))\}^{1/2} )$. It is also easy to see that $H(\cdot)$ is strictly increasing function on $(0,\infty)$ and $H(1) = q$. Hence $\lambda_{\delta_0}: = H^{-1}(q/2) \in (0,1)$. Now $G_{n_2}^c = \{W_{n_2} \leq q/2\}\subseteq \{ H(\{V_{n_2}/(n(1-\gamma))\}^{1/2}) \leq q/2 \}$. Thus a straightforward calculation shows that
 \[
 \pr (G_{n_2}^c )  \leq \pr \left(  \frac{V_{n_2}}{n_2} \leq \frac{n (1-\gamma)}{n_2}\lambda_{\delta_0}^2
 \right). 
 \]
Choose $p$ large enough such that $n_2/n > \lambda_{\delta_0} (1-\gamma)$ and hence we have,
\[
\pr (G_{n_2}^c)   \leq \pr \left( \frac{V_{n_2}}{n_2} \leq\lambda_{\delta_0}\right) \lesssim \exp(- K_{\delta_0} n_2),
\]
where $K_{\delta_0} = (1- \lambda_{\delta_0})^2/8$. Note that $\barPhi(t) \leq e^{-t^2/2}$ for all $t>0$. Using this fact we have the following:
\[
\pr (\eta_j = 1, \his \vert \D_1)\leq \E \left[\exp \left\{- \left(1+ \frac{W_{n_2}^2}{2} \right) \log p \right\} \ind_{G_{n_2}}\right] + \pr(G_{n_2}^c) \lesssim p^{-(1+ q^2/8)} + \exp( - K_{\delta_0} n_2), 
\]
for all $j \notin \calS(\beta)$.
Similarly,
\[
\pr (\eta_j \neq  1, \his \vert \D_1) \lesssim p^{- q^2/8} + \exp( - K_{\delta_0} n_2), \quad \forall j \in \calS(\beta).
\]

% As $r>1$ there exists $\mu>0$ and $0<\gamma <1$ such that $r\geq \frac{1+ \mu}{1-\gamma}$. Now set $\epsilon  = \frac{\mu}{2A}$ and $n_1 = \floor{\gamma n}$ which ensures $r(1-\gamma)>(1 + \epsilon)$. Hence by Mill's ration we have,

% \[
% \pr(\eta_j \neq 1, \his \vert \D_1) \lesssim \frac{1}{\sqrt{\log p}} p^{- \frac{ \left(\frac{r(1-\gamma)}{\sqrt{1 + \epsilon}} -  \sqrt{1 + \epsilon} \right)^2}{4r(1- \gamma)}}.
% \]

Thus marginalizing out $\D_1$ and summing over all $j$ we get,
\[
\sup_{\beta\in \M_{a}^s} \pr_\beta (\hat{\eta}\neq \eta) \lesssim  p^{-q^2/8} \log p  + p\exp(- K_{\delta_0} n_2) + \alpha_p.
\]
The result follows by taking $p\to \infty$. 

\begin{remark}
   Note that $q^2 = \Omega(\frac{\delta_0^2}{1 + \delta_0^2})$ and it shows that the upper bound in the above display deteriorates as $\delta_0 \to 0$. Also, as $\delta_0$ approaches 0, the term $K_{\delta_0}$ also approaches 0. Hence, the rate of decay worsens as $\delta_0 \to 0$, and ETS continues to lose statistical power. 
\end{remark}

\subsection{Proof of Corollary 5.3} % 5.2
\label{sec: Proof of claim ets_tilde}
Similar to previous proofs, we reparametrize $\delta$ by $8 \delta_0$  and set $\epsilon = 6 \delta_0, \varsigma = (1 + \epsilon)^{1/2}$. Now note that it is enough to prove the following:
\[
\lim_{p \to \infty}\inf_{\beta \in \M_s^a}\pr_\beta \left( \max_{j \notin \calS(\beta)} \abs{\Delta_j} < \min_{j \in \calS(\beta)} \abs{\Delta_j} \right)\to 1
\]
as $p \to \infty$. To this end first define the following quantity:
\[
t_p := \frac{\left(\frac{2r n_2 \log p}{n} \right)^{1/2}}{2} + \frac{\varsigma^2\log p}{\left(\frac{2r n_2 \log p}{n} \right)^{1/2}}.
\]
We will show that $\lim_{p \to \infty} \inf_{\beta \in \M_s^a}\pr_\beta \left(\min_{j \in \calS(\beta) } \abs{\Delta_j} > t_p, \max_{j \notin \calS(\beta) } \abs{\Delta_j} \leq t_p\right) \to 1$ as $p \to \infty$. For convenience let us define the events $G_{\min} : = \{\min_{j \in \calS(\beta) } \abs{\Delta_j} > t_p \}$ and $G_{\max} : = \{ \max_{j \notin \calS(\beta) } \abs{\Delta_j} \leq t_p \}$. Let $\his$ be the event as defined in Section \ref{sec: Proof of main ETS}. First we will analyze $\pr_\beta (G_{\min}^c)$. Note that $\pr_\beta (G_{\min}^c ) \leq \pr_\beta (G_{\min}^c \cap \his) + \pr_\beta(\his^c)$. Now the second term goes to $0$ uniformly over $\beta \in \M_{s}^a$. Also using Equation (11) under the event $\his$ we get
%\newpage
\begin{align*}
& \sup_{\beta \in \M_{s}^a}\pr_\beta (G_{\min}^c \cap \his)\\
&\leq \sup_{\beta \in \M_{s}^a} \pr_\beta \left( 
\min_{j \in \calS(\beta)}\abs{\beta_j  \Vert X_j^{(2)}\Vert_2 + (1 + \epsilon)^{1/2} g_j} \leq t_p
\right)\\
& \leq \sup_{\beta \in \M_{s}^a}\pr_\beta \left(
\max_{j \in \calS(\beta) } \frac{\abs{g_j}}{(\log p)^{1/2}} \geq \frac{1}{(1 + \epsilon)^{1/2} (\log p)^{1/2}} \left\{ a \min_{j \in \calS(\beta)} \Vert X_j^{(2)}\Vert_2 - t_p\right\}
\right)\\
& \leq  \sup_{\beta \in \M_{s}^a} \pr_\beta \left(
\max_{j \in \calS(\beta) } \frac{\abs{g_j}}{(\log p)^{1/2}} \geq \frac{1}{(1 + \epsilon)^{1/2} (\log p)^{1/2}} \left\{ a \min_{j \in [p]} \Vert X_j^{(2)}\Vert_2 - t_p\right\}
\right)
\end{align*}
where $\{g_j\}_{j \in \calS(\beta)}$ are non i.i.d. standard Gaussian. Note that $\abs{\calS(\beta)} = O(\log p)$. Hence 
\[\max_{j \in \calS(\beta)} \abs{g_j} = O_{\pr}(\log \log p),\]
which tells that $$\max_{j \in \calS(\beta) } \frac{\abs{g_j}}{(\log p)^{1/2}}\overset{\rm p}{\to} 0.$$
Also using lemma 3 from \cite{fletcher2009necessary} we have 
\[
\frac{1}{(1 + \epsilon)^{1/2} (\log p)^{1/2}}\left(a \min_{j \in [p]} \Vert X_j^{(2)}\Vert_2 - t_p \right) \overset{\rm p}{\to} \frac{1}{(2r)^{1/2}} \left\{ r  \left(\frac{1 - \gamma}{ 1 + \epsilon}\right) ^{1/2} - \left( \frac{1 +  \epsilon}{ 1 - \gamma}\right)^{1/2} \right\}.
\]
The right-hand side of the above display is at least $q( \epsilon, \delta_0, \gamma)$ (defined in Section \ref{sec: Proof of main ETS}) which is strictly positive. Again for compactness we use $q$ instead of $q(\epsilon, \delta_0, \gamma)$. The above display motivates us to define the following event:
\[
\calE_p = \left\{ \frac{1}{(1 + \epsilon)^{1/2} (\log p)^{1/2}} \left( a \min_{j \in [p]} \Vert X_j^{(2)}\Vert_2 - t_p\right) \geq q/2\right\},
\]
and it follows that $\pr(\calE_p^c)\to 0 $ as $p \to \infty$. This leads to the following inequality:
%\newpage
\begin{align*}
\sup_{\beta \in \M_{s}^a}\pr_\beta (G_{\min}^c \cap \his) \leq & \sup_{\beta \in \M_{s}^a} \pr_\beta \left(
\max_{j \in \calS(\beta) } \frac{\abs{g_j}}{(\log p)^{1/2}} \geq q/2\right) + \pr(\calE_p^c)\\
& \lesssim p^{-q^2/8} \log p + \pr(\calE_p^c) \to 0.
\end{align*}

Thus we have $\sup_{\beta \in \M_s^a}\pr_\beta (G_{\min}^c)\to 0$. Similarly it can be shown that $\sup_{\beta \in \M_s^a} \pr_\beta (G_{\max}^c)\to 0$ as $p \to \infty$. These two claims together complete the proof.

\subsection{Discussion on Remark 5.4} %5.3
As $r>1 +  \delta_*$, by reparameterizing $\delta_*$ by $8 \Tilde{\delta}$, we have $r > 1 + 8 \Tilde{\delta}$. Now we are basically back to the setting of the proof of Theorem 5.1 and all of the proof steps are exactly the same as that of Theorem 5.1 with $\tilde{\delta}$ in place of $\delta_0$. This allows us to choose the threshold using the knowledge of $\delta_*$, and we do not need the knowledge of $a$. In particular, one can construct the threshold $\kappa_\varsigma(X_i^{(2)})$ with $\varsigma = (1+ A_2 \delta_*)^{1/2}$, where $A_2$ is the same universal constant as described in Theorem 5.1

\subsection{Discussion on examples of ETS}
%\label{sec: model consistency of CCD/FISTA}
\subsubsection*{\underline{\textbf{Solving $\ell_0$-constrained problem}:}}
\subsubsection*{Proofs for ETS-IHT:}
We first introduce some standard assumptions for analyzing ETS-IHT.
\begin{definition}[RSC property]
A differentiable function $F: \R^p \to \R$ is said to satisfy restricted strong convexity (RSC) at sparsity level $s= s_1+s_2$ with strong convexity constraint $\ell_s$ if the following holds for all $\theta_1,\theta_2$ s.t. $\norm{\theta_1}_0\leq s_1$ and $\norm{\theta_2}_0\leq s_2$:
\[
F(\theta_1) - F(\theta_2)\geq \innerprod{\theta_1- \theta_2, \nabla_\theta F(\theta_2)} + \frac{\ell_s}{2}\norm{\theta_1- \theta_2}_2^2.
\]
\label{def: RSC_property}
\end{definition}

\begin{definition}[RSS property]
A differentiable function $F: \R^p \to \R$ is said to satisfy restricted strong smoothness (RSS) at sparsity level $s= s_1+s_2$ with strong smoothness constraint $L_s$ if the following holds for all $\theta_1,\theta_2$ s.t. $\norm{\theta_1}_0\leq s_1$ and $\norm{\theta_2}_0\leq s_2$:
\[
F(\theta_1) - F(\theta_2)\leq \innerprod{\theta_1- \theta_2, \nabla_\theta F(\theta_2)} + \frac{L_s}{2}\norm{\theta_1- \theta_2}_2^2.
\]
\label{def: RSS_property}
\end{definition}

Now we quote an important theorem from \cite{Jain2014iterative} which quantifies the sub-optimality gap of Algorithm 2.
\begin{theorem}[\cite{Jain2014iterative}]
Let $F$ has RSC and RSS parameters given by $\ell_{2 \hat{s} + s}(F)= \alpha$ and $L_{2 \hat{s} + \hat{\pi}}(F)= L$ respectively. Call Algorithm 2 with $\hat{s} \geq 32 L^2 \ell^{-2} s$ and $h = 2/(3L)$. Also let $\hat{\beta}= \argmin _{\theta, \norm{\theta}_0\leq s} F(\theta)$. Then $t{\rm th}$ iterate of Algorithm 2 for $t = O(L \ell^{-1} \log (F(\beta^{(0)})/\epsilon))$ satisfies:
\[
    F(\beta^{(t)}) - F(\hat{\beta}) \leq \epsilon.
\]
\label{thm: two_stage_IHT_sub_optimality}
\end{theorem}

In our setup, the observations $\{x_i\}_{i=1}^n$ are coming from i.i.d. mean zero isotropic Gaussian distribution. Thus, lemma 6 from \cite{agarwal2012fast} immediately tells that RSC and RSS at any sparsity level $m$ hold for $f_{n_1}(\cdot)$ with probability at least $1- \exp(- c_0 n_1)$ with $\ell_m = \frac{1}{2} - c_1 (m \log p)/n_1$ and $L_m = 2 + c_1 (m \log p)/n_1$, where $c_0,c_1$ are universal constants. Now set $m = 2 \hat{s} + s$ and recall that $n_1  \sim \gamma p^k$. If $n_1 > 4c_1 (2 \hat{s} + s) \log p$ then we have $\ell_{m} \geq 1/4$ and $L_{m}\leq 9/4$, which means that $L_m/(9\ell_m)\leq 1$. Thus to apply Theorem \ref{thm: two_stage_IHT_sub_optimality} it is enough to choose $\hat{s} = 2592 s$.
%Thus to apply Theorem \ref{thm: two_stage_IHT_sub_optimality} it is enough to ensure $\hat{\pi}\geq s$ and $\hat{s}\geq 323 \hat{\pi} + s$. Hence set $\hat{\pi} = 2s$ and $\hat{s} = 647 s$. 
Also by the assumption on $n$ for large $p$ we have $n_1> 4c_1 (2 \hat{s}+ s) \log p$. Let $f_{n_1}(\theta): = n_1^{-1} \Vert Y^{(1)} -X^{(1)}\theta \Vert_2^2$ for $\theta \in \R^{p}$. Note that $f_{n_1}(0) = n_1^{-1}\Vert Y^{(1)}\Vert_2^2 \overset{\rm d}{=} (1 + \norm{\beta}_2^2) V_{n_1}/n_1$, where $V_{n_1}$ is chi-square random variable with $n_1$ degrees of freedom. Also by Bernstein's type inequality it follows that $\abs{(V_{n_1}/n_1) -1}\leq 1/2$ with probability at least $1 - \exp(-c_4 n_1)$, where $c_4$ is a universal positive constant. Thus if $t = O(L_m \ell_m^{-1} \log ((1 + \norm{\beta}_2^2)/\epsilon_0)) = O( \log p + \log((1 +  \norm{\beta}_\infty)/\epsilon_0))$, then we have $f_{n_1}(\beta^{(t)}) - f_{n_1}(\hat{\beta})\leq \epsilon_0$. Thus by Theorem 3 of \cite{Jain2014iterative} it follows that with probability at lest $1 - \exp(-c_0 n_1) - \exp(- c_4 n_1) - c_2 p^{-c_3}$ ($c_2, c_3$ are universal constants) we have 
\[
\Vert\beta^{(t)} - \beta\Vert_2 \leq C \left(\frac{s \log p}{n_1} \right)^{1/2} + (8\epsilon_0)^{1/2}  \leq (9\epsilon_0)^{1/2}, \quad \text{for large $p$}.
\]
 $C$ is a positive universal constant in the above inequality. If we set $\epsilon = 9 \epsilon_0$, then Assumption 5.1 holds with $\alpha_p$ equal to   $\exp(-c_0 n_1) + \exp(- c_4 n_1) + c_2 p^{-c_3}$ and $T(\epsilon, p, \beta) = O( \log p + \log((1 +  \norm{\beta}_\infty)/\epsilon)))$.
% In practice, one can use the full data in both steps of ETS-IHT but for theoretical convenience, we use data splitting in the algorithm. The key idea here is to use one subsample for obtaining $\hat{\beta}^{\rm iht}$ and the remaining subsample in the screening step. Due to this $\hat{\beta}^{\rm iht}$ becomes independent of the second subsample, thus making the analysis of ETS-IHT easy.  

\subsubsection*{\underline{\textbf{Solving $\ell_1$-regularized problem}:}}
We start by recalling the definition of the loss function
$
f_{n_1, \lambda}(\theta) = n_1^{-1} \Vert Y^{(1)} - X^{(1)} \theta\Vert_2^2 + \lambda \norm{\theta}_1,
$
and define the minimizer of the loss function $\hat{\beta}_L:= \argmin_\theta f_{n_1, \lambda}(\theta)$. First, we will prove Proposition 5.5.

\subsubsection*{Proof of Proposition 5.5:}
Let us assume
\[
f_{n_1,\lambda}(\hat{\beta}) - f_{n_1, \lambda}(\hat{\beta}_L) \leq \epsilon_0 <1.
\]

Now, we will establish the $\ell_2$-error rate between $\hat{\beta}$ and $\beta$. We write $\widehat{b} = \hat{\beta} - \beta$. Since, $\hat{\beta}_L$ is optimal, we have 
\[
f_{n_1, \lambda}(\hat{\beta}) \leq f_{n_1,\lambda}(\hat{\beta}_L) + \epsilon_0 \leq f_{n_1, \lambda}(\beta) + \epsilon_0.
\]
Rearrangement of the above inequality yields

\begin{equation}
    \label{eq: rearranged lagrangian}
    0 \leq \frac{1}{n_1}\norm{X^{(1)} \hat{b}}_2^2 \leq \frac{2 E^\top X^{(1)} \hat{b}}{n_1} + \lambda (\norm{\beta}_1 - \Vert\hat{\beta}\Vert_1) + \epsilon_0.
\end{equation}
Since, $\Norm{X_j^{(1)}}_2^2 \sim \chi^2_{n_1}$, we have 
\[
 \pr \left(\underbrace{\max_{j \in [p]}n_1^{-1/2}\norm{X_j^{(1)}}_2 < \sqrt{4/3} }_{:=\cE} \right) \geq 1- 2p^{-2}, \quad \text{for large $p$.}
 \]
 Using Gaussianity of $X_j^{(1)\top}E/ \Norm{X_j^{(1)}}_2$, we have
 \begin{align*}
    \pr \left( \frac{1}{n_1} \norm{X^{(1)\top} E}_\infty > 2 \sqrt{\frac{\log p}{n_1}}\right) 
    &\le \pr \left( \max_{j \in [p]} \abs{\frac{X_j^{(1)\top} E}{\norm{X_j^{(1)}}_2}} >  \sqrt{3 \log p}\right) + \pr(\cE^c)\\
    & \leq 2 p^{-0.5} + 2 p^{-2}.
\end{align*}
 Setting $\lambda = 8 \{(\log p)/n_1\}^{1/2}$, we have $\lambda \geq 2 \norm{X^{(1)\top}E}_\infty/n_1$ with probability at least $1 - 2 p^{-0.5} - 2 p^{-2}$. 
Now, since $\beta$ is $s$-sparse with support on $\cS$, we have
\[
\norm{\beta}_1 - \Vert \hat{\beta}\Vert_1  = \norm{\beta_\cS}_1 - \norm{\beta_\cS + \hat{b}_\cS}_1 - \norm{\hat{b}_{\cS^c}} \leq \norm{\hat{b}_\cS}_1 - \norm{\hat{b}_{\cS^c}}_1.
\]
Substituting this in the basic inequality \eqref{eq: rearranged lagrangian} and using $\lambda \geq 2 \norm{X^{(1)\top}E}_\infty/n_1$, we get 
\begin{equation}
    \label{eq: rearranged lagrangian 2}
    \begin{aligned}
    0 \leq \frac{1}{n_1}\norm{X^{(1)}\hat{b}}_2^2  & \leq 2 \norm{\frac{ E^\top X^{(1)} }{n_1}}_\infty \norm{\hat{b}}_1 + \lambda(  \norm{ \hat{b}_\cS}_1 - \norm{\hat{b}_{\cS^c}}_1) + \epsilon_0\\
    & \leq  (\lambda/2) \norm{\hat{b}}_1 + \lambda(  \norm{ \hat{b}_\cS}_1 - \norm{\hat{b}_{\cS^c}}_1) + \epsilon_0\\
    & \leq (\lambda/2)\{3 \norm{\hat{b}_\cS}_1 - \norm{\hat{b}_{\cS^c}}_1\} + \epsilon_0.
    \end{aligned}
\end{equation}

Hence, we have 
\begin{align*}
\norm{\hat{b}}_1^2 &= (\norm{\hat{b}_\cS}_1 + \norm{\hat{b}_{\cS^c}}_1)^2 \leq (4 \norm{\hat{b}_\cS}_1 + (2 \epsilon_0/\lambda))^2  \leq 32 \norm{\hat{b}_\cS}_1^2 + \frac{8 \epsilon_0^2}{\lambda^2}.
\end{align*}
Now by Theorem 7.16 of \cite{wainwright2019high}, we have the following with probability at least $1 - 2 \exp(-n_1/32)$:
\[
\frac{\norm{X^{(1)} \theta}_2^2}{n_1} \geq c_1 \norm{\theta}_2^2 - c_2 \frac{\log p}{n_1} \norm{\theta}_1^2 \quad \text{for every $\theta \in \bbR^p$},
\]
where $c_1, c_2>0$ are universal constants.
Using the above fact we have
\begin{equation*}
    \begin{aligned}
    \frac{\norm{X^{(1)}\hat{b}}_2^2}{n_1} &\geq c_1 \norm{\hat{b}}_2^2 - 32c_2 \frac{s\log p}{n_1} \norm{\hat{b}}_2^2 - 8c_2 \frac{\log p}{n_1 \lambda^2}\epsilon_0^2  \geq \frac{c_1}{2} \norm{\hat{b}}_2^2 - \epsilon_0^2,
    \end{aligned}
\end{equation*}
when $32c_2 s \log p/(n_1) < c_1/2$ and $8c_2 \log p/(n_1\lambda^2)<1$. This is possible for large enough values of $p$ and $\lambda$.

\textbf{Case 1:} If $(c_1/4) \norm{\hat{b}}_2^2 > \epsilon_0^2$, then using \eqref{eq: rearranged lagrangian 2}, we get 
\[
\frac{c_1}{4} \norm{\hat{b}}_2^2 \leq \frac{3 \lambda \sqrt{s}}{2} \norm{\hat{b}}_2 + \epsilon_0.
\]
This bound involves a quadratic form of $\norm{\hat{b}}_2$; computing the roots of the quadratic form we get the following bound:
\[
\norm{b}_2 \leq \underbrace{\frac{6 \lambda \sqrt{s}}{c_1}}_{= O(\sqrt{(s \log p)/n_1})} + \quad \frac{2 \sqrt{\epsilon_0}}{\sqrt{c_1}}.
\]

\textbf{Case 2:} If $(c_1/4) \norm{\hat{b}}_2^2 \leq \epsilon_0^2$, then $\norm{\hat{b}}_2 \leq 2 \epsilon_0/\sqrt{c_1} < 2 \sqrt{\epsilon_0/c_1}$. The last inequality uses the fact that $\epsilon_0<1$.

Combining the bounds obtained in Case 1 and Case 2 and with probability at least $1 - 2 p^{-0.5} - 2 p^{-2} - e^{-n_1/32}(\text{which is $\ge 1 - 3p^{-0.5}$ for large $p$})$, we finally have $\Vert\hat{\beta} - \beta\Vert_2 \leq (C_3\epsilon_0)^{1/2}$ for large enough $p$ and $C_3$ being an absolute constant. Now, using the reparameterization $\epsilon = C_3 \epsilon_0$, we have $\epsilon <C_3$, and the initial sub-optimality gap turns out to be $\epsilon/C_3$. This finishes the proof.

\subsubsection*{Proof for ETS-PICASSO:}

First, we will show that the estimator generated by the PICASSO algorithm is a good estimate of $\beta$. 
%Specifically, we will use Lemma E.1 in the supplementary material of \cite{Zhao2018Pathwise} (DOI: \href{https://projecteuclid.org/journals/annals-of-statistics/volume-46/issue-1/Pathwise-coordinate-optimization-for-sparse-learning-Algorithm-and-theory/10.1214/17-AOS1547.full?tab=ArticleLinkSupplemental}{10.1214/17-AOS1547SUPP}). 
Let $\betapic$ be the estimator obtained by applying PICASSO for minimizing $f_{n_1}(\theta; \lambda, \cD_1)$. Now define the largest and smallest $s_0$-sparse eigenvalues of $G:=X^{(1)\top} X^{(1)}/n_1$ as:
\[
\rho_+(s_0) := \max_{v: \norm{v}_0\le s} \frac{v^\top G v}{\norm{v}_2^2}; \quad \text{and} \quad \rho_-(s_0) := \min_{v: \norm{v}_0\le s} \frac{v^\top G v}{\norm{v}_2^2}.
\]

Let $\tilde{s}_\psi = (484 \psi^2 + 100 \psi) s$, where $\psi>0$ is constant. 
An application of union bound and Equation 4.22 in \cite{vershynin_2018} yields that
\begin{equation}
    \pr \left\{
    \max_{\cD: \abs{\cD} \leq s+2\tilde{s}_2}\norm{\frac{X^{(1)\top}_\cD X^{(1)}_\cD}{n_1} - \sfI_p}_\op
    \lesssim \frac{\log p}{p^{k/2}}\right\} \geq 1 - 2 p^{-2},
    \label{eq: spectral concentration}
\end{equation}
 for large values of $p$. This ensures that for large values of $p$, we have 
 \begin{equation}
 0.99 \leq \rho_-(s + 2\tilde{s}_2) \leq \rho_+(s + 2\tilde{s}_2) \leq 1.1.
 \label{eq: SRC inequality}
 \end{equation}
 Defining $\kappa:= \rho_+(s + 2\tilde{s}_2)/\rho_-(s + 2\tilde{s}_2)$, we have $\kappa <2$. Thus, Assumption 3.5 of \cite{Zhao2018Pathwise} holds with $\tilde{s} = \tilde{s}_2 > \tilde{s}_\kappa$. Also, \eqref{eq: spectral concentration} shows that for large $p$ 
 \[
 \pr \left(\underbrace{\max_{j \in [p]}n_1^{-1/2}\norm{X_j^{(1)}}_2 < \sqrt{4/3} }_{:=\cE} \right) \geq 1- 2p^{-2}.
 \]
% \newpage
Hence, we have 
\begin{align*}
    \pr \left( \frac{1}{n_1} \norm{X^{(1)\top} E}_\infty > 2 \sqrt{\frac{\log p}{n_1}}\right) 
    &\le \pr \left( \max_{j \in [p]} \abs{\frac{X_j^{(1)\top} E}{\norm{X_j^{(1)}}_2}} >  \sqrt{3 \log p}\right) + \pr(\cE^c)\\
    & \leq 2 p^{-0.5} + 2 p^{-2}.
\end{align*}
Hence, Assumption 3.1 of \cite{Zhao2018Pathwise} holds with high probability when $\lambda_N = \lambda \ge 8 \{\log p/n_1\}^{1/2}$, where $N$ denotes the final iteration count of the outermost loop of PICASSO and $\lambda_K$ denotes the regularization parameter at the $K$th iteration of the outer loop. Also, in this case, $N = O(\log (\Norm{X^{(1)^\top }Y^{(1)}/n_1}_\infty \sqrt{n_1/\log p}))$ which follows from the description of PICASSO (see Algorithm 3 in \cite{Zhao2018Pathwise}). From triangle inequality, it follows that 
\[
\norm{\frac{X^{(1)\top} Y^{(1)}}{n_1} }_\infty
 \leq \norm{G \beta}_\infty + \norm{X^{(1)\top} E/n_1}_\infty,
\]
and $\norm{G \beta}_\infty \leq \norm{G}_{\infty , \infty} \norm{\beta}_\infty \leq \sqrt{p} \norm{G}_\op \norm{\beta}_\infty$. Note that 
\[
\norm{G}_{\op} = \frac{p}{n_1} \norm{\frac{X^{(1)}X^{(1)\top}}{p}}_\op.
\]
Thus, applying Equation (4.22) in \cite{vershynin_2018}, we get $\norm{G\beta}_\infty \lesssim (p^{1.5}/n) \norm{\beta}_\infty$ with probability at least $2 \exp(-n_1)$. Hence, we have $N = O(\log p)$ with probability at least $1 - 2p^{-0.5} - 2p^{-2} - \exp(n_1)$.

Now will make sure that Assumption 3.7 of \cite{Zhao2018Pathwise} also holds. Before, going any further let us clarify some notations. At $K$th outer iteration of PICASSO, \cite{Zhao2018Pathwise} denotes the inner and middle loop precision parameters as $\tau_K$ and $\delta_K$, and the active set initialization parameter is $\varphi$. To be consistent with our notation we set $\epsilon_0 = \delta_K$ for every $K$. Also, we choose the parameters in such a way so that
\begin{equation}
\label{eq: picasso tuning params}
\epsilon_0 \leq \min\{1/8, C_3\}, \quad \tau_K \le \frac{\epsilon_0}{\rho_+(s + 2 \tilde{s})} \sqrt{\frac{\rho_-(1)}{\rho_+(1) (s + 2 \tilde{s})}} , \quad \varphi \leq 1/8,
\end{equation}
$C_3$ is the same constant ad in Proposition 5.5.
Using \eqref{eq: SRC inequality}, one can set $\tau_K = O(\epsilon_0/\sqrt{\log p})$ which will be less than 1 for large $p$. So, under the above conditions, Assumption 3.7 in \cite{Zhao2018Pathwise} holds. 
Hence, part (iii) of Theorem 3.12 in \cite{Zhao2018Pathwise} tells that 
\[
f_{n_1}(\betapic; \lambda, \cD_1) - f_{n_1}(\hat{\beta}_L; \lambda , \cD_1) \leq \epsilon_0 \frac{500 \lambda^2 s}{11}.
\]
If $\lambda = 8 \sqrt{(\log p)/n_1}$, then for large $p$ 
\[
f_{n_1}(\betapic; \lambda, \cD_1) - f_{n_1}(\hat{\beta}_L; \lambda , \cD_1) \leq \epsilon_0 /C_3,
\]
 Hence, due to Proposition 5.5, we have $\Norm{\betapic - \beta}_2^2 \leq \epsilon_0$. Also, using Theorem 3.12  and Lemma 3.13 of \cite{Zhao2018Pathwise}, we get that PICASSO needs no more than 
\[
T(\epsilon_0, \beta, p) = O\bigg((\log p + \log \norm{\beta}_\infty) (\log p)^3 \{\log \log p +  \log (\epsilon_0^{-1}) \}\bigg).
\]
 But the above facts are true when $\epsilon_0 \leq \min\{1/8, C_3\}$. In order to extend the above results to a bigger range of $\epsilon_0$ we first define $C_\epsilon:= \min\{1/8, C_3\}/2$. If $\epsilon \leq C_\epsilon$, then one can get the same result by setting $\delta_K$ and $\tau_K$ appropriately as prescribed in \eqref{eq: picasso tuning params}. If $\epsilon> C_\epsilon$, then setting $\epsilon_0 = C_\epsilon$ and using \eqref{eq: picasso tuning params}, we again get  $\Norm{\betapic - \beta}_2^2 \leq \epsilon_0 < \epsilon$ within $O((\log p + \log \norm{\beta}_\infty) (\log p)^3 \{\log \log p +  \log (C_\epsilon^{-1}))$ iterations. Thus, conditions in Assumption 5.1 is met with \[
 T(\epsilon, p, \beta) = O\bigg((\log p + \log \norm{\beta}_\infty) (\log p)^3 \{\log \log p +  \log (\epsilon^{-1} \vee C_\epsilon^{-1})\bigg),
 \]
 and with probability at least $1 - O(p^{-0.5})$.

%Hence, due to Lemma E.1 in the supplementary material of \cite{Zhao2018Pathwise}, we have 
% $
% \Vert\betapic - \beta\Vert_2 \leq (A_4 \epsilon)^{1/2},
% $ for large $p$.  

% Recall that $\delta = r-1$.
% Now, again we can reparameterize $\delta $ as $8 \delta_0$ and set $\epsilon = \min \{6 \delta_0/ A_4, 1/64\}$. Using Theorem 3.12  and Lemma 3.13 of \cite{Zhao2018Pathwise}, we get that PICASSO needs no more than $A_3 (\log p + \log \norm{\beta}_\infty) (\log p)^3 \{\log \log p +  \log (\delta^{-1} \vee A_3) \}^2$ for some universal positive constant. Here, we used the fact that $\norm{G}_\infty \leq \sqrt{p}\norm{G}_\op = O_{\pr}(p^{1.5}/n)$.
% The rest of the proof follows exactly the same steps as in Section \ref{sec: Proof of main ETS} and hence omitted.

\subsubsection*{Proof for ETS-PGH:}
Let $\betapgh$ be the solution obtained by minimizing $f_{n_1}(\theta; \lambda, \cD_1)$ via PGH method. For the proof of this part we will use Theorem 3.2 of \cite{xiao2013proximal}. To apply that theorem we need to make sure that Assumption 3.2 of \cite{xiao2013proximal} is satisfied. To avoid notational confusion, we use $\tilde{\gamma}$ and $\tilde{\delta}$ to denote the parameters $\gamma, \delta^\prime$ considered in \cite{xiao2013proximal} respectively. Let $\tilde{\delta} = 0.1$ and $\Tilde{\gamma} = 2$. By a similar argument as before, it can be shown that with probability at least $1-2p^2$
\[
\kappa(G,s_0):= \frac{\rho_+(s_0)}{\rho_-(s_0)} \leq \frac{1+ \nu}{1-\nu},
\]
where $s_0 = \floor{46(1 + \tilde{\gamma})s}, \nu = 0.1$ and $p$ is sufficiently large. Also, we choose 
\begin{equation}
\begin{aligned}
\lambda & = 8 \max\{2, \frac{\Tilde{\gamma}+1}{\Tilde{\gamma}(1- \Tilde{\delta}) - (1 + \Tilde{\delta})}\} \{\log p /n_1\}^{1/2}\\
& \geq 4 \max\{2, \frac{\Tilde{\gamma}+1}{\Tilde{\gamma}(1- \Tilde{\delta} - (1 + \Tilde{\delta}))}\}\norm{X^{(1)^\top}E/n_1}_\infty.
\end{aligned}
\label{eq: lambda_pgh range}
\end{equation}
The last inequality of the above display shows that $\lambda \geq 8 \Vert E^\top X^{(1)}/n_1\Vert_\infty$ and it happens with the probability at least $1 - O( p^{-0.5})$.
Hence, following the arguments of the second bullet point on page 10 of \cite{xiao2013proximal}, we can conclude that the Assumption 3.2 of \cite{xiao2013proximal} holds with $\tilde{s} = \floor{22(1+\tilde{\gamma})s}$, $\gamma_{\rm inc} = 1.2$ (see \cite{xiao2013proximal}), $L_{\min}  = 1.32$ and $\lambda $. Using part 3 of Theorem 3.2 in \cite{xiao2013proximal}, for a given precision level $\epsilon_0$, we have 
\[
f_{n_1}(\betapgh; \lambda, \cD_1) - f_{n_1}(\hat{\beta}_L;  \lambda, \cD_1) \leq O(\epsilon_0 s \sqrt{(\log p)/n_1}) < \epsilon_0/C_3 , \quad \text{for large $p$.}
\]
$C_3$ is the same universal constant as in Proposition 5.5. If $\epsilon_0 = \min\{1, C_3\}/2=: C_4$,
then
$\Norm{\betapgh - \beta}_2^2 \leq \epsilon_0$, and
the total iteration complexity is 
\[
O \left(
\log p \log \log p + \log \max\{1, (\lambda^2/ \epsilon_0^2) \log p\}
\right),
\]
which for large value of of $p$, the the form $$O((\log p + \log \norm{\beta}_\infty)\log \log p + \log \max\{1, (1/\epsilon_0^2)\log p\}) \quad \text{(as $\lambda<1$).}$$
If $\epsilon_0 \le C_4$, then the order becomes  $$O((\log p + \log \norm{\beta}_\infty) \log\log p + \log (1/\epsilon_0)).$$ Otherwise, i.e., if $\epsilon_0>C_4
$ one can set the tolerance level at $\epsilon = C_4$ and the overall order in that case is $O((\log p + \log \norm{\beta}_\infty) \log \log p + \log (1/C_4))$. Thus, for a given tolerance level $\epsilon$, the total iteration complexity is $O((\log p + \log \norm{\beta}_\infty) \log\log p + \log (\epsilon^{-1}\vee C_4^{-1}))$. Thus, Assumption 5.1 holds with probability at least $1 - O(p^{-0.5})$ and $T(\epsilon, p, \beta ) = O((\log p + \log \norm{\beta}_\infty) \log\log p + \log (\epsilon^{-1}\vee C_4^{-1}))$.
\end{appendices}